\documentclass{amsart}
\usepackage{amsfonts,amssymb,amsmath,amsthm}
\usepackage{url}
\usepackage{enumerate}
\usepackage{xcolor}

\newtheorem{thm}{Theorem}[section]
\newtheorem{cor}[thm]{Corollary}
\newtheorem{lem}[thm]{Lemma}

\theoremstyle{definition}
\newtheorem{ex}[thm]{Example}

\newcommand{\N}{\mathbb{N}}

\newcommand{\R}{\mathbb{R}}

\newcommand{\M}{M}

\newcommand{\D}{\operatorname{D}}

\newcommand{\eps}{\varepsilon}

\newcommand{\rank}{\operatorname{rank}}
\newcommand{\sgn}{\operatorname{sgn}}
\newcommand{\pszm}[1]{\frac{\partial}{\partial{#1}}}

\begin{document}


\title[Affine hypersurfaces of arbitrary signature...]{Affine hypersurfaces of arbitrary signature with an almost symplectic form}

\author[M. Szancer]{Michal Szancer}

\author[Z. Szancer]{Zuzanna Szancer}

\begin{abstract}
 In this paper we study affine hypersurfaces with non-degenerate second fundamental form of arbitrary signature additionally equipped with an almost symplectic structure $\omega$.
 We prove that if $R^p\omega=0$ or $\nabla^p\omega =0$ for some positive integer $p$ then the rank of the shape operator is at most one.
The results provide complete classification of affine hypersurfaces with higher order parallel almost symplectic forms and are generalization of recently obtained results for
Lorentzian affine hypersurfaces.
\end{abstract}

\maketitle

\section{Introduction}
\label{intro}
\par Parallel structures are of the great interest in the classical Riemannian geometry (see \cite{D, Lumiste, BDI})
as well as in affine differential geometry (\cite{BNS, DV2, DVY, HLLVran, HLVran, Hildebrand2015, HLLV2}).
Higher order parallel structures are the natural generalization of parallel structures and are widely studied as well
(\cite{D, D2, Vrancken, LucVrancken2, VanderVeken}).

\par In \cite{BC} O. Baues and V. Cort\'{e}s  studied affine hypersurfaces equipped with an almost complex structure.
They showed that there is direct relation between simply connected special K\"{a}hler manifolds (\cite{Freed}) and improper affine hyperspheres.
Later  V. Cort\'{e}s together with M.-A. Lawn and L. Sch\"{a}fer  proved a similar result for special para-K\"{a}hler manifolds (\cite{CLS}).
In both cases an important role was played by the K\"{a}hlerian (resp. para-K\"{a}hlerian) symplectic form $\omega$.
The concept of special affine hyperspheres was generalized by the first author in \cite{MSz3}. Some other results related to affine hypersurfaces with
almost complex structures can be also found in the paper of M. Kon (\cite{Kon}). In all the above cases the important role was played
by an (almost) symplectic structure related in some way to the induced affine structure on a hypersurface. In particular, relation between an almost symplectic
structure $\omega$ and the induced affine connection $\nabla$ and its curvature $R$ seemed crucial.
\par The above results motivated the first author to study non-degenerate affine hypersurfaces
$f\colon M\rightarrow \R^{2n+1}$ with a transversal vector field $\xi$ additionally equipped with an almost symplectic structure $\omega$ in a more general setting.
More precisely affine hypersurfaces with the property $\nabla^p\omega=0$ or even more general $R^p\omega=0$.
In \cite{MSz} it was shown that if $\dim M\geq 4$ condition $R\omega=0$ implies that $\nabla$ must be flat (what is generalization of result obtained in \cite{Kon}) and the result
generalizes to an arbitrary power of $R$ under additional assumption that the second fundamental form is positive definite and the transversal vector field $\xi$ is locally equiaffine (i.e. $d\tau=0$).
Namely, we have
\begin{thm}[\cite{MSz}]\label{tw::RKOmega0}
Let $f\colon \M\rightarrow\R^{2n+1}$ be a non-degenerate affine hypersurface ($\dim M\geq 4$) with a locally equiaffine transversal vector field  $\xi$
and an almost symplectic form $\omega$. Additionally assume that the second fundamental form is positive definite on $M$.
If $R^p\omega=0$ for some positive integer $p$ then $\nabla$ is flat.
\end{thm}
From the above theorem it follows that
\begin{thm}[\cite{MSz}]\label{tw::NablaKOmega0}
Let $f\colon \M\rightarrow\R^{2n+1}$ be a non-degenerate affine hypersurface ($\dim M\geq 4$) with a locally equiaffine transversal vector field  $\xi$
and an almost symplectic form $\omega$. Additionally assume that the second fundamental form is positive definite on $M$.
If $\nabla^p\omega=0$ for some positive integer $p$ then $\nabla$ is flat.
\end{thm}
In the very same paper it was shown that the condition that the second fundamental form is positive definite cannot be relaxed since there exists affine hypersurface
with property $R^2\omega=0$ which is not flat.
\par Later in \cite{MSz2} it was shown that for  affine hypersurfaces with Lorentzian second fundamental form we still have
strong constrains on the shape operator $S$ if only $\dim M\geq 6$. More precisely, although it cannot be shown that $S=0$ (what is equivalent with $\nabla$ being flat)
one may still show that $\rank S\leq 1$. Recently (\cite{MSzZSz}) it was shown that the same constrains apply when $\dim M = 4$.
Combining results from \cite{MSz2} and \cite{MSzZSz} we have the following theorems:
\begin{thm}[\cite{MSz2}, \cite{MSzZSz}]\label{tw::RKomega0Lorentzian}
Let $f\colon \M\rightarrow\R^{2n+1}$ ($\dim M\geq 4$) be a non-degenerate affine hypersurface with a locally equiaffine transversal vector field  $\xi$
and an almost symplectic form $\omega$. If $R^p\omega=0$ for some $p\geq 1$ and the second fundamental form is Lorentzian on $M$ (that is has signature $(2n-1,1)$)
then the shape operator $S$ has the rank $\leq 1$.
\end{thm}

\begin{thm}[\cite{MSz2},\cite{MSzZSz}]\label{tw::NablaKOmega0Lorentzian}
Let $f\colon \M\rightarrow\R^{2n+1}$ ($\dim M\geq 4$) be a non-degenerate affine hypersurface with a locally equiaffine transversal vector field  $\xi$
and an almost symplectic form $\omega$. If $\nabla^p\omega=0$ for some $p\geq 1$ and the second fundamental form is Lorentzian on $M$ (that is has signature $(2n-1,1)$)
then the shape operator $S$ has the rank $\leq 1$.
\end{thm}

\par The main purpose of the present paper is to prove that the condition that the second fundamental form is Lorentzian can be dropped and
Theorem \ref{tw::RKomega0Lorentzian} and Theorem \ref{tw::NablaKOmega0Lorentzian} stay true for an arbitrary non-degenerate second fundamental form.
The main difficulty of this generalization comes from the fact that due to exponential grow of different possible scenarios methods from  \cite{MSz2} and \cite{MSzZSz}
cannot be easily repeated. For this reasons we need to change approach and first focus more on particular types of Jordan blocks rather than on all possible configurations.
Reducing number of allowed Jordan blocks dramatically decrease number of different configurations we need to consider when proving main theorems of this paper.

\par The paper is organized as follows. In Section 2 we briefly recall the basic formulas of affine differential
geometry, Jordan decomposition and some basic definitions from symplectic geometry. In this section we also prove a simple but important lemma about simultaneous decomposition of
the shape operator $S$ and the second fundamental form $h$.
The Section 3 is devoted to real Jordan blocks. The main result of this section is that the Jordan decomposition of the shape operator cannot contain
real Jordan blocks of dimension grater than or equal $3$ and it contains at most one block of dimension $2$.
In section 4 we study  complex Jordan blocks. It is shown that condition $R^p\omega=0$ implies that the Jordan decomposition of the shape operator cannot contain complex Jordan blocks.
Sections 5 contains the main results of this paper. Basing on results from sections 3 and 4, we show that if there exists an almost
symplectic structure $\omega$ satisfying condition $R^p\omega=0$ or $\nabla^p\omega=0$ for some positive integer $p$ then the rank of the shape operator must be $\leq 1$.
We conclude the section with some general example.


\section{Preliminaries}
\label{sec:1}
We briefly recall the basic formulas of affine differential
geometry. For more details, we refer to \cite{NomSas}. Let $f\colon M\rightarrow\R^{n+1}$ be an orientable
connected differentiable $n$-dimensional hypersurface immersed in
the affine space $\R^{n+1}$ equipped with its usual flat connection
$\D$. Then for any transversal vector field $\xi$ we have
\begin{equation}\label{eq::FormulaGaussa}
\D_Xf_\ast Y=f_\ast(\nabla_XY)+h(X,Y)\xi
\end{equation}
and
\begin{equation}\label{eq::FormulaWeingartena}
\D_X\xi=-f_\ast(SX)+\tau(X)\xi,
\end{equation}
where $X,Y$ are vector fields tangent to $M$. It is known that $\nabla$ is a torsion-free connection, $h$ is a symmetric
bilinear form on $M$, called \emph{the second
fundamental form}, $S$ is a tensor of type $(1,1)$, called \emph{the
shape operator}, and $\tau$ is a 1-form, called \emph{the transversal connection form}.
The vector field $\xi$ is called \emph{equiaffine} if $\tau=0$.  When $d\tau=0$ the vector field $\xi$ is called \emph{locally equiaffine}.

\par When $h$ is non-degenerate then $h$ defines a
pseudo-Rie\-man\-nian metric on $M$. In this case we say that the hypersurface or the
hypersurface immersion is \emph{non-degenerate}. In this paper we always assume that $f$ is
non-degenerate.  We have the following
\begin{thm}[\cite{NomSas}, Fundamental equations]\label{tw::FundamentalEquations}
For an arbitrary transversal vector field $\xi$ the induced
connection $\nabla$, the second fundamental form $h$, the shape
operator $S$, and the 1-form $\tau$ satisfy
the following equations:
\begin{align}
\label{eq::Gauss}&R(X,Y)Z=h(Y,Z)SX-h(X,Z)SY,\\
\label{eq::Codazzih}&(\nabla_X h)(Y,Z)+\tau(X)h(Y,Z)=(\nabla_Y h)(X,Z)+\tau(Y)h(X,Z),\\
\label{eq::CodazziS}&(\nabla_X S)(Y)-\tau(X)SY=(\nabla_Y S)(X)-\tau(Y)SX,\\
\label{eq::Ricci}&h(X,SY)-h(SX,Y)=2d\tau(X,Y).
\end{align}
\end{thm}
The equations (\ref{eq::Gauss}), (\ref{eq::Codazzih}),
(\ref{eq::CodazziS}), and (\ref{eq::Ricci}) are called the
equations of Gauss, Codazzi for $h$, Codazzi for $S$ and Ricci,
respectively.

\par Let $\omega$ be a non-degenerate 2-form on manifold $M$. The form $\omega$
we call an \emph{almost symplectic structure}. It is easy to see that if a manifold $M$ admits some almost symplectic structure then M is orientable manifold of even dimension.
Structure $\omega$ is called a \emph{symplectic structure}, if it is almost symplectic and additionally satisfies $d\omega=0$.
Pair $(M,\omega)$ we call \emph{a (an almost) symplectic manifold}, if $\omega$ is a (an almost) symplectic structure on $M$.

\par Recall (\cite{AlbPic}) that an affine connection $\nabla$ on an almost symplectic manifold $(M,\omega)$ we call an \emph{almost symplectic connection} if $\nabla\omega=0$.
An affine connection $\nabla$ on an almost symplectic manifold $(M,\omega)$ we call a \emph{symplectic connection} if it is almost symplectic and torsion-free.

\par Now we recall a well-know theorem about Jordan normal form (See eg. Th. A.2.6 in \cite{GohLanRod}).
\begin{thm}[\cite{GohLanRod}, Jordan]\label{tw::Jordana}
If $A\colon V\rightarrow V$ is an endomorphism of real finite dimensional vector space $V$ then there exists a basis of $V$ such that the matrix of
the endomorphism $A$ in this basis has a form
\begin{equation}
\left[\begin{matrix}
L_1 & 0 & \ldots & 0 \\
0 & L_2 & \ldots & 0 \\
\vdots & \vdots & \ddots & \ldots \\
0 & 0 & \ldots & L_s
\end{matrix}\right],
\end{equation}
where $L_i$ is the Jordan block corresponding to the eigenvalue $\lambda_i$ and given by the formula
\begin{equation}\label{eq::JordanRealEigenvalue}
\left[\begin{matrix}
\lambda_i & 0 & 0 & \ldots & 0 & 0 \\
1 & \lambda_i &  0 &\ldots & 0  & 0\\
0 & 1 & \lambda_i  &\ldots & 0  & 0\\
\vdots & \vdots & \vdots & \ddots & \ldots & \ldots\\
0 & 0 & 0 &\ldots & \lambda_i & 0\\
0 & 0 & 0 &\ldots & 1 &\lambda_i
\end{matrix}\right]\in M(k_i,k_i,\R),
\end{equation}
when $\lambda_i$ is a real number,
or by the formula
\begin{equation}\label{eq::JordanComplexEigenvalue}
\left[\begin{matrix}
B_i & 0 & 0 & \ldots & 0 & 0 \\
I & B_i &  0 &\ldots & 0  & 0\\
0 & I & B_i  &\ldots & 0  & 0\\
\vdots & \vdots & \vdots & \ddots & \ldots & \ldots\\
0 & 0 & 0 &\ldots & B_i & 0\\
0 & 0 & 0 &\ldots & I &B_i
\end{matrix}\right]\in M(2k_i,2k_i,\R),
\end{equation}
where
$$
B_i=
\left[
\begin{matrix}
\alpha_i & \beta_i\\
-\beta_i & \alpha_i
\end{matrix}
\right],\qquad
I=
\left[
\begin{matrix}
1 & 0\\
0 & 1
\end{matrix}
\right] \in M(2,2,\R),
$$
when $\lambda_i=\alpha_i+i\beta_i$ ($\beta_i\neq 0$) is a complex number.
\end{thm}

\par A square matrix $P$ of dimension $n$ is called \emph{the sip matrix (the standard involutory permutation)} (\cite{GohLanRod}) if it has a form:
\begin{equation}
\left[\begin{matrix}
0 & 0 & \cdots & 0 & 1 \\
0 & 0 & \cdots & 1 & 0 \\
\vdots & \vdots & \ddots & \cdots & \cdots \\
0 & 1 & \cdots & 0 & 0 \\
1 & 0 & \cdots & 0 & 0
\end{matrix}\right].
\end{equation}
Note that $P$ is non-singular symmetric matrix and $P^2=I$. In particular all its eigenvalues are equal $\pm 1$.
Moreover, it is easy to verify that we have the following formula for signature of $P$:
\begin{equation}\label{eq::sygnMacirzySip}
\operatorname{sig}P=
\begin{cases}
(\frac{n}{2},\frac{n}{2}), & \text{if $n$ is even} \\
(\frac{n+1}{2},\frac{n-1}{2}), & \text{if $n$ is odd.}
\end{cases}
\end{equation}
\begin{thm}[Th. 6.1.5 \cite{GohLanRod}]\label{tw::RozkladZSipMacirzami}
Let $H$ be a real invertible and symmetric matrix of dimension $n$. Then for every square $n$ dimensional
and $H$-selfadjoint matrix $A$ (i.e. $A^{T}H=HA$) there exists a basis  $\{e_1,\ldots,e_n\}$ such that
\begin{equation} \label{eq::RozkladA}
A=J_1\oplus\cdots \oplus J_t\oplus J_{t+1}\oplus\cdots\oplus J_{t+s},
\end{equation}
where $J_1,\ldots,J_t$ are Jordan blocks of type {\rm{(\ref{eq::JordanRealEigenvalue})}} and $J_{t+1},\ldots,J_{t+s}$ are Jordan blocks of type {\rm{(\ref{eq::JordanComplexEigenvalue})}}. Moreover
\begin{equation} \label{eq::RozkladH}
H=\eps_1 P_1\oplus\cdots \oplus \eps_t P_t\oplus P_{t+1}\oplus\cdots\oplus P_{t+s},
\end{equation}
where $P_j$ is a sip matrix of dimension equal to dimension of matrix $J_j$ for $j=1,\ldots,t+s$ and
$\eps_j=\pm 1$ for $j=1,\ldots,t$. The signs $\eps_j$ are determined uniquely by $(A,H)$ up to permutation of signs in the blocks of {\rm(\ref{eq::RozkladH})}
corresponding to the Jordan blocks of $A$ with the same real eigenvalue and the
same size.
\end{thm}

\par For a tensor field $T$ of type $(0,p)$ its covariant derivation $\nabla T$ is a tensor field of type $(0,p+1)$ given by the formula:
\begin{align*}
(\nabla T)(X_1,X_2,\ldots,X_{p+1}):=X_1(T(X_2,\ldots,X_{p+1}))\\ -\sum_{i=2}^{p+1}T(X_2,\ldots,\nabla_{X_1}X_i,\ldots,X_{p+1}).
\end{align*}
Higher order covariant derivatives of $T$ can be defined by recursion:
\begin{align*}
(\nabla^{k+1} T)=\nabla(\nabla^kT).
\end{align*}
To simplify computation it is often convenient to define $\nabla^0T:= T$.
\par If $R$ is a curvature tensor for an affine connection  $\nabla$, one can define a new tensor  $R\cdot T$ of type $(0,p+2)$ by the formula
\begin{align*}
(R\cdot T)(X_1,X_2,\ldots,X_{p+2}):= -\sum_{i=3}^{p+2}T(X_3,\ldots,R(X_1,X_2)X_i,\ldots,X_{p+2}).
\end{align*}
Analogously to the previous case, we may define a tensor $R^k\cdot T$ of type $(0,2k+p)$ using the following recursive formula:
$$
R^k\cdot T=R\cdot (R^{k-1}\cdot T)
$$
and additionally $R^0\cdot T:=T$.
\par In order to simplify the notation, we will be often omitting "$\cdot$" in $R^k\cdot T$ when no confusion arises.
Thus we will be writing often $R^k T$ instead of $R^k\cdot T$.

We conclude this section with the following lemma:

\begin{lem}\label{lm::JordanDecompositionOfS}
Let $f\colon \M\rightarrow\R^{2n+1}$ be a non-degenerate affine hypersurface with a locally equiaffine transversal vector field  $\xi$.
Then for every point $x\in M$ there exists a basis ${e_1,\ldots,e_{2n}}$ of $T_xM$
such that the shape operator $S$ and the second fundamental form  $h$ can be expressed in this basis in the block matrix form
\begin{equation}\label{eq::ShDecomposition}
S=\left[\begin{matrix}
S_1 & 0 & \ldots & 0 \\
0 & S_2 & \ldots & 0 \\
\vdots & \vdots & \ddots & \ldots \\
0 & 0 & \ldots & S_{q+r}
\end{matrix}\right],
h=\left[\begin{matrix}
H_1 & 0 & \ldots & 0 \\
0 & H_2 & \ldots & 0 \\
\vdots & \vdots & \ddots & \ldots \\
0 & 0 & \ldots & H_{q+r}
\end{matrix}\right],
\end{equation}
and $S_i$, $H_i$ satisfy the following conditions:
\begin{itemize}
  \item For $i=1,\ldots,q+r$ $\dim S_i=\dim H_i$.
  \item For $i=1,\ldots,q$ $S_i$ is a Jordan block of type (\ref{eq::JordanRealEigenvalue}) and $\dim S_i\geq \dim S_{i+1}$ for $i=1,\ldots,q-1$.
  \item For $i=q+1,\ldots,q+r$ $S_i$ is a Jordan block of type (\ref{eq::JordanComplexEigenvalue}) and $\dim S_i\geq \dim S_{i+1}$ for $i=q+1,\ldots,q+r-1$.
  \item For $i=1,\ldots,q$ $H_i$ is up to a sign sip matrix.
  \item For $i=q+1,\ldots,q+r$ $H_i$ is a sip matrix.
\end{itemize}
\begin{proof}
Since $\xi$ is locally equiaffine we have $d\tau=0$ and in consequence $h(SX,Y)=h(X,SY)$ for all $X,Y\in T_xM$.
Now thesis immediately follows from Theorem \ref{tw::RozkladZSipMacirzami} and the fact that we can rearrange Jordan blocks $S_i$ and matrices $H_i$
in desired order (rearranging vectors ${e_1,\ldots,e_{2n}}$ if needed).
\end{proof}
\end{lem}


\section{Real Jordan blocks}
In this chapter we study properties of real Jordan blocks of the shape operator $S$.

In all the below lemmas we assume that $f\colon \M\rightarrow\R^{2n+1}$ is a non-degenerate affine hypersurface with a locally equiaffine transversal vector field  $\xi$ and an almost symplectic form $\omega$.
About objects $\nabla$, $h$, $S$ and $\tau$ we assume that they are induced by $\xi$.

In the following lemmas (if not stated otherwise)
we assume that $S_1$ (from Lemma \ref{lm::JordanDecompositionOfS}) is a $k$-dimensional block of the form
\begin{align}\label{eq::S1withAlpha}
S_1=\left[\begin{matrix}
\alpha & 0 & 0 & \ldots & 0 & 0 \\
1 & \alpha &  0 &\ldots & 0  & 0\\
0 & 1 & \alpha  &\ldots & 0  & 0\\
\vdots & \vdots & \vdots & \ddots & \ldots & \ldots\\
0 & 0 & 0 &\ldots & \alpha & 0\\
0 & 0 & 0 &\ldots & 1 &\alpha
\end{matrix}\right]\in M(k,k,\R),
\end{align}
where $\alpha \in \R$
and
\begin{align*}
H_1=\left[\begin{matrix}
0 & 0 & \cdots & 0 & \varepsilon \\
0 & 0 & \cdots & \varepsilon & 0 \\
\vdots & \vdots & \ddots & \cdots & \cdots \\
0 & \varepsilon & \cdots & 0 & 0 \\
\varepsilon & 0 & \cdots & 0 & 0
\end{matrix}\right]\in M(k,k,\R),
\end{align*}
where $\varepsilon \in \{-1,1\}$.
By $\{e_1,\ldots,e_{2n}\}$ we will be denoting basis of $T_xM$
such that $\{e_1,\ldots,e_{k}\}$ is a basis for $S_1$.

\begin{lem}\label{lm::withPiX}
Let $p\geq 1$. If $k\geq 2$ and $X_1,\ldots,X_p\in \operatorname{span} \{e_1,\ldots,e_k\}=:V$  then for every $i\in \{2,3,\ldots,2n\}$
\begin{align}\label{eq::withPiX}
R^p\omega(\underbrace{X_1,e_k,X_2,e_k,\ldots,X_p,e_k}_{2p},e_i,e_k)=\pi(X_1)\cdot \ldots \pi(X_p)\varepsilon^p\alpha^p\omega(e_i,e_k),
\end{align}
where $\pi\colon V\rightarrow \mathbb{R}$ is a projection defined as follows:
$$
\pi(X)=\lambda_1
$$
if $X=\lambda_1e_1+\ldots +\lambda_ke_k\in V$ for some $\lambda_1,\ldots,\lambda_k\in\R$.

\begin{proof}
First we shall prove the following formulas:
\begin{align}
\label{eq::RXekek} R(X,e_k)e_k&=-\pi(X)\varepsilon\alpha e_k,\\
\label{eq::RXekei} R(X,e_k)e_i&=-h(X,e_i)\alpha e_k,\\
\label{eq::RXekY} \pi (R(X,e_k)Y)&=\varepsilon\alpha \pi(X)\pi(Y),
\end{align}
for all $X,Y \in V$.
To prove (\ref{eq::RXekek}) we compute
\begin{align*}
R(X,e_k)e_k&=\underbrace{h(e_k,e_k)}_{0}S_1X-h(X,e_k)S_1e_k=-h(\alpha_1e_1+\ldots+\alpha_ke_k,e_k)\alpha e_k\\
&=-\alpha_1\underbrace{h(e_1,e_k)}_{\varepsilon}\alpha e_k=-\pi(X)\varepsilon\alpha e_k.
\end{align*}
To prove (\ref{eq::RXekei}) we compute
\begin{align*}
R(X,e_k)e_i&=\underbrace{h(e_k,e_i)}_{0}S_1X-h(X,e_i)S_1e_k=-h(X,e_i)\alpha e_k.
\end{align*}
To prove (\ref{eq::RXekY}) we compute
\begin{align*}
\pi(R(X,e_k)Y)&=\pi(h(e_k,Y)S_1X-h(X,Y)S_1e_k)=h(e_k,Y)\pi(S_1X)-h(X,Y)\underbrace{\pi(S_1e_k)}_{0}\\
&=\varepsilon\pi(Y)\pi(X)\alpha.
\end{align*}
Now using formulas (\ref{eq::RXekek}), (\ref{eq::RXekei}), (\ref{eq::RXekY}) we shall prove the thesis of the lemma.
For $p=1$ we have
\begin{align*}
(R\omega )(X_1,e_k,e_i,e_k)&=-\omega (R(X_1,e_k)e_i,e_k)-\omega (e_i,R(X_1,e_k)e_k)\\
&=\underbrace{\omega (h(X_1,e_i)\alpha e_k,e_k)}_{0}+\omega (e_i,\pi(X_1)\varepsilon\alpha e_k)\\
&=\pi (X_1)\varepsilon\alpha \omega (e_i,e_k).
\end{align*}
Assume that formula (\ref{eq::withPiX}) is true for some $p\geq 1$.Then for $p+1$ we get
\begin{align*}
R^{p+1}\omega &(\underbrace{X_1,e_k,\ldots,X_{p+1},e_k}_{2p+2},e_i,e_k)\\
&=-R^p\omega (R(X_1,e_k)X_2,e_k,\ldots )\\
&\phantom{=}-R^p\omega (X_2,R(X_1,e_k)e_k,\ldots )\\
&\phantom{=}-R^p\omega (X_2,e_k,R(X_1,e_k)X_3,\ldots)\\
&\phantom{=}-R^p\omega (X_2,e_k,X_3,R(X_1,e_k)e_k,\ldots)\\
&\phantom{=} \cdots \\
&\phantom{=}-R^p\omega (X_2,e_k\ldots,R(X_1,e_k)e_i,e_k)\\
&\phantom{=}-R^p\omega (X_2,e_k,\ldots,e_i,R(X_1,e_k)e_k)\\
&=-\pi ((R(X_1,e_k)X_2)\pi (X_3)\ldots \pi(X_{p+1})\varepsilon^p\alpha^p\omega (e_i,e_k)\\
&\phantom{=}+R^p\omega (X_2,\pi(X_1)\varepsilon\alpha e_k,X_3,\ldots)\\
&\phantom{=}-\pi (X_2)\pi ((R(X_1,e_k)X_3)\pi (X_4)\ldots \pi(X_{p+1})\varepsilon^p\alpha^p\omega (e_i,e_k)\\
&\phantom{=}+R^p\omega (X_2,e_k,X_3,\pi(X_1)\varepsilon\alpha e_k,\ldots)\\
&\phantom{=} \cdots \\
&\phantom{=}-\pi(X_2)\ldots \pi(X_{p+1})\underbrace{\pi(R(X_1,e_k)e_i)}_{0}\varepsilon^p\alpha ^p \omega(e_i,e_k)\\
&\phantom{=} +R^p\omega (X_2,e_k,\ldots,e_i,\pi(X_1)\varepsilon\alpha e_k)\\
&=-\varepsilon\alpha\pi(X_1)\pi(X_2)\pi(X_3)\ldots \pi(X_{p+1})\varepsilon^p\alpha^p\omega(e_i,e_k)\\
&\phantom{=}+\varepsilon\alpha\pi(X_1)\pi(X_2)\pi(X_3)\ldots \pi(X_{p+1})\varepsilon^p\alpha^p\omega(e_i,e_k)\\
&\phantom{=}-\varepsilon\alpha\pi(X_1)\pi(X_2)\pi(X_3)\ldots \pi(X_{p+1})\varepsilon^p\alpha^p\omega(e_i,e_k)\\
&\phantom{=}+\varepsilon\alpha\pi(X_1)\pi(X_2)\pi(X_3)\ldots \pi(X_{p+1})\varepsilon^p\alpha^p\omega(e_i,e_k)\\
&\phantom{=} \cdots \\
&\phantom{=}-0\\
&\phantom{=} +\varepsilon\alpha\pi(X_1) \pi(X_2)\ldots \pi(X_{p+1})\varepsilon^p\alpha^p\omega(e_i,e_k)\\
&=\varepsilon^{p+1}\alpha^{p+1}\pi(X_1)\pi(X_2)\ldots \pi(X_{p+1})\omega(e_i,e_k).
\end{align*}
Now, by the induction principle the formula (\ref{eq::withPiX}) holds for every $p\geq 1$.
\end{proof}
\end{lem}

As an immediate consequence of Lemma \ref{lm::withPiX} (setting $X_1=X_2=\ldots=X_p=e_1$) we obtain:
\begin{cor}\label{cor::rpomegaeiek}
If $k\geq 2$ then for every $i\in\{2,\ldots,2n\}$ we have
\begin{align*}
R^p\omega (e_1,e_k,\ldots,e_1,e_k,e_i,e_k)=\varepsilon^p\alpha ^p\omega (e_i,e_k).
\end{align*}
\end{cor}

In the next few lemmas we shall obtain some properties of $R^p\omega$ under the assumption that $k>3$.


\begin{lem}\label{lm::basicForkGT3}
If $k>3$ then
\begin{align}
\label{eq::kGT3Basic1} R(e_1,e_{k-1})e_1&=R(e_1,e_{k-1})e_{k-1}=0\\
\label{eq::kGT3Basic2} R(e_1,e_{k-1})e_2&=\varepsilon S_1e_1=\varepsilon\alpha e_1+\varepsilon e_2\\
\label{eq::kGT3Basic3} R(e_1,e_{k-1})e_k&=-\varepsilon S_1e_{k-1}=-\varepsilon\alpha e_{k-1}-\varepsilon e_k\\
\label{eq::kGT3Basic4} R(e_{k-1},e_{k})e_1&=\varepsilon S_1e_{k-1}=\varepsilon\alpha e_{k-1}+\varepsilon e_k\\
\label{eq::kGT3Basic5} R(e_{k-1},e_{k})e_2&=-\varepsilon S_1e_{k}=-\varepsilon\alpha e_{k}\\
\label{eq::kGT3Basic6} R(e_{k-1},e_{k})e_{k-1}&=R(e_{k-1},e_{k})e_{k}=0
\end{align}
\begin{proof}
The proof is an immediate consequence of the Gauss equation and the fact that
\begin{align*}
h(e_1,e_k)=h(e_2,e_{k-1})=\varepsilon
\end{align*}
and
\begin{align*}
h(e_1,e_1)&=h(e_1,e_2)=h(e_1,e_{k-1})=h(e_{k-1},e_{k-1})\\&
=h(e_{k},e_{k-1})=h(e_{k},e_{k})=h(e_{2},e_{k})=0,
\end{align*}
if only $k>3$.
\end{proof}
\end{lem}


\begin{lem}\label{lm::Lemma34}
If $k>3$ and $p\geq 1$ we have
\begin{align}
\label{eq::Rpomegae1ekmniej1} R^p\omega (e_1,e_{k-1},e_1,e_{k-1},\ldots,e_1,e_{k-1})&=0,\\
\label{eq::Rpomegae2kmniej1} R^p\omega(e_1,e_{k-1},\ldots,e_1,e_{k-1},e_2,e_{k-1})&=
(-1)^p\varepsilon^p(\alpha\omega(e_1,e_{k-1})+\omega(e_2,e_{k-1})),\\
\label{eq::Rpomegaeiek} R^p\omega (e_1,e_{k-1},\ldots,e_1,e_{k-1},e_i,e_k)&=\varepsilon^p(\alpha\omega(e_i,e_{k-1})+\omega(e_i,e_k))\\
\nonumber \text{for $i\in\{1,\ldots,2n\}\setminus\{2,k\}$.}
\end{align}
\end{lem}
\begin{proof}
\par First note that the formula (\ref{eq::Rpomegae1ekmniej1}) easily follows from (\ref{eq::kGT3Basic1}) for every $p\geq 1$.
\par In order to prove (\ref{eq::Rpomegae2kmniej1}) let us note that for $p=1$ we have
\begin{align*}
R\omega (e_1,e_{k-1},e_2,e_{k-1})&=-\omega (R(e_1,e_{k-1})e_2,e_{k-1}))-\omega(e_2,R(e_1,e_{k-1})e_{k-1})\\
&=-\omega(\varepsilon\alpha e_1+\varepsilon e_2,e_{k-1})=-\varepsilon\alpha\omega(e_1,e_{k-1})-\varepsilon\omega (e_2,e_{k-1})
\end{align*}
thanks to (\ref{eq::kGT3Basic1}) and (\ref{eq::kGT3Basic2}).
Now assume that (\ref{eq::Rpomegae2kmniej1}) is true for some $p\geq 1$. Using (\ref{eq::kGT3Basic1}), (\ref{eq::kGT3Basic2}) and (\ref{eq::Rpomegae1ekmniej1}) we compute
\begin{align*}
R^{p+1}\omega &(\underbrace{e_1,e_{k-1},e_1,e_{k-1},\ldots,e_1,e_{k-1}}_{2p+2},e_2,e_{k-1})\\
&=(R(e_1,e_{k-1})\cdot R^p\omega)(\underbrace{e_1,e_{k-1},\ldots,e_1,e_{k-1}}_{2p},e_2,e_{k-1})\\
&=-R^p\omega (e_1,e_{k-1},\ldots,e_1,e_{k-1},R(e_1,e_{k-1})e_2,e_{k-1})\\
&=-R^p\omega (e_1,e_{k-1},\ldots,e_1,e_{k-1},\varepsilon S_1e_1,e_{k-1})\\
&=-\varepsilon R^p\omega (e_1,e_{k-1},\ldots,e_1,e_{k-1},\alpha e_1+e_2,e_{k-1})\\
&=-\varepsilon\alpha \underbrace{R^p\omega (e_1,e_{k-1},\ldots,e_1,e_{k-1},e_1,e_{k-1})}_{0}\\
&\phantom{=}-\varepsilon R^p\omega (e_1,e_{k-1},\ldots,e_1,e_{k-1},e_2,e_{k-1})\\
&=-\varepsilon(-1)^p\varepsilon^p(\alpha\omega(e_1,e_{k-1})+\omega(e_2,e_{k-1}))\\
&=(-1)^{p+1}\varepsilon^{p+1}(\alpha\omega(e_1,e_{k-1})+\omega(e_2,e_{k-1})).
\end{align*}
Thus, by the induction principle (\ref{eq::Rpomegae2kmniej1}) is true for all $p\geq 1$.
\par In order to prove (\ref{eq::Rpomegaeiek}) first note that if $e_i\bot\{e_1,e_{k-1}\}$ we have
\begin{align}
\label{eq::Re1ekeiFormula}
R(e_1,e_{k-1})e_i=h(e_{k-1},e_i)S_1e_1-h(e_1,e_i)S_1e_{k-1}=0.
\end{align}
In particular the above holds for all $i\in\{1,\ldots,2n\}\setminus\{2,k\}$.
Using (\ref{eq::kGT3Basic3}) and (\ref{eq::Re1ekeiFormula}) we get that for every $i\in\{1,\ldots,2n\}\setminus\{2,k\}$
\begin{align*}
R\omega(e_1,e_{k-1},e_i,e_k)&=-\omega(\underbrace{R(e_1,e_{k-1})e_i}_{0},e_k)-\omega(e_i,R(e_1,e_{k-1})e_k)\\
&=\varepsilon(\alpha \omega(e_i,e_{k-1})+\omega(e_i,e_k))
\end{align*}
thus (\ref{eq::Rpomegaeiek}) is true for $p=1$.
Now assume that (\ref{eq::Rpomegaeiek}) holds for some $p\geq 1$ and all $i\in\{1,\ldots,2n\}\setminus\{2,k\}$. We have
\begin{align*}
R^{p+1}\omega &(\underbrace{e_1,e_{k-1},e_1,e_{k-1},\ldots,e_1,e_{k-1}}_{2p+2},e_i,e_k)\\
&=(R(e_1,e_{k-1})\cdot R^p\omega)(\underbrace{e_1,e_{k-1},\ldots,e_1,e_{k-1}}_{2p},e_i,e_k)\\
&=-R^p\omega (e_1,e_{k-1},\ldots,e_1,e_{k-1},e_i,R(e_1,e_{k-1})e_k)\\
&=\varepsilon\alpha R^p\omega (e_1,e_{k-1},\ldots,e_1,e_{k-1},e_i,e_{k-1})\\
&\phantom{=}+\varepsilon R^p\omega (e_1,e_{k-1},\ldots,e_1,e_{k-1},e_i,e_k)\\
&=\varepsilon R^p\omega (e_1,e_{k-1},\ldots,e_1,e_{k-1},e_i,e_k)\\
&=\varepsilon^{p+1}(\alpha\omega(e_i,e_{k-1})+\omega(e_i,e_k))
\end{align*}
since
\begin{align*}
R(e_1,e_{k-1})e_1=R(e_1,e_{k-1})e_{k-1}=R(e_1,e_{k-1})e_i=0
\end{align*}
by (\ref{eq::kGT3Basic1}) and (\ref{eq::Re1ekeiFormula}).
Now, by the induction principle (\ref{eq::Rpomegaeiek}) holds for all $p\geq 1$.
\end{proof}


\begin{lem}\label{lm::Revenandodd}
If $k>3$ and $p\geq 0$ we have
\begin{align}
\label{eq::R2pplus1omega} &R^{2p+1}\omega (\underbrace{e_1,e_{k-1},\ldots,e_1,e_{k-1}}_{4p+2},e_2,e_k)=-\varepsilon\alpha (\omega (e_1,e_k)-\omega (e_2,e_{k-1})),\\
\label{eq::R2pplus2omega} &R^{2p+2}\omega (\underbrace{e_1,e_{k-1},\ldots,e_1,e_{k-1}}_{4p+4},e_2,e_k)\\
\nonumber &\quad\qquad\qquad\qquad\qquad=-\alpha (2\alpha\omega (e_1,e_{k-1})+\omega (e_2,e_{k-1})+\omega (e_1,e_{k})).
\end{align}
\begin{proof}
First note that by straightforward computations, using (\ref{eq::kGT3Basic1}), (\ref{eq::kGT3Basic2}) and (\ref{eq::kGT3Basic3}), one may easily check that (\ref{eq::R2pplus1omega}) and (\ref{eq::R2pplus2omega}) are true for $p=0$.
\par Let us assume that $p>0$.
In order to prove (\ref{eq::R2pplus1omega}), using Lemma \ref{lm::basicForkGT3}, we compute
\begin{align*}
 R^{2p+1}\omega & (\underbrace{e_1,e_{k-1},\ldots,e_1,e_{k-1}}_{4p+2},e_2,e_k)\\
&=(R(e_1,e_{k-1})\cdot R^{2p}\omega)(\underbrace{e_1,e_{k-1},\ldots,e_1,e_{k-1}}_{4p},e_2,e_k)\\
&=-R^{2p}\omega(e_1,e_{k-1},\ldots,e_1,e_{k-1},R(e_1,e_{k-1})e_2,e_k)\\
&\phantom{=}-R^{2p}\omega(e_1,e_{k-1},\ldots,e_1,e_{k-1},e_2,R(e_1,e_{k-1})e_k)\\
&=-\varepsilon R^{2p}\omega(e_1,e_{k-1},\ldots,e_1,e_{k-1},\alpha e_1+e_2,e_k)\\
&\phantom{=}-\varepsilon R^{2p}\omega(e_1,e_{k-1},\ldots,e_1,e_{k-1},e_2,-\alpha e_{k-1}-e_k)\\
&=-\varepsilon \alpha R^{2p}\omega(e_1,e_{k-1},\ldots,e_1,e_{k-1},e_1,e_k)\\
&\phantom{=}-\varepsilon R^{2p}\omega(e_1,e_{k-1},\ldots,e_1,e_{k-1},e_2,e_k)\\
&\phantom{=}+\varepsilon\alpha R^{2p}\omega(e_1,e_{k-1},\ldots,e_1,e_{k-1},e_2,e_{k-1})\\
&\phantom{=}+\varepsilon R^{2p}\omega(e_1,e_{k-1},\ldots,e_1,e_{k-1},e_2,e_k)\\
&=-\varepsilon \alpha R^{2p}\omega(e_1,e_{k-1},\ldots,e_1,e_{k-1},e_1,e_k)\\
&\phantom{=}+\varepsilon\alpha R^{2p}\omega(e_1,e_{k-1},\ldots,e_1,e_{k-1},e_2,e_{k-1}).
\end{align*}
Now using (\ref{eq::Rpomegaeiek}) (for $i=1$) and (\ref{eq::Rpomegae2kmniej1}) we obtain
\begin{align*}
R^{2p+1}\omega & (\underbrace{e_1,e_{k-1},\ldots,e_1,e_{k-1}}_{4p+2},e_2,e_k)\\
&=-\varepsilon \alpha R^{2p}\omega(e_1,e_{k-1},\ldots,e_1,e_{k-1},e_1,e_k)\\
&\phantom{=}+\varepsilon\alpha R^{2p}\omega(e_1,e_{k-1},\ldots,e_1,e_{k-1},e_2,e_{k-1})\\
&=-\varepsilon\alpha \varepsilon^{2p}(\alpha \omega(e_1,e_{k-1})+\omega (e_1,e_k))+\varepsilon\alpha\cdot (-1)^{2p}\varepsilon^{2p} (\alpha \omega(e_1,e_{k-1})+\omega (e_2,e_{k-1}))\\
&=-\varepsilon\alpha (\omega(e_1,e_k)-\omega(e_2,e_{k-1}))
\end{align*}
that is  (\ref{eq::R2pplus1omega}) is true for all $p\geq 0$.
\par In order to prove (\ref{eq::R2pplus2omega}) we compute
\begin{align*}
R^{2p+2}\omega &(\underbrace{e_1,e_{k-1},\ldots,e_1,e_{k-1}}_{4p+4},e_2,e_k) \\
&= -R^{2p+1}\omega (R(e_1,e_{k-1})e_1,e_{k-1},\ldots,e_1,e_{k-1},e_2,e_k)\\
&\phantom{=} -R^{2p+1}\omega (e_1,R(e_1,e_{k-1})e_{k-1},\ldots,e_1,e_{k-1},e_2,e_k)\\
&\phantom{=}\cdots\\
&\phantom{=} -R^{2p+1}\omega (e_1,e_{k-1},\ldots,e_1,e_{k-1},R(e_1,e_{k-1})e_2,e_k)\\
&\phantom{=} -R^{2p+1}\omega (e_1,e_{k-1},\ldots,e_1,e_{k-1},e_2,R(e_1,e_{k-1})e_k)\\
&= -R^{2p+1}\omega (e_1,e_{k-1},\ldots,e_1,e_{k-1},R(e_1,e_{k-1})e_2,e_k)\\
&\phantom{=} -R^{2p+1}\omega (e_1,e_{k-1},\ldots,e_1,e_{k-1},e_2,R(e_1,e_{k-1})e_k)
\end{align*}
where the last equality follows from (\ref{eq::kGT3Basic1}). Now using (\ref{eq::kGT3Basic2}) and (\ref{eq::kGT3Basic3}) we get
\begin{align*}
R^{2p+2}\omega &(\underbrace{e_1,e_{k-1},\ldots,e_1,e_{k-1}}_{4p+4},e_2,e_k) \\
&=-\varepsilon\alpha R^{2p+1}\omega (e_1,e_{k-1},\ldots,e_1,e_{k-1},e_1,e_k)\\
&\phantom{=}-\varepsilon R^{2p+1}\omega (e_1,e_{k-1},\ldots,e_1,e_{k-1},e_2,e_k)\\
&\phantom{=}+\varepsilon\alpha R^{2p+1}\omega (e_1,e_{k-1},\ldots,e_1,e_{k-1},e_2,e_{k-1})\\
&\phantom{=}+\varepsilon R^{2p+1}\omega (e_1,e_{k-1},\ldots,e_1,e_{k-1},e_2,e_k)\\
&=-\varepsilon\alpha R^{2p+1}\omega (e_1,e_{k-1},\ldots,e_1,e_{k-1},e_1,e_k)\\
&\phantom{=}+\varepsilon\alpha R^{2p+1}\omega (e_1,e_{k-1},\ldots,e_1,e_{k-1},e_2,e_{k-1}).
\end{align*}
Again using (\ref{eq::Rpomegaeiek}) (for $i=1$) and (\ref{eq::Rpomegae2kmniej1}) we obtain
\begin{align*}
R^{2p+2}\omega &(\underbrace{e_1,e_{k-1},\ldots,e_1,e_{k-1}}_{4p+4},e_2,e_k) \\
&= -\varepsilon\alpha\varepsilon^{2p+1}(\alpha\omega(e_1,e_{k-1})+\omega(e_1,e_k))\\
&\phantom{=}+\varepsilon\alpha\cdot (-1)^{2p+1}\varepsilon^{2p+1}(\alpha\omega(e_1,e_{k-1})+\omega(e_2,e_{k-1}))\\
&=-\alpha (2\alpha\omega (e_1,e_{k-1})+\omega (e_2,e_{k-1})+\omega (e_1,e_{k})),
\end{align*}
what completes the proof of (\ref{eq::R2pplus2omega}).
\end{proof}
\end{lem}


\begin{lem}\label{lm::Lemma36}
If $k>3$ and $p\geq 1$ we have
\begin{align}\label{lem::rpomegae1e2}
R^p\omega &(e_{k-1},e_k,\underbrace{e_1,e_{k-1},\ldots,e_1,e_{k-1}}_{2p-2},e_1,e_2) \\
\nonumber &= \varepsilon^p\alpha (\omega (e_1,e_k)+\omega(e_2,e_{k-1}))+\varepsilon^p\omega(e_2,e_k)
\end{align}
\begin{proof}
For $p=1$, using (\ref{eq::kGT3Basic4}) and (\ref{eq::kGT3Basic5}), we directly check that
\begin{align*}
R\omega(e_{k-1},e_k,e_1,e_2)=\varepsilon\alpha (\omega (e_1,e_k)+\omega(e_2,e_{k-1}))+\varepsilon\omega(e_2,e_k).
\end{align*}
Now assume that (\ref{lem::rpomegae1e2}) is true for some $p\geq 1$. First we compute that
\begin{align*}
R^{p+1}\omega &(e_{k-1},e_k,\underbrace{e_1,e_{k-1},\ldots,e_1,e_{k-1}}_{2p},e_1,e_2)\\
&=(R(e_{k-1},e_k)\cdot R^p\omega)(\underbrace{e_1,e_{k-1},\ldots,e_1,e_{k-1}}_{2p},e_1,e_2)\\
&=-R^p\omega(R(e_{k-1},e_k)e_1,e_{k-1},\ldots,e_1,e_{k-1},e_1,e_2)\\
&\phantom{=}-R^p\omega(e_1,R(e_{k-1},e_k)e_{k-1},\ldots,e_1,e_{k-1},e_1,e_2)\\
&\phantom{=}-\ldots \\
&\phantom{=}-R^p\omega(e_1,e_{k-1},\ldots,R(e_{k-1},e_k)e_1,e_{k-1},e_1,e_2)\\
&\phantom{=}-R^p\omega(e_1,e_{k-1},\ldots,e_1,R(e_{k-1},e_k)e_{k-1},e_1,e_2)\\
&\phantom{=}-R^p\omega(e_1,e_{k-1},\ldots,e_1,e_{k-1},R(e_{k-1},e_k)e_1,e_2)\\
&\phantom{=}-R^p\omega(e_1,e_{k-1},\ldots,e_1,e_{k-1},e_1,R(e_{k-1},e_k)e_2)\\
&=-\varepsilon R^p\omega(e_k,e_{k-1},\ldots,e_1,e_{k-1},e_1,e_2)\\
&\phantom{=}-\varepsilon R^p\omega(e_1,e_{k-1},e_k,e_{k-1},\ldots,e_1,e_{k-1},e_1,e_2)\\
&\phantom{=}-\ldots \\
&\phantom{=}-\varepsilon R^p\omega(e_1,e_{k-1},\ldots,e_k,e_{k-1},e_1,e_2)\\
&\phantom{=}-R^p\omega(e_1,e_{k-1},\ldots,e_1,e_{k-1},\varepsilon \alpha e_{k-1}+\varepsilon e_k,e_2)\\
&\phantom{=}-R^p\omega(e_1,e_{k-1},\ldots,e_1,e_{k-1},e_1,-\varepsilon \alpha e_k).
\end{align*}
Before we proceed we shall show that for $p\geq 2$
\begin{align*}
R^p\omega(\overbrace{e_1,e_{k-1}}^{(1)},\ldots,\overbrace{e_k,e_{k-1}}^{(i)},\ldots,\overbrace{e_1,e_{k-1}}^{(p)},e_1,e_2)=0
\end{align*}
if only $i\in\{2,\ldots,p\}$. Indeed we have
\begin{align*}
R^p\omega&(\overbrace{e_1,e_{k-1}}^{(1)},\ldots,\overbrace{e_k,e_{k-1}}^{(i)},\ldots,\overbrace{e_1,e_{k-1}}^{(p)},e_1,e_2)\\
&=-R^{p-1}\omega(R(e_1,e_{k-1})e_1,e_{k-1},\ldots)\\
&\phantom{=}-R^{p-1}\omega(e_1,R(e_1,e_{k-1})e_{k-1},\ldots)\\
&\phantom{=}-\cdots\\
&\phantom{=}-R^{p-1}\omega(\ldots,\overbrace{R(e_1,e_{k-1})e_k,e_{k-1}}^{(i-1)},\ldots)\\
&\phantom{=}-R^{p-1}\omega(\ldots,\overbrace{e_k,R(e_1,e_{k-1})e_{k-1}}^{(i-1)},\ldots)\\
&\phantom{=}-\cdots\\
&\phantom{=}-R^{p-1}\omega(\ldots,R(e_1,e_{k-1})e_1,e_2)\\
&\phantom{=}-R^{p-1}\omega(\ldots,e_1,R(e_1,e_{k-1})e_2)\\
&=-R^{p-1}\omega(\ldots,\overbrace{R(e_1,e_{k-1})e_k,e_{k-1}}^{(i-1)},\ldots)\\
&\phantom{=}-R^{p-1}\omega(\ldots,e_1,R(e_1,e_{k-1})e_2)
\end{align*}
thanks to (\ref{eq::kGT3Basic1}). Now formulas (\ref{eq::kGT3Basic2}) and (\ref{eq::kGT3Basic3}) imply that
\begin{align}\label{eq::RpOmega_1ip}
R^p\omega(\overbrace{e_1,e_{k-1}}^{(1)},\ldots,\overbrace{e_k,e_{k-1}}^{(i)},\ldots,\overbrace{e_1,e_{k-1}}^{(p)},e_1,e_2)=0.
\end{align}
Now applying (\ref{eq::RpOmega_1ip}) if needed (that is if $p>1$) we conclude that
\begin{align*}
R^{p+1}\omega &(e_{k-1},e_k,\underbrace{e_1,e_{k-1},\ldots,e_1,e_{k-1}}_{2p},e_1,e_2)\\
&=\varepsilon R^p\omega(e_{k-1},e_k,e_1,e_{k-1},\ldots,e_1,e_{k-1},e_1,e_2)\\
&\phantom{=}-\varepsilon \alpha R^p\omega(e_1,e_{k-1},\ldots,e_1,e_{k-1},e_{k-1},e_2)\\
&\phantom{=}-\varepsilon R^p\omega(e_1,e_{k-1},\ldots,e_1,e_{k-1},e_k,e_2)\\
&\phantom{=}+\varepsilon \alpha R^p\omega (e_1,e_{k-1},\ldots,e_1,e_{k-1},e_1,e_k)\\
&=\varepsilon R^p\omega(e_{k-1},e_k,e_1,e_{k-1},\ldots,e_1,e_{k-1},e_1,e_2)\\
&\phantom{=}+\varepsilon \alpha R^p\omega(e_1,e_{k-1},\ldots,e_1,e_{k-1},e_2,e_{k-1})\\
&\phantom{=}+\varepsilon R^p\omega(e_1,e_{k-1},\ldots,e_1,e_{k-1},e_2,e_k)\\
&\phantom{=}+\varepsilon \alpha R^p\omega (e_1,e_{k-1},\ldots,e_1,e_{k-1},e_1,e_k).\\
\end{align*}
By the induction principle, using formula (\ref{eq::Rpomegae2kmniej1}) and formula (\ref{eq::Rpomegaeiek}) (for $i=1$) we obtain
\begin{align*}
R^{p+1}\omega &(e_{k-1},e_k,\underbrace{e_1,e_{k-1},\ldots,e_1,e_{k-1}}_{2p},e_1,e_2)\\
&=\varepsilon\Big(\varepsilon^p\alpha (\omega (e_1,e_k)+\omega(e_2,e_{k-1}))+\varepsilon^p\omega(e_2,e_k)\Big)\\
&\phantom{=}+\varepsilon\alpha\Big( (-1)^p\varepsilon^p(\alpha\omega(e_1,e_{k-1})+\omega(e_2,e_{k-1}))\Big)\\
&\phantom{=}+\varepsilon R^p\omega(e_1,e_{k-1},\ldots,e_1,e_{k-1},e_2,e_k)\\
&\phantom{=}+\varepsilon\alpha \Big(\varepsilon^p(\alpha\omega(e_1,e_{k-1})+\omega(e_1,e_k))\Big)\\
&=\varepsilon^{p+1}\Big(2\alpha\omega(e_1,e_k)+(\alpha+(-1)^p\alpha)\omega(e_2,e_{k-1})\\
&\qquad\qquad+(\alpha^2+(-1)^p\alpha^2)\omega(e_1,e_{k-1})+\omega(e_2,e_{k})\Big)\\
&\phantom{=}+\varepsilon R^p\omega(e_1,e_{k-1},\ldots,e_1,e_{k-1},e_2,e_k).
\end{align*}
Finally we need to consider two cases. If $p$ is odd, using Lemma \ref{lm::Revenandodd}, we obtain
\begin{align*}
R^{p+1}\omega &(e_{k-1},e_k,\underbrace{e_1,e_{k-1},\ldots,e_1,e_{k-1}}_{2p},e_1,e_2)\\
&=\varepsilon^{p+1}\Big(2\alpha\omega(e_1,e_k)+(\alpha+(-1)^p\alpha)\omega(e_2,e_{k-1})\\
&\qquad\qquad+(\alpha^2+(-1)^p\alpha^2)\omega(e_1,e_{k-1})+\omega(e_2,e_{k})\Big)\\
&\phantom{=}+\varepsilon (-\varepsilon\alpha (\omega (e_1,e_k)-\omega (e_2,e_{k-1})))\\
&=\alpha(\omega(e_1,e_{k})+\omega(e_2,e_{k-1}))+\omega(e_2,e_{k}).
\end{align*}
If $p$ is even, again from Lemma \ref{lm::Revenandodd} we get
\begin{align*}
R^{p+1}\omega &(e_{k-1},e_k,\underbrace{e_1,e_{k-1},\ldots,e_1,e_{k-1}}_{2p},e_1,e_2)\\
&=\varepsilon^{p+1}\Big(2\alpha\omega(e_1,e_k)+(\alpha+(-1)^p\alpha)\omega(e_2,e_{k-1})\\
&\qquad\qquad+(\alpha^2+(-1)^p\alpha^2)\omega(e_1,e_{k-1})+\omega(e_2,e_{k})\Big)\\
&\phantom{=}+\varepsilon (-\alpha (2\alpha\omega (e_1,e_{k-1})+\omega (e_2,e_{k-1})+\omega (e_1,e_{k})))\\
&=\varepsilon \alpha (\omega (e_1,e_k)+\omega(e_2,e_{k-1}))+\varepsilon \omega(e_2,e_k).
\end{align*}
The proof is completed.
\end{proof}
\end{lem}


\begin{lem}\label{lm::LematD}
If $k>3$ then for every $p\geq 1$ we have
\begin{align}
\label{eq::Lem::D} R^{p}\omega (e_1,e_{k-1},e_1,e_{k-1},\ldots,e_1,e_{k-1},e_1,e_3)=0
\end{align}
\end{lem}
\begin{proof}
If $p=1$ we have
\begin{align*}
R\omega &(e_1,e_{k-1},e_1,e_3)=-\omega (R(e_1,e_{k-1})e_1,e_3)-\omega (e_1,R(e_1,e_{k-1})e_3).
\end{align*}
For $p>1$ we have in general
\begin{align*}
R^{p}\omega &(e_1,e_{k-1},e_1,e_{k-1},\ldots,e_1,e_{k-1},e_1,e_3)\\
&=-R^{p-1}\omega (R(e_1,e_{k-1})e_1,e_{k-1},\ldots,e_1,e_{k-1},e_1,e_3)\\
&\phantom{=}-R^{p-1}\omega (e_1,R(e_1,e_{k-1})e_{k-1},\ldots,e_1,e_{k-1},e_1,e_3)\\
&\phantom{=}-\ldots \\
&\phantom{=}-R^{p-1}\omega (e_1,e_{k-1},\ldots,e_1,R(e_1,e_{k-1})e_{k-1},e_1,e_3)\\
&\phantom{=}-R^{p-1}\omega (e_1,e_{k-1},\ldots,e_1,e_{k-1},R(e_1,e_{k-1})e_1,e_3)\\
&\phantom{=}-R^{p-1}\omega (e_1,e_{k-1},\ldots,e_1,e_{k-1},e_1,R(e_1,e_{k-1})e_3).
\end{align*}
Now the thesis follows immediately from (\ref{eq::kGT3Basic1}) and the fact that for $k>3$ we always have $R(e_1,e_{k-1})e_3=0$.
\end{proof}


\begin{lem}\label{lm::S1dimMax3}
Let $f\colon \M\rightarrow\R^{2n+1}$ ($\dim M\geq 4$) be a non-degenerate affine hypersurface with a locally equiaffine transversal vector field  $\xi$
and an almost symplectic form $\omega$. If $R^p\omega =0$ for some $p\geq 1$ and
\begin{align}\label{eq::SDecompositionComplex}
S=\left[\begin{matrix}
S_1 & 0 & \ldots & 0 \\
0 & S_2 & \ldots & 0 \\
\vdots & \vdots & \ddots & \ldots \\
0 & 0 & \ldots & S_{q+r}
\end{matrix}\right]
\end{align}
is the Jordan decomposition of $S$ as stated in the Lemma \ref{lm::JordanDecompositionOfS} then
$\dim S_1\leq 3$.

\begin{proof}
Let us assume that $S_1$ has the form (\ref{eq::S1withAlpha}) and $k>3$.
If $\alpha \neq 0$,
from Corollary \ref{cor::rpomegaeiek} we obtain
\begin{align*}
\omega (e_2,e_k)=\ldots =\omega (e_{2n},e_k)=0.
\end{align*}
From Lemma \ref{lm::Revenandodd}, formula (\ref{eq::R2pplus1omega}) we have
\begin{align*}
\omega (e_1,e_k)=\omega (e_2,e_{k-1}).
\end{align*}
From Lemma \ref{lm::Lemma36}, formula (\ref{lem::rpomegae1e2}) we get
\begin{align*}
\omega (e_1,e_k)=-\omega (e_2,e_{k-1}),
\end{align*}
since $\omega (e_2,e_k)=0$.
Therefore $\omega (e_i,e_k)=0$ for $i\in \{1,\ldots,2n\}$
what contradicts assumption that $\omega$ is non-degenerate.
Thus it must be $\alpha =0$.
\par When $\alpha =0$ from Lemma \ref{lm::Lemma34}, formula (\ref{eq::Rpomegaeiek}) we get $\omega(e_i,e_k)=0$ for $i\in \{1,2,\ldots,2n\}\setminus\{2,k\}$.
Of course $\omega(e_k,e_k)=0$ as well.
From Lemma \ref{lm::Lemma36}, formula (\ref{lem::rpomegae1e2}) we have that also $\omega(e_2,e_k)=0$, so again we obtain that $\omega$ is degenerate
thus $\dim S_1$ cannot exceed $3$.
\end{proof}
\end{lem}


\begin{lem}\label{lm::Block3x3omega12}
If $S_1$ is a 3-dimensional real Jordan block, $p\geq 1$ then
\begin{align}\label{eq::Block3x3omega12}
R^p\omega(e_1,e_2,e_1,e_2,\ldots,e_1,e_2) =(-1)^p\varepsilon^p p!\omega(e_1,e_2)
\end{align}
\begin{proof}
From the Gauss equation we have
\begin{align}\label{eq::RforS1dim3}
R(e_1,e_2)e_1=0, \quad R(e_1,e_2)e_2=\varepsilon S_1e_1=\varepsilon\alpha e_1+\varepsilon e_2.
\end{align}
Now, for $p=1$ we easily check that
$$
R\omega(e_1,e_2,e_1,e_2)=-\varepsilon\omega(e_1,e_2).
$$
Let us assume that (\ref{eq::Block3x3omega12}) is true for some $p\geq 1$.
Using (\ref{eq::RforS1dim3}) we compute
\begin{align*}
R^{p+1}\omega(e_1&,e_2,e_1,e_2,\ldots,e_1,e_2)\\
&=(R(e_1,e_{2})\cdot R^p\omega)(\underbrace{e_1,e_{2},\ldots,e_1,e_{2}}_{2p},e_1,e_2)\\
&=-R^p\omega (e_1,\varepsilon S_1e_1,\ldots,e_1,e_{2},e_1,e_2)\\
&\phantom{=}-R^p\omega (e_1,e_2,e_1,\varepsilon S_1e_1,\ldots,e_1,e_2)\\
&\phantom{=}\ldots \\
&\phantom{=}-R^p\omega (e_1,e_2,e_1,e_2,\ldots,e_1,\varepsilon S_1e_1,e_1,e_2)\\
&\phantom{=}-R^p\omega (e_1,e_2,e_1,e_2,\ldots,e_1,e_2,e_1,\varepsilon S_1e_1)\\
&=-\varepsilon (p+1)R^p\omega (e_1,e_2,e_1,e_2,\ldots,e_1,e_2)\\
&=-\varepsilon (p+1)\cdot (-1)^p\varepsilon^p p!\omega (e_1,e_2)=(-1)^{p+1}\varepsilon^{p+1} (p+1)!\omega (e_1,e_2)
\end{align*}
Now by the induction principle (\ref{eq::Block3x3omega12}) holds for all $p\geq 1$.
\end{proof}
\end{lem}


\begin{lem}\label{lm::Block3x3omega12ij}
If $S_1$ is a 3-dimensional real Jordan block, $p\geq 2$ and $i,j>3$ then
\begin{align}\label{eq::Block3x3omega12ij}
R^p\omega&(e_1,e_2,e_1,e_2,\ldots,e_1,e_2,e_2,e_i,e_1,e_j)\\
\nonumber &=(-1)^p\varepsilon^{p-1}  (p-1)!h(e_i,e_j)(2\alpha\omega(e_1,e_2)+\omega(e_1,e_3))
\end{align}
\begin{proof}
Using the Gauss equetion let us note  that
\begin{align}\label{eq::Block3x3proof1}
R^2\omega(e_1,e_2,e_2,e_i,e_1,e_j)&=\varepsilon h(e_i,e_j)(2\alpha\omega(e_1,e_2)+\omega(e_1,e_3)),\\
\label{eq::Block3x3proof2} R^2\omega(e_1,e_2,e_1,e_i,e_1,e_j)&=0
\end{align}
Now, for $p\geq 2$ we have
$$
R^{p+1}\omega(e_1,e_2,\ldots,e_1,e_2,e_1,e_i,e_1,e_j)=-\varepsilon(p-1)R^{p}\omega(e_1,e_2,\ldots,e_1,e_2,e_1,e_i,e_1,e_j)
$$
since $R(e_1,e_2)e_1=R(e_1,e_2)e_i=R(e_1,e_2)e_j=0$ and $R(e_1,e_2)e_2=\varepsilon S_1 e_1=\varepsilon\alpha e_1+\varepsilon e_2$.
In consequence, by the induction principle
\begin{align}\label{eq::Block3x3proof3}
R^{p}\omega(e_1,e_2,\ldots,e_1,e_2,e_1,e_i,e_1,e_j)=0
\end{align}
for all $p\geq 2$.
We also compute
\begin{align*}
R^{p+1}\omega(e_1,e_2,\ldots,e_1,e_2,e_2,e_i,e_1,e_j)&=-\varepsilon(p-1)R^{p}\omega(e_1,e_2,\ldots,e_1,e_2,e_2,e_i,e_1,e_j)\\
&\phantom{=}-R^{p}\omega(e_1,e_2,\ldots,e_1,e_2,R(e_1,e_2)e_2,e_i,e_1,e_j)\\
&=-\varepsilon p R^{p}\omega(e_1,e_2,\ldots,e_1,e_2,e_2,e_i,e_1,e_j)\\
&\phantom{=}-\varepsilon\alpha R^{p}\omega(e_1,e_2,\ldots,e_1,e_2,e_1,e_i,e_1,e_j)\\
&=-\varepsilon p R^{p}\omega(e_1,e_2,\ldots,e_1,e_2,e_2,e_i,e_1,e_j),
\end{align*}
where the last equality follows from (\ref{eq::Block3x3proof3}).
The induction principle implies that (\ref{eq::Block3x3omega12ij}) holds for all $p\geq 2$.
\end{proof}
\end{lem}


\begin{lem}\label{lm::Block3x3omega122i}
If $S_1$ is a 3-dimensional real Jordan block, $p\geq 1$ and $\alpha =0$ then
\begin{align}\label{eq::Block3x3omega122i}
R^p\omega(e_1,e_2,e_1,e_2,\ldots,e_1,e_2,e_2,e_i) =(-1)^p\varepsilon^p p!\omega (e_2,e_i)
\end{align}
for any $i \in \{1,\ldots,2n\}\setminus \{3\}$.
\begin{proof}
First notice that (\ref{eq::Block3x3omega122i}) holds trivially  for $i=2$. Now we assume that $i \in \{1,\ldots,2n\}\setminus \{2,3\}$.
Note that from  the Gauss equation and assumption $\alpha=0$ we have
\begin{align}\label{eq::RfromGausswithAlpha0}
R(e_1,e_2)e_1=R(e_1,e_2)e_i=0, \quad R(e_1,e_2)e_2=\varepsilon S_1e_1=\varepsilon e_2.
\end{align}
\par For $p=1$ we easily get
\begin{align*}
R\omega(e_1,e_2,e_2,e_i)=-\varepsilon\omega (e_2,e_i).
\end{align*}
Let us assume that (\ref{eq::Block3x3omega122i}) holds for some $p\geq 1$, we shall show that it also holds for $p+1$.
Indeed we obtain
\begin{align*}
R^{p+1}&\omega(\underbrace{e_1,e_2,\ldots,e_1,e_2}_{2p+2},e_2,e_i)\\
&=(R(e_1,e_2)\cdot R^{p}\omega)(\underbrace{e_1,e_2,\ldots,e_1,e_2}_{2p},e_2,e_i)\\
&=-\varepsilon(p+1)R^p\omega (e_1,e_2,\ldots,e_1,e_2,e_2,e_i)\\
&=-\varepsilon(p+1)\cdot (-1)^p\varepsilon^p p!\omega (e_2,e_i)=(-1)^{p+1}\varepsilon^{p+1}(p+1)!\omega (e_2,e_i)
\end{align*}
thanks to (\ref{eq::RfromGausswithAlpha0}).
Now by the induction principle (\ref{eq::Block3x3omega122i}) holds for all $p\geq 1$.
\end{proof}
\end{lem}


\begin{lem}\label{lm::Block3x3omega122312}
If $S_1$ is a 3-dimensional real Jordan block, $p\geq 1$ and $\alpha =0$ then
\begin{align}\label{eq::Block3x3omega123212}
R^p\omega(e_1,e_2,\ldots,e_1,e_2,e_2,e_3,e_1,e_2) =(-1)^{p+1}\varepsilon^p (p-1)!\omega (e_2,e_3).
\end{align}
\begin{proof}
First we notice that (\ref{eq::Block3x3omega123212}) is satisfied for $p=1$.
Indeed, since
$$
R(e_2,e_3)e_1=\varepsilon S_1e_2=\varepsilon e_3,\quad R(e_2,e_3)e_2=-\varepsilon S_1e_3=0
$$
we easily obtain
\begin{align*}
R\omega(e_2,e_3,e_1,e_2)=-\omega (R(e_2,e_3)e_1,e_2)-\omega (e_1,R(e_2,e_3)e_2)=\varepsilon\omega (e_2,e_3).
\end{align*}
Now assume that (\ref{eq::Block3x3omega123212}) holds for some $p\geq 1$, we shall show that it also holds for $p+1$. From  the Gauss equation and our assumptions we have
\begin{align*}
R(e_1,e_2)e_1=0, \quad R(e_1,e_2)e_2=\varepsilon S_1e_1=\varepsilon e_2, \quad R(e_1,e_2)e_3=-\varepsilon S_1e_2=-\varepsilon e_3.
\end{align*}
Now we obtain
\begin{align*}
R^{p+1}&\omega(\underbrace{e_1,e_2,\ldots,e_1,e_2}_{2p},e_2,e_3,e_1,e_2)\\
&=(R(e_1,e_2)\cdot R^{p}\omega)(\underbrace{e_1,e_2,\ldots,e_1,e_2}_{2p-2},e_2,e_3,e_1,e_2)\\
&=-\varepsilon(p+1)R^p\omega (e_1,e_2,\ldots,e_1,e_2,e_2,e_3,e_1,e_2)\\
&\phantom{=}-R^p\omega (e_1,e_2,\ldots,e_1,e_2,e_2,R(e_1,e_2)e_3,e_1,e_2)\\
&=-\varepsilon pR^p\omega (e_1,e_2,\ldots,e_1,e_2,e_2,e_3,e_1,e_2)\\
&=-\varepsilon p\cdot (-1)^{p+1}\varepsilon^p (p-1)!\omega (e_2,e_3)=(-1)^{p+2}\varepsilon^{p+1}p!\omega (e_2,e_3).
\end{align*}
Now by the induction principle (\ref{eq::Block3x3omega123212}) holds for all $p\geq 1$.
\end{proof}
\end{lem}


\begin{lem}\label{lm::s1s2real2}
Let us assume that $S_1$, $S_2$ (from Lemma \ref{lm::JordanDecompositionOfS}) are 2-dimensional real Jordan blocks. That is
$$
S_1=\left[
\begin{matrix}
\alpha & 0\\
1 & \alpha
\end{matrix}\right], \quad
S_2=\left[
\begin{matrix}
\beta & 0\\
1 & \beta
\end{matrix}\right],
$$
where $\alpha, \beta \in \R$.
We also assume that $H_1, H_2$ have the form
$$
H_1=\left[
\begin{matrix}
0 & \varepsilon\\
\varepsilon & 0
\end{matrix}\right], \quad
H_2=\left[
\begin{matrix}
0 & \eta\\
\eta & 0
\end{matrix}\right],
$$
where $\varepsilon,\eta\in\{-1,1\}$.
\par Then for every $i\in \{1,\ldots,2n\}\setminus \{2,4\}$ we have
\begin{align}
\label{eq::s1s2real2a}&R^p\omega (e_1,e_3,\ldots,e_1,e_3,e_i,e_3)=0 \quad \text{for }p\geq 1,
\end{align}
\begin{align}
\label{eq::s1s2real2b}R^{2p+1}\omega &(e_1,e_3,\ldots,e_1,e_3,e_i,e_4)\\
\nonumber &=(-1)^{p+1}\eta^{p+1}\varepsilon^p(\alpha \omega (e_i,e_1)+\omega (e_i,e_2)) \quad \text{for }p\geq 0.
\end{align}
\begin{proof}
First note that from the Gauss equation we easily obtain that
$R(e_1,e_3)e_i=0$ for $i\in \{1,\ldots,2n\}\setminus \{2,4\}$.
In particular
\begin{align}\label{eq::Re1e3withe1e3}
R(e_1,e_3)e_1=R(e_1,e_3)e_3=0.
\end{align}
Thus (\ref{eq::s1s2real2a}) follows immediately.
\par In order to prove (\ref{eq::s1s2real2b}) first notice that (\ref{eq::s1s2real2b}) is satisfied for $p=0$, since $R(e_1,e_3)e_i=0$ and
\begin{align}\label{eq::Re1e3withe4}
R(e_1,e_3)e_4=\eta S_1e_1=\eta\alpha e_1+\eta e_2.
\end{align}
Now assume that (\ref{eq::s1s2real2b}) holds for some $p\geq 0$, we shall show that it also holds for $p+1$. From the Gauss equation we also have
\begin{align}\label{eq::Re1e3withe2}
R(e_1,e_3)e_2=-\varepsilon S_2e_3=-\varepsilon\beta e_3-\varepsilon e_4.
\end{align}
Now using (\ref{eq::Re1e3withe1e3})--(\ref{eq::Re1e3withe2}) we obtain
\begin{align*}
R^{2p+3}&\omega(e_1,e_3,\ldots,e_1,e_3,e_i,e_4)\\
&=-\eta\alpha \underbrace{R^{2p+2}\omega(e_1,e_3,\ldots,e_1,e_3,e_i,e_1)}_{0}\\
&\phantom{=}-\eta R^{2p+2}\omega (e_1,e_3,\ldots,e_1,e_3,e_i,e_2)\\
&=\eta R^{2p+1}\omega (e_1,e_3,\ldots,e_1,e_3,e_i,R(e_1,e_3)e_2)\\
&=-\eta\varepsilon\beta R^{2p+1}\omega (e_1,e_3,\ldots,e_1,e_3,e_i,e_3)\\
&\phantom{=}-\eta\varepsilon R^{2p+1}\omega (e_1,e_3,\ldots,e_1,e_3,e_i,e_4)\\
&=-\eta\varepsilon R^{2p+1}\omega (e_1,e_3,\ldots,e_1,e_3,e_i,e_4),
\end{align*}
where the last equality follows from (\ref{eq::s1s2real2a}).
Finally we have
\begin{align*}
R^{2p+3}&\omega(e_1,e_3,\ldots,e_1,e_3,e_i,e_4)\\
&=-\eta\varepsilon R^{2p+1}\omega (e_1,e_3,\ldots,e_1,e_3,e_i,e_4)\\
&=-\eta\varepsilon(-1)^{p+1}\eta^{p+1}\varepsilon^{p}(\alpha \omega (e_i,e_1)+\omega (e_i,e_2))\\
&=(-1)^{p+2}\eta^{p+2}\varepsilon^{p+1}(\alpha \omega (e_i,e_1)+\omega (e_i,e_2)).
\end{align*}
By the induction principle (\ref{eq::s1s2real2b}) holds for all $p\geq 0$.
\end{proof}
\end{lem}


\begin{lem}\label{lm::rw12121414}
Let $S_1$, $S_2$ be 2-dimensional real Jordan blocks like in Lemma \ref{lm::s1s2real2}. If $\alpha = \beta =0$ then
\begin{align}
\label{eq::rw1212} R^{2p}\omega (e_1,e_3,\ldots,e_1,e_3,e_1,e_2,e_1,e_2)&=(-1)^p(\varepsilon\eta)^{p-1}\cdot 2^{2p-2}\omega (e_2,e_4),\\
\label{eq::rw1414} R^{2p}\omega (e_1,e_3,\ldots,e_1,e_3,e_1,e_4,e_1,e_4)&=(-1)^{p+1}(\varepsilon\eta)^{p}\cdot 2^{2p-2}\omega (e_2,e_4)
\end{align}
for every $p\geq 1$.
\begin{proof}
First note that formulas (\ref{eq::Re1e3withe1e3}), (\ref{eq::Re1e3withe4}) and (\ref{eq::Re1e3withe2}) from the proof of Lemma \ref{lm::s1s2real2} are still valid in our case.
Using them and taking into account that $\alpha = \beta =0$  we easily compute that
\begin{align*}
  R^{2}\omega(e_1,e_3,e_1,e_2,e_1,e_2)=-\varepsilon^2\omega (e_2,e_4)=-\omega (e_2,e_4)
\end{align*}
and
\begin{align*}
  R^{2}\omega(e_1,e_3,e_1,e_4,e_1,e_4)=\varepsilon\eta\omega (e_2,e_4).
\end{align*}
That is (\ref{eq::rw1212}) and (\ref{eq::rw1414}) are true for $p=1$.
Now assume that (\ref{eq::rw1212}) and (\ref{eq::rw1414}) hold for some $p\geq 1$.
Again using (\ref{eq::Re1e3withe1e3})--(\ref{eq::Re1e3withe2}) and the fact that $\alpha = \beta =0$ we obtain
\begin{align*}
  R^{2p+2}\omega &(e_1,e_3,\ldots,e_1,e_3,e_1,e_2,e_1,e_2)\\
  &=-R^{2p+1}\omega (e_1,e_3,\ldots,e_1,e_3,e_1,R(e_1,e_3)e_2,e_1,e_2)\\
  &\phantom{=}-R^{2p+1}\omega (e_1,e_3,\ldots,e_1,e_3,e_1,e_2,e_1,R(e_1,e_3)e_2)\\
  &=\varepsilon R^{2p+1}\omega (e_1,e_3,\ldots,e_1,e_3,e_1,e_4,e_1,e_2)\\
  &\phantom{=}+\varepsilon R^{2p+1}\omega (e_1,e_3,\ldots,e_1,e_3,e_1,e_2,e_1,e_4)\\
  &=-\varepsilon R^{2p}\omega (e_1,e_3,\ldots,e_1,e_3,e_1,R(e_1,e_3)e_4,e_1,e_2)\\
  &\phantom{=}-\varepsilon R^{2p}\omega (e_1,e_3,\ldots,e_1,e_3,e_1,e_4,e_1,R(e_1,e_3)e_2)\\
  &\phantom{=}-\varepsilon R^{2p}\omega (e_1,e_3,\ldots,e_1,e_3,e_1,R(e_1,e_3)e_2,e_1,e_4)\\
  &\phantom{=}-\varepsilon R^{2p}\omega (e_1,e_3,\ldots,e_1,e_3,e_1,e_2,e_1,R(e_1,e_3)e_4)\\
  &=-2\varepsilon\eta R^{2p}\omega (e_1,e_3,\ldots,e_1,e_3,e_1,e_2,e_1,e_2)\\
  &\phantom{=}+2\varepsilon^2 R^{2p}\omega (e_1,e_3,\ldots,e_1,e_3,e_1,e_4,e_1,e_4)\\
  &=2\varepsilon^2R^{2p}\omega (e_1,e_3,\ldots,e_1,e_3,e_1,e_4,e_1,e_4)\\
  &\phantom{=}-2\varepsilon\eta R^{2p}\omega(e_1,e_3,\ldots,e_1,e_3,e_1,e_2,e_1,e_2)\\
  &=2(\varepsilon^2 (-1)^{p+1}(\varepsilon\eta)^{p}\cdot 2^{2p-2}\omega (e_2,e_4)-\varepsilon\eta(-1)^{p}(\varepsilon\eta)^{p-1}\cdot 2^{2p-2}\omega (e_2,e_4))\\
  &=(-1)^{p+1}(\varepsilon\eta)^p\cdot 2^{2p}\omega (e_2,e_4).
\end{align*}
In a similar way we show that
\begin{align*}
  R^{2p+2}\omega &(e_1,e_3,\ldots,e_1,e_3,e_1,e_4,e_1,e_4)\\
  &=2\eta^2 R^{2p}\omega(e_1,e_3,\ldots,e_1,e_3,e_1,e_2,e_1,e_2)\\
  &\phantom{=}-2\varepsilon\eta R^{2p}\omega (e_1,e_3,\ldots,e_1,e_3,e_1,e_4,e_1,e_4)\\
  &=2\eta^2 (-1)^p(\varepsilon\eta)^{p-1}\cdot 2^{2p-2}\omega (e_2,e_4)\\
  &\phantom{=}-2\varepsilon\eta (-1)^{p+1}(\varepsilon\eta)^{p}\cdot 2^{2p-2}\omega (e_2,e_4)\\
  &=(-1)^{p+2}2^{2p}(\varepsilon\eta)^{p+1}\omega(e_2,e_4)
\end{align*}
Now by the induction principle (\ref{eq::rw1212}) and (\ref{eq::rw1414}) hold for all $p\geq 1$.
\end{proof}
\end{lem}

\begin{thm}\label{tw::MaxOneRealOfDim2}
Let $f\colon \M\rightarrow\R^{2n+1}$ ($\dim M\geq 4$)  be a non-degenerate affine hypersurface with a locally equiaffine transversal vector field  $\xi$
and an almost symplectic form $\omega$. Let
\begin{align}\label{eq::SDecomposition}
S=\left[\begin{matrix}
S_1 & 0 & \ldots & 0 \\
0 & S_2 & \ldots & 0 \\
\vdots & \vdots & \ddots & \ldots \\
0 & 0 & \ldots & S_{q+r}
\end{matrix}\right]
\end{align}
be the Jordan decomposition of $S$ as stated in the Lemma \ref{lm::JordanDecompositionOfS}.
If $R^p\omega=0$ for some $p\geq 1$ then $\dim S_1\leq 2$ and $\dim S_i=1$ for $i=2,\ldots,q$.
\begin{proof}
By Lemma \ref{lm::S1dimMax3} $\dim S_1\leq 3$. If $\dim S_1=3$, using Corollary
\ref{cor::rpomegaeiek} we obtain
\begin{align}\label{eq::e3eiZero}
\varepsilon^p\alpha^p\omega (e_2,e_3)=\ldots =\varepsilon^p\alpha^p\omega (e_{2n},e_3)=0.
\end{align}
By Lemma \ref{lm::Block3x3omega12} $\omega (e_1,e_2)=0$. Now using Lemma \ref{lm::Block3x3omega12ij}
we obtain
$$
h(e_i,e_j)\omega(e_1,e_3)=0
$$
for $i,j>3$. Since $\dim M\geq 4$ and $h$ is non-degenerate there exist $i,j>3$ such that $h(e_i,e_j)\neq 0$ and in consequence $\omega(e_1,e_3)=0$.
If $\alpha\neq 0$ then from (\ref{eq::e3eiZero}) we get $\omega (e_2,e_3)=\ldots =\omega (e_{2n},e_3)=0$, so $\omega$
is degenerate. That is we must have $\alpha=0$. Now using Lemma \ref{lm::Block3x3omega122i} and Lemma \ref{lm::Block3x3omega122312}
we show that $\omega (e_2,e_i)=0$ for $i=1,\ldots,2n$ that is $\omega$ is degenerate again. In consequence the case $\dim S_1=3$ is not possible.
\par Assume now that $\dim S_1=2$. If $\dim S_2=2$ then
$$
S_1=\left[
\begin{matrix}
\alpha & 0\\
1 & \alpha
\end{matrix}\right], \quad
S_2=\left[
\begin{matrix}
\beta & 0\\
1 & \beta
\end{matrix}\right],
$$ where $\alpha, \beta \in \R$. If $R^p\omega=0$ for some $p\geq 1$ then also $R^{2p+1}\omega=0$. Now using Lemma \ref{lm::s1s2real2}
we get
$$
\alpha \omega (e_i,e_1)+\omega (e_i,e_2)=0
$$
for $i\in \{1,\ldots,2n\}\setminus \{2,4\}$. In particular, (for $i=1$) we get that $\omega (e_1,e_2)=0$.
On the other hand by Corollary \ref{cor::rpomegaeiek} we have
$$
\varepsilon^p\alpha^p\omega (e_3,e_2)=\ldots =\varepsilon^p\alpha^p\omega (e_{2n},e_2)=0.
$$
Note that case $\alpha\neq 0$ is not possible since then $\omega (e_i,e_2)=0$ for all $i\in\{1,\ldots,2n\}$ and $\omega$ is degenerate.
Thus we must have $\alpha=0$. In this case Lemma \ref{lm::s1s2real2} implies that
$\omega (e_i,e_2)=0$
for $i\in \{1,\ldots,2n\}\setminus \{2,4\}$.
Since, without loss of generality, we can exchange $S_1$ with $S_2$ we also get that $\beta=0$. Now by Lemma \ref{lm::rw12121414}
we get that $\omega (e_2,e_4)=0$ and $\omega$ is degenerate, so this case is also not possible. Summarising we must have
$\dim S_1\leq 2$ and $\dim S_i=1$ for $i=2,\ldots,q$, what completes the proof.
\end{proof}
\end{thm}


\section{Complex Jordan blocks}
In this chapter we study properties of complex Jordan blocks of the shape operator $S$.
Before we proceed, to simplify proofs in this chapter, we need to do slight modification in the notation of Lemma \ref{lm::JordanDecompositionOfS}.
\par Let $\{e_1,\ldots,e_{2n}\}$ be the basis of $T_xM$ from Lemma \ref{lm::JordanDecompositionOfS}.
Without loss of generality rearranging and renaming vectors $e_1,\ldots,e_{2n}$ we can
change order of $S_i$ and $H_i$ in such way that $S_{1},\ldots,S_{r}$ will be complex blocks and $S_{r+1},\ldots,S_{r+q}$ will be
real blocks. If we assume that $S_1$ is a $2k$-dimensional block, $k\geq 1$ in the new notation we will have
\begin{align*}
S_1=\left[\begin{matrix}
\alpha & \beta  &   0    & 0       & \ldots & 0      & 0 \\
-\beta & \alpha &   0    & 0       & \ldots & 0      & 0\\
1      & 0      & \alpha &\beta    & \ldots & 0      & 0\\
0      & 1      & -\beta & \alpha  & \ldots & 0      & 0\\
\vdots & \vdots & \vdots & \vdots  & \ddots & \ldots & \ldots\\
0      & 0      & 0      & 0       & \ldots & \alpha & \beta\\
0      & 0      & 0      & 0       & \ldots & -\beta &\alpha
\end{matrix}\right]\in M(2k,2k,\R),
\end{align*}
where $\alpha, \beta \in \R$, $\beta\neq 0$
and
\begin{align*}
H_1=\left[\begin{matrix}
0 & 0 & \cdots & 0 & 1 \\
0 & 0 & \cdots & 1 & 0 \\
\vdots & \vdots & \ddots & \cdots & \cdots \\
0 & 1 & \cdots & 0 & 0 \\
1 & 0 & \cdots & 0 & 0
\end{matrix}\right]\in M(2k,2k,\R).
\end{align*}
Moreover, vectors $\{e_1,\ldots,e_{2k}\}$ will be a basis for $S_1$.
\par In all the below lemmas (if not stated otherwise) we always assume that
$S_1$ and $H_1$ are as above.

\par Let us start with the following three lemmas related to 2-dimensional complex Jordan blocks.
\begin{lem}\label{lm::ComplexBlocks2x2BasicProperties}
If $S_1$ is a 2-dimensional complex Jordan block then for every $i\in\{3,\ldots,2n\}$ we have
\begin{align*}
R(e_1,e_2)e_1&=Se_1=\alpha e_1 -\beta e_2,\\
R(e_1,e_2)e_2&=-Se_2=-\beta e_1-\alpha e_2,\\
R(e_1,e_2)e_i&=0.
\end{align*}
\begin{proof}
Proof easily follows from the Gauss equation and the fact that $h(e_1,e_i)=h(e_2,e_i)=0$ for $i>2$.
\end{proof}
\end{lem}

\begin{lem}\label{lm::ComplexBlocks1}
If $S_1$ is a 2-dimensional complex Jordan block then for every $p\geq 1$, $i\in\{3,\ldots,2n\}$ we have
\begin{align}
\label{eq::Complex::R2p} R^{2p}\omega (e_1,e_{2},e_1,e_{2},\ldots,e_1,e_{2},e_1,e_i)=(\det S_1)^p\omega(e_1,e_i)\\
\label{eq::ComplexR2pe2}R^{2p}\omega (e_1,e_{2},e_1,e_{2},\ldots,e_1,e_{2},e_2,e_i)=(\det S_1)^p\omega(e_2,e_i)
\end{align}
\begin{proof}
We shall prove (\ref{eq::Complex::R2p}). For $p=1$,
using Lemma \ref{lm::ComplexBlocks2x2BasicProperties} we compute
\begin{align*}
R^2\omega(e_1,e_2,&e_1,e_2,e_1,e_i)\\
&=-R\omega(R(e_1,e_2)e_1,e_2,e_1,e_i)-R\omega (e_1,R(e_1,e_2)e_2,e_1,e_i)\\
&\phantom{=}-R\omega(e_1,e_2,R(e_1,e_2)e_1,e_i)-R\omega(e_1,e_2,e_1,R(e_1,e_2)e_i)\\
&=-\alpha R\omega(e_1,e_2,e_1,e_i)+\alpha R\omega(e_1,e_2,e_1,e_i)\\
&\phantom{=}-R\omega(e_1,e_2,\alpha e_1-\beta e_2,e_i)\\
&=-\alpha R\omega(e_1,e_2,e_1,e_i)+\beta R\omega(e_1,e_2,e_2,e_i)\\
&=-\alpha (-\omega(R(e_1,e_2)e_1,e_i)-\omega(e_1,R(e_1,e_2)e_i))\\
&\phantom{=}+\beta (-\omega(R(e_1,e_2)e_2,e_i)-\omega(e_2,R(e_1,e_2)e_i))\\
&=\alpha ^2\omega(e_1,e_i)-\alpha \beta\omega(e_2,e_i)+\beta^2\omega(e_1,e_i)+\beta\alpha\omega(e_2,e_i)\\
&=(\alpha^2+\beta^2)\omega(e_1,e_i).
\end{align*}
Assume that (\ref{eq::Complex::R2p}) holds for some $p\geq 1$ we shall show that it also holds for $p+1$.
Indeed, using Lemma \ref{lm::ComplexBlocks2x2BasicProperties} we have
\begin{align*}
R^{2p+2}&\omega(e_1,e_2,\ldots,e_1,e_2,e_1,e_i)\\
&=-R^{2p+1}\omega(R(e_1,e_2)e_1,e_2,\ldots,e_1,e_2,e_1,e_i)\\
&\phantom{=}-R^{2p+1}\omega(e_1,R(e_1,e_2)e_2,\ldots,e_1,e_2,e_1,e_i)\\
&\phantom{=}-\ldots \\
&\phantom{=}-R^{2p+1}\omega(e_1,e_2,\ldots,R(e_1,e_2)e_1,e_2,e_1,e_i)\\
&\phantom{=}-R^{2p+1}\omega(e_1,e_2,\ldots,e_1,R(e_1,e_2)e_2,e_1,e_i)\\
&\phantom{=}-R^{2p+1}\omega(e_1,e_2,\ldots,e_1,e_2,R(e_1,e_2)e_1,e_i)\\
&\phantom{=}-R^{2p+1}\omega(e_1,e_2,\ldots,e_1,e_2,e_1,R(e_1,e_2)e_i)\\
&=- \alpha R^{2p+1}\omega(e_1,e_2,\ldots,e_1,e_2,e_1,e_i)+\alpha R^{2p+1}\omega(e_1,e_2,\ldots,e_1,e_2,e_1,e_i)\\
&\phantom{=}-\ldots\\
&\phantom{=}-\alpha R^{2p+1}\omega(e_1,e_2,\ldots,e_1,e_2,e_1,e_i)+\alpha R^{2p+1}\omega(e_1,e_2,\ldots,e_1,e_2,e_1,e_i)\\
&\phantom{=}-R^{2p+1}\omega(e_1,e_2,\ldots,e_1,e_2,\alpha e_1 -\beta e_2,e_i)-0\\
&=-R^{2p+1}\omega(e_1,e_2,\ldots,e_1,e_2,\alpha e_1 -\beta e_2,e_i)\\
&=-\alpha R^{2p+1}\omega(e_1,e_2,\ldots,e_1,e_2,e_1,e_i)+\beta R^{2p+1}\omega(e_1,e_2,\ldots,e_1,e_2,e_2,e_i).
\end{align*}
Now we compute that
\begin{align*}
R^{2p+1}&\omega(e_1,e_2,\ldots,e_1,e_2,e_1,e_i)\\
&=-R^{2p}\omega(R(e_1,e_2)e_1,e_2,\ldots,e_1,e_2,e_1,e_i)\\
&\phantom{=}-R^{2p}\omega(e_1,R(e_1,e_2)e_2,\ldots,e_1,e_2,e_1,e_i)\\
&\phantom{=}\ldots\\
&\phantom{=}-R^{2p}\omega(e_1,e_2,\ldots,e_1,e_2,R(e_1,e_2)e_1,e_i)\\
&\phantom{=}-R^{2p}\omega(e_1,e_2,\ldots,e_1,e_2,e_1,R(e_1,e_2)e_i)\\
&=-R^{2p}\omega(e_1,e_2,\ldots,e_1,e_2,R(e_1,e_2)e_1,e_i).
\end{align*}
Similarly
\begin{align*}
R^{2p+1}&\omega(e_1,e_2,\ldots,e_1,e_2,e_2,e_i)=-R^{2p}\omega(e_1,e_2,\ldots,e_1,e_2,R(e_1,e_2)e_2,e_i).
\end{align*}
Thus we obtain
\begin{align*}
R^{2p+2}&\omega(e_1,e_2,\ldots,e_1,e_2,e_1,e_i)\\
&=\alpha R^{2p}\omega(e_1,e_2,\ldots,e_1,e_2,\alpha e_1-\beta e_2,e_i)\\
&\phantom{=}-\beta R^{2p}\omega(e_1,e_2,\ldots,e_1,e_2,-\beta e_1-\alpha e_2,e_i)\\
&=\alpha^2R^{2p}\omega(e_1,e_2,\ldots,e_1,e_2,e_1,e_i)+\beta^2 R^{2p}\omega(e_1,e_2,\ldots,e_1,e_2,e_1,e_i)\\
&=\det S_1\cdot R^{2p}\omega(e_1,e_2,\ldots,e_1,e_2,e_1,e_i)=(\det S_1)^{p+1}\omega(e_1,e_i).
\end{align*}
Now by the induction principle (\ref{eq::Complex::R2p}) holds for all $p\geq 1$.
The formula (\ref{eq::ComplexR2pe2}) can be shown in a similar way.
\end{proof}
\end{lem}


\begin{lem}\label{lm::OneComplexAndOther}
If $S_1$ is a 2-dimensional complex Jordan block then for every $p\geq 1$, $i,j >2$ we have
\begin{align}
\label{eq::ComplexAndOther::R2p} R^{2p}\omega (e_1,e_{2},\ldots,e_1,e_{2},e_1,e_{i},e_2,e_j)&=2^{2p-1}\beta^2(\det S_1)^{p-1}h(e_i,e_j)\omega(e_1,e_2)\\
\label{eq::ComplexAndOther::R2p2} R^{2p}\omega (e_1,e_{2},\ldots,e_1,e_{2},e_2,e_{i},e_1,e_j)&=2^{2p-1}\beta^2(\det S_1)^{p-1}h(e_i,e_j)\omega(e_1,e_2)\\
\label{eq::ComplexAndOther::R2p3} R^{2p}\omega (e_1,e_{2},\ldots,e_1,e_{2},e_1,e_{i},e_1,e_j)&-R^{2p}\omega (e_1,e_{2},\ldots,e_1,e_{2},e_2,e_{i},e_2,e_j)\\
\nonumber &=-2^{2p}\alpha\beta(\det S_1)^{p-1}h(e_i,e_j)\omega(e_1,e_2)
\end{align}
\begin{proof}
To simplify computations let us denote
\begin{align}\label{eq::Brsijp}
B^p(X,Y,i,j):=R^p\omega(e_1,e_{2},e_1,e_{2},\ldots,e_1,e_{2},X,e_{i},Y,e_{j})
\end{align}
where $X,Y\in\operatorname{span}\{e_1,e_2\}$, $i,j\in\{3,\dots,2n\}$ and $p\geq 1$.
\par First we shall show that (\ref{eq::ComplexAndOther::R2p})--(\ref{eq::ComplexAndOther::R2p3}) are true for $p=1$.
Indeed, using Lemma \ref{lm::ComplexBlocks2x2BasicProperties}, we obtain
\begin{align*}
R^{2}&\omega (e_1,e_{2},e_1,e_{i},e_2,e_j)\\
&=-R\omega (S e_1,e_{i},e_2,e_j)+R\omega (e_1,e_{i},S e_2,e_j)\\
&=-\alpha R\omega (e_1,e_{i},e_2,e_j)+\beta R\omega (e_2,e_{i},e_2,e_j)\\
&\phantom{=}+\beta R\omega (e_1,e_{i},e_1,e_j)+\alpha R\omega (e_1,e_{i},e_2,e_j)\\
&=\beta (R\omega (e_2,e_{i},e_2,e_j)+R\omega (e_1,e_{i},e_1,e_j))\\
&=\beta (-\omega(R(e_2,e_{i})e_2,e_j)-\omega(e_2,R(e_2,e_{i})e_j)\\
&\phantom{=}-\omega(R(e_1,e_{i})e_1,e_j)-\omega(e_1,R(e_1,e_{i})e_j))\\
&=-\beta (\omega(e_2,R(e_2,e_{i})e_j)+\omega(e_1,R(e_1,e_{i})e_j))\\
&=-\beta h(e_i,e_j)\Big(\omega(e_2,S e_2)+\omega(e_1, S e_1)\Big)\\
&=2\beta^2 h(e_i,e_j)\omega(e_1,e_2).
\end{align*}
that is
\begin{align}
B^2(e_1,e_2,i,j)=2\beta^2 h(e_i,e_j)\omega(e_1,e_2).
\end{align}
Exactly in the same way we show that
\begin{align}
B^2(e_2,e_1,i,j)=2\beta^2 h(e_i,e_j)\omega(e_1,e_2).
\end{align}
For (\ref{eq::ComplexAndOther::R2p3}) we compute
\begin{align*}
B^2(e_1,e_1,i,j)-B^2(e_2,e_2,i,j)&=-B^1(R(e_1,e_2)e_1,e_1,i,j)-B^1(e_1,R(e_1,e_2)e_1,i,j)\\
&\phantom{=}+B^1(R(e_1,e_2)e_2,e_2,i,j)+B^1(e_2,R(e_1,e_2)e_2,i,j)\\
&=-B^1(S e_1,e_1,i,j)-B^1(e_1,S e_1,i,j)\\
&\phantom{=}-B^1(S e_2,e_2,i,j)-B^1(e_2,S e_2,i,j)\\
&=-2\alpha(B^1(e_1,e_1,i,j)+B^1(e_2,e_2,i,j))\\
&=-2\alpha(-h(e_i,e_j)\omega(e_1, S e_1)-h(e_i,e_j)\omega(e_2, S e_2))\\
&=2\alpha h(e_i,e_j)(-\beta\omega(e_1,e_2)+\beta\omega(e_2,e_1))\\
&=-4\alpha\beta h(e_i,e_j)\omega(e_1,e_2).
\end{align*}
\par Now assume that (\ref{eq::ComplexAndOther::R2p})--(\ref{eq::ComplexAndOther::R2p3}) are all true for some $p\geq 1$.
We compute
\begin{align*}
B^{2p+2}(e_1,e_2,i,j)&=R^{2p+2}\omega (e_1,e_{2},\ldots,e_1,e_{2},e_1,e_{i},e_2,e_j)\\
&=-R^{2p+1}\omega (S e_1,e_{2},\ldots,e_1,e_{2},e_1,e_{i},e_2,e_j)\\
&\phantom{=}+R^{2p+1}\omega (e_1,S e_{2},\ldots,e_1,e_{2},e_1,e_{i},e_2,e_j)\\
&\phantom{=}\cdots\\
&\phantom{=}-R^{2p+1}\omega (e_1,e_{2},\ldots,S e_1,e_{2},e_1,e_{i},e_2,e_j)\\
&\phantom{=}+R^{2p+1}\omega (e_1,e_{2},\ldots,e_1,S e_{2},e_1,e_{i},e_2,e_j)\\
&\phantom{=}-R^{2p+1}\omega (e_1,e_{2},\ldots,e_1,e_{2},S e_1,e_{i},e_2,e_j)\\
&\phantom{=}+R^{2p+1}\omega (e_1,e_{2},\ldots,e_1,e_{2},e_1,e_{i},S e_2,e_j)\\
&=-R^{2p+1}\omega (e_1,e_{2},\ldots,e_1,e_{2},S e_1,e_{i},e_2,e_j)\\
&\phantom{=}+R^{2p+1}\omega (e_1,e_{2},\ldots,e_1,e_{2},e_1,e_{i},S e_2,e_j)\\
&=-B^{2p+1}(S e_1,e_2,i,j)+B^{2p+1}(e_1,S e_2,i,j),
\end{align*}
since for $l=1,\ldots,2p$, terms $2l-1$ and $2l$ cancel each other.
Now we easily compute that
\begin{align*}
B^{2p+1}(S e_1,e_2,i,j)=-B^{2p}(R(e_1,e_2)S e_1,e_2,i,j)+B^{2p}(S e_1,S e_2,i,j)
\end{align*}
and
\begin{align*}
B^{2p+1}(e_1,S e_2,i,j)=-B^{2p}(S e_1,S e_2,i,j)-B^{2p}(e_1,R(e_1,e_2)S e_2,i,j)
\end{align*}
In consequence
\begin{align*}
B^{2p+2}(e_1,e_2,i,j)&=-2 B^{2p}(S e_1,S e_2,i,j)+B^{2p}(R(e_1,e_2)S e_1,e_2,i,j)\\
&\phantom{=}-B^{2p}(e_1,R(e_1,e_2)S e_2,i,j)\\
&=-2\alpha\beta B^{2p}(e_1,e_1,i,j)+2\alpha\beta B^{2p}(e_2,e_2,i,j)\\
&\phantom{=}-2\alpha^2 B^{2p}(e_1,e_2,i,j)+2\beta^2 B^{2p}(e_2,e_1,i,j)\\
&\phantom{=}+(\alpha^2+\beta^2)B^{2p}(e_1,e_2,i,j)+(\alpha^2+\beta^2)B^{2p}(e_1,e_2,i,j)\\
&=-2\alpha\beta (B^{2p}(e_1,e_1,i,j)-B^{2p}(e_2,e_2,i,j))\\
&\phantom{=}+2\beta^2(B^{2p}(e_1,e_2,i,j)+B^{2p}(e_2,e_1,i,j)),
\end{align*}
since $R(e_1,e_2)S e_1=(\alpha^2+\beta^2)e_1$ and $R(e_1,e_2)S e_2=-(\alpha^2+\beta^2)e_2$.
Now using assumptions (\ref{eq::ComplexAndOther::R2p})--(\ref{eq::ComplexAndOther::R2p3}) we obtain
\begin{align*}
B^{2p+2}(e_1,e_2,i,j)&=-2\alpha\beta (-2^{2p}\alpha\beta(\det S_1)^{p-1}h(e_i,e_j)\omega(e_1,e_2))\\
&\phantom{=}+4\beta^2(2^{2p-1}\beta^2(\det S_1)^{p-1}h(e_i,e_j)\omega(e_1,e_2))\\
&=2^{2p+1}(\alpha^2\beta^2+\beta^4)(\det S_1)^{p-1}h(e_i,e_j)\omega(e_1,e_2)\\
&=2^{2p+1}\beta^2(\det S_1)^{p}h(e_i,e_j)\omega(e_1,e_2)
\end{align*}
In a similar way we show that
\begin{align*}
B^{2p+2}(e_2,e_1,i,j)=2^{2p+1}\beta^2(\det S_1)^{p}h(e_i,e_j)\omega(e_1,e_2).
\end{align*}
Eventually
\begin{align*}
B^{2p+2}(e_1,e_1,i,j)&-B^{2p+2}(e_2,e_2,i,j)\\
&=-B^{2p+1}(R(e_1,e_2)e_1,e_1,i,j)-B^{2p+1}(e_1,R(e_1,e_2)e_1,i,j)\\
&\phantom{=}+B^{2p+1}(R(e_1,e_2)e_2,e_2,i,j)+B^{2p+1}(e_2,R(e_1,e_2)e_2,i,j)\\
&=-B^{2p+1}(S e_1,e_1,i,j)-B^{2p+1}(e_1,S e_1,i,j)\\
&\phantom{=}-B^{2p+1}(S e_2,e_2,i,j)-B^{2p+1}(e_2,S e_2,i,j)\\
&=-2\alpha(B^{2p+1}(e_1,e_1,i,j)+B^{2p+1}(e_2,e_2,i,j))\\
&=-2\alpha(-B^{2p}(S e_1,e_1,i,j)-B^{2p}(e_1,S e_1,i,j)\\
&\phantom{=}+B^{2p}(S e_2,e_2,i,j)+B^{2p}(e_2,S e_2,i,j))\\
&=-2\alpha(-2\alpha(B^{2p}(e_1,e_1,i,j)-B^{2p}(e_2,e_2,i,j))+4\beta B^{2p}(e_1,e_2,i,j))\\
&=4\alpha^2(B^{2p}(e_1,e_1,i,j)-B^{2p}(e_2,e_2,i,j))-8\alpha\beta B^{2p}(e_1,e_2,i,j).
\end{align*}
Using (\ref{eq::ComplexAndOther::R2p})--(\ref{eq::ComplexAndOther::R2p3}) we get
\begin{align*}
B^{2p+2}(e_1,e_1,i,j)&-B^{2p+2}(e_2,e_2,i,j)\\
&=4\alpha^2(B^{2p}(e_1,e_1,i,j)-B^{2p}(e_2,e_2,i,j))-8\alpha\beta B^{2p}(e_1,e_2,i,j)\\
&=4\alpha^2(-2^{2p}\alpha\beta(\det S_1)^{p-1}h(e_i,e_j)\omega(e_1,e_2))\\
&\phantom{=}-8\alpha\beta\cdot 2^{2p-1}\beta^2(\det S_1)^{p-1}h(e_i,e_j)\omega(e_1,e_2)\\
&=-2^{2p+2}(\alpha^3\beta+\alpha\beta^3)(\det S_1)^{p-1}h(e_i,e_j)\omega(e_1,e_2)\\
&=-2^{2p+2}\alpha\beta(\det S_1)^{p}h(e_i,e_j)\omega(e_1,e_2).
\end{align*}
Now by the induction principle (\ref{eq::ComplexAndOther::R2p})--(\ref{eq::ComplexAndOther::R2p3}) hold for all $p\geq 1$.
\end{proof}
\end{lem}

In the next three lemmas we study properties of complex Jordan block of dimension greater than 2 in relation to other Jordan blocks from the decomposition.
Thus in these lemmas, we implicitly assume that the Jordan decomposition contains more than one (not necessarily complex) block.

\begin{lem}\label{lm::Rpomegac1}
Let $S_1$ be a $2k$-dimensional complex Jordan block, $k\geq 2$, $i\in \{1,\ldots,2k-1\}\setminus \{3\}$, $s\in \{1,\ldots,2k\}$, $j>2k$, $p\geq 1$.
Then
\begin{align}\label{eq::lm::Rpomegac1}
R^p\omega(e_1&,e_{2k-2},\ldots,e_1,e_{2k-2},e_i,e_s,e_1,e_j)\\
\nonumber&=
\begin{cases}
0,                  &\quad \text{for $s<2k$}\\
-\omega(S_1e_i,e_j) &\quad \text{for $s=2k$.}
\end{cases}
\end{align}
\begin{proof}
First let us  notice that
\begin{align*}
R(e_i,e_s)e_1=
\begin{cases}
0, &\quad \text{for}\quad s<2k\\
S_1e_i, &\quad \text{for}\quad s=2k
\end{cases}
\end{align*}
and
\begin{align*}
R(e_i,e_s)e_j=0.
\end{align*}
From the above properties we obtain
\begin{align*}
R\omega(e_i,e_s,e_1,e_j)&=-\omega (R(e_i,e_s)e_1,e_j)-\omega(e_1,R(e_i,e_s)e_j)\\
&=\begin{cases}
0,                   &\quad \text{for $s<2k$}\\
-\omega(S_1e_i,e_j), &\quad \text{for $s=2k$}
  \end{cases}
\end{align*}
thus (\ref{eq::lm::Rpomegac1}) is true for $p=1$.
Moreover we have that
\begin{align*}
R(e_1,e_{2k-2})e_s&=
\begin{cases}
0,            &\quad \text{for $s\in\{1,\ldots,2k-1\}\setminus \{3\}$} \\
S_1e_1,       &\quad \text{for $s=3$} \\
-S_1e_{2k-2}, &\quad \text{for $s=2k$,}
\end{cases}\\
\end{align*}
and
\begin{align*}
R(e_1,e_{2k-2})e_1=R(e_1,e_{2k-2})e_{2k-2}=R(e_1,e_{2k-2})e_i=R(e_1,e_{2k-2})e_j=0.
\end{align*}
Now, assume that (\ref{eq::lm::Rpomegac1}) is true for some $p\geq1$. Using the above formulas we easily get
\begin{align*}
&R^{p+1}\omega(e_1,e_{2k-2},\ldots,e_1,e_{2k-2},e_i,e_s,e_1,e_j)\\
&=(R(e_1,e_{2k-2})\cdot R^p\omega)(e_1,e_{2k-2},\ldots,e_1,e_{2k-2},e_i,e_s,e_1,e_j)\\
&=-R^p\omega(e_1,e_{2k-2},\ldots,e_1,e_{2k-2},e_i,R(e_1,e_{2k-2})e_s,e_1,e_j)\\
&=
\begin{cases}
0,                                                                    & \hspace{-3cm} \text{for }s\in\{1,\ldots,2k-1\}\setminus \{3\}\\
-R^p\omega(e_1,e_{2k-2},\ldots,e_1,e_{2k-2},e_i,S_1e_1,e_1,e_j),      & \text{for } s=3\\
R^p\omega(e_1,e_{2k-2},\ldots,e_1,e_{2k-2},e_i,S_1e_{2k-2},e_1,e_j),  & \text{for } s=2k
\end{cases}\\
&=
\begin{cases}
0, &\hspace{-3cm} \text{for }s\in\{1,\ldots,2k-1\}\setminus \{3\}\\
-R^p\omega(e_1,e_{2k-2},\ldots,e_1,e_{2k-2},e_i,\alpha e_1-\beta e_2+e_3,e_1,e_j), &\text{for }s=3\\
R^p\omega(e_1,e_{2k-2},\ldots,e_1,e_{2k-2},e_i,\beta e_{2k-3}+\alpha e_{2k-2}+e_{2k},e_1,e_j), &\text{for } s=2k
\end{cases}\\
&=
\begin{cases}
0, &\hspace{-3cm} \text{for } s\in\{1,\ldots,2k-1\}\setminus \{3\}\\
0,&\text{for } s=3\\
R^p\omega(e_1,e_{2k-2},\ldots,e_1,e_{2k-2},e_i,e_{2k},e_1,e_j),&\text{for }s=2k.
\end{cases}\\
&=
\begin{cases}
0, &\text{for }s\in\{1,\ldots,2k-1\}\setminus \{3\}\\
0, &\text{for }s=3\\
-\omega(S_1e_i,e_j), &\text{for } s=2k.
\end{cases}
\end{align*}
Now by the induction principle (\ref{eq::lm::Rpomegac1}) holds for all $p\geq 1$.
\end{proof}
\end{lem}

\begin{lem}\label{lm::Rpomegac2}
Let $S_1$ be a $2k$-dimensional complex Jordan block, $k\geq 2$, $j>2k$, $p\geq 1$.
Then
\begin{align}\label{eq::lm::Rpomegac2}
R^p\omega(e_1&,e_{2k-1},\ldots,e_1,e_{2k-1},e_3,e_{2k},e_1,e_j)=-(-\beta )^{p-1}\omega(S_1e_3,e_j).
\end{align}
\begin{proof}
For $p=1$, by direct computations, we obtain
\begin{align*}
R\omega(e_3,e_{2k},e_1,e_j)&=-\omega (R(e_3,e_{2k})e_1,e_j)-\omega(e_1,R(e_3,e_{2k})e_j)\\
&=-\omega(S_1e_3,e_j).
\end{align*}
Now let us assume that (\ref{eq::lm::Rpomegac2}) holds for some $p\geq 1$.
Since
\begin{align*}
R(e_1,e_{2k-1})e_1&=R(e_1,e_{2k-1})e_3=R(e_1,e_{2k-1})e_{2k-1}=R(e_1,e_{2k-1})e_j=0,\\
R(e_1,e_{2k-1})e_{2k}&=-S_1e_{2k-1}=-\alpha e_{2k-1}+\beta e_{2k}
\end{align*}
we compute that
\begin{align*}
(R^{p+1}\omega)(e_1&,e_{2k-1},\ldots,e_1,e_{2k-1},e_3,e_{2k},e_1,e_j)\\
&=(R(e_1,e_{2k-1})\cdot R^p\omega)(e_1,e_{2k-1},\ldots,e_1,e_{2k-1},e_3,e_{2k},e_1,e_j)\\
&=-R^p\omega(e_1,e_{2k-1},\ldots,e_1,e_{2k-1},e_3,-\alpha e_{2k-1}+\beta e_{2k},e_1,e_j)\\
&=\alpha R^p\omega(e_1,e_{2k-1},\ldots,e_1,e_{2k-1},e_3,e_{2k-1},e_1,e_j)\\
&\phantom{=}-\beta R^p\omega(e_1,e_{2k-1},\ldots,e_1,e_{2k-1},e_3,e_{2k},e_1,e_j).
\end{align*}
Since
$$
R(e_3,e_{2k-1})e_1=R(e_3,e_{2k-1})e_j=0
$$
one may easily deduce that $R\omega(e_3,e_{2k-1},e_1,e_j)=0$ and more general
$$
R^p\omega(e_1,e_{2k-1},\ldots,e_1,e_{2k-1},e_3,e_{2k-1},e_1,e_j)=0
$$
for all $p\geq 1$.
Thus we have
\begin{align*}
(R^{p+1}\omega)(e_1&,e_{2k-1},\ldots,e_1,e_{2k-1},e_3,e_{2k},e_1,e_j)\\
&=-\beta R^p\omega(e_1,e_{2k-1},\ldots,e_1,e_{2k-1},e_3,e_{2k},e_1,e_j)\\
&=-\beta\cdot (-(-\beta)^{p-1})\omega (S_1e_3,e_j)=-(-\beta)^p\omega(S_1e_3,e_j),
\end{align*}
where the last equality follows from (\ref{eq::lm::Rpomegac2}).
Now by the induction principle (\ref{eq::lm::Rpomegac2}) holds for all $p\geq 1$.
\end{proof}
\end{lem}

\begin{lem}\label{lm::Rpomegac3}
Let $S_1$ be a $2k$-dimensional complex Jordan block, $k\geq 2$, $i\in \{1,\ldots,2k\}$, $j>2k$.
Then for every $p\geq 1$ we have
\begin{align}\label{eq::lm::Rpomegac3}
R^p\omega(e_2&,e_{2k},\ldots,e_2,e_{2k},e_i,e_{2k},e_{2k},e_j)\\
&=
\begin{cases}
0, \quad \text{for}\quad i\in \{2,\ldots,2k\}\\ \nonumber
(-\beta)^{p-1}\omega(S_1e_{2k},e_j) \quad \text{for}\quad i=1.
\end{cases}
\end{align}
\begin{proof}
First notice that
\begin{align*}
R(e_i,e_{2k})e_{2k}&=
\begin{cases}
0, \quad \text{for}\quad i\in \{2,\ldots,2k\}\\
-S_1e_{2k}, \quad \text{for}\quad i=1
\end{cases}\\
R(e_i,e_{2k})e_j&=0.
\end{align*}
From the above we obtain
\begin{align*}
R\omega(e_i,e_{2k},e_{2k},e_j)&=-\omega (R(e_i,e_{2k})e_{2k},e_j)-\omega(e_{2k},R(e_i,e_{2k})e_j)\\
&=\begin{cases}
0, \quad \text{for}\quad i\in \{2,\ldots,2k\}\\
\omega(S_1e_{2k},e_j), \quad \text{for}\quad i=1,
  \end{cases}
\end{align*}
thus (\ref{eq::lm::Rpomegac3}) is true for $p=1$.
Moreover, we have
\begin{align*}
R(e_2,e_{2k})e_i=
\begin{cases}
0, \quad \text{for}\quad i\in\{2,\ldots,2k\}\\
S_1e_2,\quad \text{for}\quad i=1
\end{cases}
\end{align*}
Now let us assume that (\ref{eq::lm::Rpomegac3}) holds for some $p\geq 1$.
Using the above formulas we obtain
\begin{align*}
(R^{p+1}\omega)(e_2&,e_{2k},\ldots,e_2,e_{2k},e_i,e_{2k},e_{2k},e_j)\\
&=-R^p\omega(e_2,e_{2k},\ldots,e_2,e_{2k},R(e_2,e_{2k})e_i,e_{2k},e_{2k},e_j)\\
&=
\begin{cases}
0, \quad \text{for}\quad i\in\{2,\ldots,2k\}\\
-R^p\omega(e_2,e_{2k},\ldots,e_2,e_{2k},S_1e_2,e_{2k},e_{2k},e_j),\quad \text{for}\quad i=1
\end{cases}\\
&=
\begin{cases}
0, \quad \text{for}\quad i\in\{2,\ldots,2k\}\\
-R^p\omega(e_2,e_{2k},\ldots,e_2,e_{2k},\beta e_1+\alpha e_2+e_4,e_{2k},e_{2k},e_j),\quad \text{for}\quad i=1
\end{cases}\\
&=
\begin{cases}
0, \quad \text{for}\quad i\in\{2,\ldots,2k\}\\
-\beta R^p\omega(e_2,e_{2k},\ldots,e_2,e_{2k},e_1,e_{2k},e_{2k},e_j),\quad \text{for}\quad i=1
\end{cases}\\
&=
\begin{cases}
0, \quad \text{for}\quad i\in\{2,\ldots,2k\}\\
-\beta \cdot (-\beta)^{p-1}\omega (S_1e_{2k},e_j)=(-\beta)^p\omega (S_1e_{2k},e_j),\quad \text{for}\quad i=1
\end{cases}
\end{align*}
Now by the induction principle (\ref{eq::lm::Rpomegac3}) holds for all $p\geq 1$.
\end{proof}
\end{lem}

Now we can prove

\begin{cor}\label{cor::omegaXej}
Let $S_1$ be a $2k$-dimensional complex Jordan block, $k\geq 1$. If $R^p\omega=0$ for some $p\geq 1$ then for every $X\in \operatorname{span} \{e_1,\ldots,e_{2k}\}:=V$ and $j>2k$
\begin{align}
\omega (X,e_j)=0.
\end{align}
\begin{proof}
For $k=1$ Corollary \ref{cor::omegaXej} is an immediate consequence of Lemma \ref{lm::ComplexBlocks1}.
For $k\geq 2$ by Lemma \ref{lm::Rpomegac1}, Lemma \ref{lm::Rpomegac2} and Lemma \ref{lm::Rpomegac3}  we get
\begin{align*}
\omega(S_1e_i,e_j)=0 \quad \text{for}\quad i\in \{1,\ldots,2k\}.
\end{align*}
Since $\det S_1\neq 0$ and $S_1\colon V\rightarrow V$ is an isomorphism (so $\{S_1e_1,\ldots,S_1e_{2k}\}$ generate $V$) we obtain that
\begin{align*}
\omega(SX,e_j)=0
\end{align*}
for all $X\in V$.
\end{proof}
\end{cor}

In the next few lemmas we study intrinsic properties of  $2k$-dimensional complex Jordan block for $k\geq 2$.


\begin{lem}\label{lm::complexblock2}
Let $S_1$ be a $2k$-dimensional complex Jordan block, $k\geq 2$, $i\in \{1,\ldots,2k-1\}\setminus \{2\}$.
Then for every $p\geq 1$
\begin{align}
\label{eq::q1}R^p\omega(e_1,e_{2k-1},\ldots,e_1,e_{2k-1},e_i,e_{2k-1})&=0,\\
\label{eq::q2}R^p\omega(e_1,e_{2k-1},\ldots,e_1,e_{2k-1},e_i,e_{2k})&=(-\beta)^{p-1}\omega(e_i,S_1e_{2k-1}).
\end{align}
\begin{proof}
The Gauss equation implies that
\begin{align*}
R(e_1,e_{2k-1})e_1=R(e_1,e_{2k-1})e_{2k-1}=R(e_1,e_{2k-1})e_i=0.
\end{align*}
By straightforward computations we get (\ref{eq::q1}) for every $p\geq 1$.
In order to prove (\ref{eq::q2}) again by the Gauss equation we have
\begin{align*}
R(e_1,e_{2k-1})e_{2k}=-S_1e_{2k-1}=-\alpha e_{2k-1}+\beta e_{2k}.
\end{align*}
In particular, for $p=1$, we get
\begin{align*}
R\omega(e_1,e_{2k-1},e_i,e_{2k})=\omega (e_i,S_1e_{2k-1}).
\end{align*}
Now, assume that the formula (\ref{eq::q2}) is true for some $p\geq 1$. Then, for $p+1$, we get
\begin{align*}
(R^{p+1}\omega)(e_1&,e_{2k-1},\ldots,e_1,e_{2k-1},e_i,e_{2k})\\
&=(R(e_1,e_{2k-1})\cdot R^p\omega)(e_1,e_{2k-1},\ldots,e_1,e_{2k-1},e_i,e_{2k})\\
&=-R^p\omega(e_1,e_{2k-1},\ldots,e_1,e_{2k-1},e_i,-S_1e_{2k-1})\\
&=R^p\omega(e_1,e_{2k-1},\ldots,e_1,e_{2k-1},e_i,\alpha e_{2k-1}-\beta e_{2k})\\
&=\alpha R^p\omega(e_1,e_{2k-1},\ldots,e_1,e_{2k-1},e_i,e_{2k-1})-\beta R^p\omega(e_1,e_{2k-1},\ldots,e_1,e_{2k-1},e_i,e_{2k})
\end{align*}
Now, by (\ref{eq::q1}) and the induction principle, the formula (\ref{eq::q2}) holds for every $p\geq 1$.
\end{proof}
\end{lem}


\begin{lem}\label{lm::complexblock3}
Let $S_1$ be a $2k$-dimensional complex Jordan block, $k\geq 2$, $p\geq 1$, $i\in \{1,\ldots, 2k\}\setminus \{1,2k-1\}$.
Then
\begin{align}
\label{eq::t1}R^p\omega(e_2,e_{2k},\ldots,e_2,e_{2k},e_i,e_{2k})&=0,\\
\label{eq::t2}R^p\omega(e_2,e_{2k},\ldots,e_2,e_{2k},e_i,e_{2k-1})&=\beta^{p-1}\omega(e_i,S_1e_{2k}).
\end{align}
\begin{proof}
The Gauss equation implies that
\begin{align*}
R(e_2,e_{2k})e_2=R(e_2,e_{2k})e_{2k}=R(e_2,e_{2k})e_i=0.
\end{align*}
By straightforward computations we get (\ref{eq::t1}) for every $p\geq 1$.
In order to prove (\ref{eq::t2}) again by the Gauss equation we have
\begin{align*}
R(e_2,e_{2k})e_{2k-1}=-S_1e_{2k}=-\beta e_{2k-1}-\alpha e_{2k}.
\end{align*}
In particular, for $p=1$, we get
\begin{align*}
R\omega(e_2,e_{2k},e_i,e_{2k-1})=\omega (e_2,S_1e_{2k}).
\end{align*}
Now, assume that the formula (\ref{eq::t2}) is true for some $p\geq 1$. Then, for $p+1$, we get
\begin{align*}
(R^{p+1}\omega)(e_2&,e_{2k},\ldots,e_2,e_{2k},e_i,e_{2k-1})\\
&=(R(e_2,e_{2k})\cdot R^p\omega)(e_2,e_{2k},\ldots,e_2,e_{2k},e_i,e_{2k-1})\\
&=-R^p\omega(e_2,e_{2k},\ldots,e_2,e_{2k},e_i,-S_1e_{2k})\\
&=R^p\omega(e_2,e_{2k},\ldots,e_2,e_{2k},e_i,\beta e_{2k-1}+\alpha e_{2k})\\
&=\alpha R^p\omega(e_2,e_{2k},\ldots,e_2,e_{2k},e_i,e_{2k})+\beta R^p\omega(e_2,e_{2k},\ldots,e_2,e_{2k},e_i,e_{2k-1})
\end{align*}
Now, by (\ref{eq::t1}) and the induction principle, the formula (\ref{eq::t2}) holds for every $p\geq 1$.
\end{proof}
\end{lem}


\begin{lem}\label{lm::complexblock4}
Let $S_1$ be a $2k$-dimensional complex Jordan block, $k\geq 2$, $p\geq 1$.
If
$$\omega(e_j,e_{2k-1})=\omega(e_j,e_{2k})=0$$
for $j\in \{3,\ldots, 2k\}$ then
\begin{align}
\label{eq::Re3}R^p\omega(e_1,e_{2k-1},\ldots,\overbrace{e_{2k-1},e_{2k}}^{(i)}\ldots,e_1,e_{2k-1},e_1,e_{3})=0
\end{align}
for $i\in \{1,\ldots, p\}$
\begin{proof}
For $p=1$ we have
\begin{align*}
R\omega(e_{2k-1},e_{2k},e_1,e_{3})&=-\omega(R(e_{2k-1},e_{2k})e_1,e_3)-\omega(e_1,R(e_{2k-1},e_{2k})e_3)\\
&=-\omega(R(e_{2k-1},e_{2k})e_1,e_3)=-\omega(S_1e_{2k-1},e_3)=0,
\end{align*}
since $\omega(e_3,e_{2k-1})=\omega(e_3,e_{2k})=0$ by assumption.
\par Assume that (\ref{eq::Re3}) holds for some $p\geq 1$ and for all $i\in \{1,\ldots, p\}$.
Let $i_0\in \{1,\ldots, p+1\}$. If $i_0>1$ then we have
\begin{align*}
R^{p+1}&\omega(e_1,e_{2k-1},\ldots,\overbrace{e_{2k-1},e_{2k}}^{(i_0)}\ldots,e_1,e_{2k-1},e_1,e_{3})\\
&=-R^{p}\omega(e_1,e_{2k-1},\ldots,\overbrace{e_{2k-1},R(e_1,e_{2k-1})e_{2k}}^{(i_0-1)}\ldots,e_1,e_{2k-1},e_1,e_{3})\\
&=R^{p}\omega(e_1,e_{2k-1},\ldots,\overbrace{e_{2k-1},S_1 e_{2k-1}}^{(i_0-1)}\ldots,e_1,e_{2k-1},e_1,e_{3})\\
&=R^{p}\omega(e_1,e_{2k-1},\ldots,\overbrace{e_{2k-1},\alpha e_{2k-1}-\beta e_{2k}}^{(i_0-1)}\ldots,e_1,e_{2k-1},e_1,e_{3})\\
&=-\beta R^{p}\omega(e_1,e_{2k-1},\ldots,\overbrace{e_{2k-1},e_{2k}}^{(i_0-1)}\ldots,e_1,e_{2k-1},e_1,e_{3})
\end{align*}
since $R(e_1,e_{2k-1})e_{1}=R(e_1,e_{2k-1})e_{2k-1}=R(e_1,e_{2k-1})e_{3}=0$ and $R(e_1,e_{2k-1})e_{2k}=-S_1e_{2k-1}=-\alpha e_{2k-1}+\beta e_{2k}$.
Now by (\ref{eq::Re3}) we obtain that for $i_0>1$
\begin{align*}
R^{p+1}&\omega(e_1,e_{2k-1},\ldots,\overbrace{e_{2k-1},e_{2k}}^{(i_0)}\ldots,e_1,e_{2k-1},e_1,e_{3})=0.
\end{align*}
If $i_0=1$ we compute
\begin{align*}
R^{p+1}&\omega(\overbrace{e_{2k-1},e_{2k}}^{(1)},e_1,e_{2k-1},\ldots,\ldots,e_1,e_{2k-1},e_1,e_{3})\\
&=-R^{p}\omega(R(e_{2k-1},e_{2k})e_1,e_{2k-1},e_1,e_{2k-1},\ldots,e_1,e_{2k-1},e_1,e_{3})\\
&\phantom{=}-R^{p}\omega(e_1,e_{2k-1},R(e_{2k-1},e_{2k})e_1,e_{2k-1},\ldots,e_1,e_{2k-1},e_1,e_{3})\\
&\phantom{=}\cdots\\
&\phantom{=}-R^{p}\omega(e_1,e_{2k-1},e_1,e_{2k-1},\ldots,R(e_{2k-1},e_{2k})e_1,e_{2k-1},e_1,e_{3})\\
&\phantom{=}-R^{p}\omega(e_1,e_{2k-1},e_1,e_{2k-1},\ldots,e_1,e_{2k-1},R(e_{2k-1},e_{2k})e_1,e_{3})\\
&=\beta R^{p}\omega(\overbrace{e_{2k},e_{2k-1}}^{(1)},e_1,e_{2k-1},\ldots,e_1,e_{2k-1},e_1,e_{3})\\
&\phantom{=}+\beta R^{p}\omega(e_1,e_{2k-1},\overbrace{e_{2k},e_{2k-1}}^{(2)},\ldots,e_1,e_{2k-1},e_1,e_{3})\\
&\phantom{=}\cdots\\
&\phantom{=}+\beta R^{p}\omega(e_1,e_{2k-1},e_1,e_{2k-1},\ldots,\overbrace{e_{2k},e_{2k-1}}^{(p)},e_1,e_{3})\\
&\phantom{=}-R^{p}\omega(e_1,e_{2k-1},e_1,e_{2k-1},\ldots,e_1,e_{2k-1},S_1 e_{2k-1},e_{3})\\
&=-R^{p}\omega(e_1,e_{2k-1},e_1,e_{2k-1},\ldots,e_1,e_{2k-1},S_1 e_{2k-1},e_{3}),
\end{align*}
since all terms but last are equal 0 thanks to (\ref{eq::Re3}). Now it is enough to show that
$$
R^{p}\omega(e_1,e_{2k-1},e_1,e_{2k-1},\ldots,e_1,e_{2k-1},S_1 e_{2k-1},e_{3})=0.
$$
Indeed, by Lemma \ref{lm::complexblock2} we have
\begin{align*}
R^{p}&\omega(e_1,e_{2k-1},e_1,e_{2k-1},\ldots,e_1,e_{2k-1},S_1 e_{2k-1},e_{3})\\
&=R^{p}\omega(e_1,e_{2k-1},e_1,e_{2k-1},\ldots,e_1,e_{2k-1},\alpha e_{2k-1}-\beta e_{2k},e_{3})\\
&=-\alpha R^{p}\omega(e_1,e_{2k-1},e_1,e_{2k-1},\ldots,e_1,e_{2k-1},e_{3},e_{2k-1})\\
&\phantom{=}+\beta R^{p}\omega(e_1,e_{2k-1},e_1,e_{2k-1},\ldots,e_1,e_{2k-1},e_{3},e_{2k})\\
&=\beta\cdot (-\beta)^{p-1}\omega(e_3,S_1 e_{2k-1})=0,
\end{align*}
where the last equality follows from the assumption that $\omega(e_3,e_{2k-1})=\omega(e_3,e_{2k})=0$.
Summarising, we have shown that
$$
R^{p+1}\omega(e_1,e_{2k-1},\ldots,\overbrace{e_{2k-1},e_{2k}}^{(i_0)}\ldots,e_1,e_{2k-1},e_1,e_{3})=0
$$
for all $i_0\in\{1,\ldots,p+1\}$. Now by the induction principle (\ref{eq::Re3}) holds for all $p\geq 1$.
\end{proof}
\end{lem}


\begin{lem}\label{lm::complexblock4.5}
Let $S_1$ be a $2k$-dimensional complex Jordan block, $k\geq 2$, $p\geq 1$.
If $\omega(e_3,e_{2k-1})=\omega(e_3,e_{2k})=0$ then
\begin{align}
\label{eq::complexblock4.5a} R^p\omega(e_1,e_{2k-1},\ldots,e_1,e_{2k-1},e_2,e_{2k-1})&=\beta ^{p-1}(-\alpha \omega (e_1,e_{2k-1})+\beta \omega(e_2,e_{2k-1}))\\
\label{eq::complexblock4.5b} R^p\omega(e_1,e_{2k-1},\ldots,e_1,e_{2k-1},e_2,e_{2k})&=\alpha \beta ^{p-2}(-\alpha \omega (e_1,e_{2k-1})+\beta \omega (e_2,e_{2k-1}))\\
\nonumber &\hspace{-1cm}-\alpha ^2(-\beta)^{p-2}\omega (e_1,e_{2k-1})-\alpha (-\beta)^{p-1}\omega (e_1,e_{2k})
\end{align}
\begin{proof}
For $p=1$ we compute
\begin{align*}
R\omega (e_1,e_{2k-1},e_2,e_{2k-1})&=-\omega (R(e_1,e_{2k-1})e_2,e_{2k-1})-\omega (e_2,R(e_1,e_{2k-1})e_{2k-1})\\
&=-\omega (\alpha e_1-\beta e_2+e_3,e_{2k-1})\\
&=-\alpha \omega (e_1,e_{2k-1})+\beta \omega (e_2,e_{2k-1}),
\end{align*}
since by assumption $\omega (e_3,e_{2k-1})=0$.
Now, assume that the formula (\ref{eq::complexblock4.5a}) is true for some $p\geq 1$. Then, for $p+1$, we get
\begin{align*}
(R^{p+1}\omega)(e_1&,e_{2k-1},\ldots,e_1,e_{2k-1},e_2,e_{2k-1})\\
&=(R(e_1,e_{2k-1})\cdot R^p\omega)(e_1,e_{2k-1},\ldots,e_1,e_{2k-1},e_2,e_{2k-1})\\
&=-R^p\omega(e_1,e_{2k-1},\ldots,e_1,e_{2k-1},R(e_1,e_{2k-1})e_2,e_{2k-1})\\
&=-R^p\omega(e_1,e_{2k-1},\ldots,e_1,e_{2k-1},S_1e_1,e_{2k-1})\\
&=-R^p\omega(e_1,e_{2k-1},\ldots,e_1,e_{2k-1},\alpha e_1-\beta e_2+e_3,e_{2k-1})\\
&=-\alpha R^p\omega(e_1,e_{2k-1},\ldots,e_1,e_{2k-1},e_1,e_{2k-1})\\
&\phantom{=}+\beta R^p\omega(e_1,e_{2k-1},\ldots,e_1,e_{2k-1},e_2,e_{2k-1})\\
&\phantom{=}-R^p\omega(e_1,e_{2k-1},\ldots,e_1,e_{2k-1},e_3,e_{2k-1})\\
&=\beta R^p\omega(e_1,e_{2k-1},\ldots,e_1,e_{2k-1},e_2,e_{2k-1})
\end{align*}
where the last equality follows from Lemma \ref{lm::complexblock2} (formula (\ref{eq::q1}), $i=1,3$).
Now, using (\ref{eq::complexblock4.5a})  we obtain
\begin{align*}
(R^{p+1}\omega)(e_1&,e_{2k-1},\ldots,e_1,e_{2k-1},e_2,e_{2k-1})\\
&=\beta R^p\omega(e_1,e_{2k-1},\ldots,e_1,e_{2k-1},e_2,e_{2k-1})\\
&=\beta (\beta ^{p-1}(-\alpha \omega (e_1,e_{2k-1})+\beta \omega(e_2,e_{2k-1})))\\
&=\beta ^p (-\alpha \omega (e_1,e_{2k-1})+\beta \omega(e_2,e_{2k-1})).
\end{align*}
By the induction principle (\ref{eq::complexblock4.5a}) holds for all $p\geq 1$.
\par The formula (\ref{eq::complexblock4.5b}) can be easily obtained in a similar way using (\ref{eq::q2}), (\ref{eq::complexblock4.5a}) and the principle of induction.
\end{proof}
\end{lem}

From Lemma \ref{lm::complexblock4.5} we immediately get


\begin{cor}\label{cor::help1}
Let $S_1$ be a $2k$-dimensional complex Jordan block, $k\geq 2$, $p\geq 1$.
If $\omega(e_3,e_{2k-1})=\omega(e_3,e_{2k})=0$ then
\begin{align}
R^p\omega(e_1,e_{2k-1},\ldots,e_1,e_{2k-1},e_2,S_1e_{2k-1})&=(-1)^p\beta ^{p-1}\alpha(\alpha \omega (e_1,e_{2k-1})-\beta \omega(e_1,e_{2k}))
\end{align}
\end{cor}


\begin{lem}\label{lm::complexblock5}
Let $S_1$ be a $2k$-dimensional complex Jordan block, $k\geq 2$, $p\geq 1$.
If $\omega(e_3,e_{2k-1})=\omega(e_3,e_{2k})=0$
then for $i\in\{1,\ldots,p\}$
\begin{align}
\label{eq::Re2}R^p&\omega(e_1,e_{2k-1},\ldots,\overbrace{e_{2k-1},e_{2k}}^{(i)}\ldots,e_1,e_{2k-1},e_1,e_{2})\\
&=\nonumber
\begin{cases}
  (-\beta)^{p-1}(-\omega(S_1 e_{2k-1},e_2)+\omega(e_1,S_1 e_{2k})), & \text{if $i=1$} \\
  0, & \text{if $i>1$} .
\end{cases}
\end{align}
\begin{proof}
We directly check that
\begin{align*}
R\omega(e_{2k-1},e_{2k},e_1,e_{2})&=-\omega(R(e_{2k-1},e_{2k})e_1,e_2)-\omega(e_1,R(e_{2k-1},e_{2k})e_2)\\
&=-\omega(S_1 e_{2k-1},e_2)+\omega(e_1,S_1 e_{2k}),
\end{align*}
so (\ref{eq::Re2}) is true for $p=1$.
Now let us assume that (\ref{eq::Re2}) is true for some $p\geq 1$ and for all $i\in\{1,\ldots,p\}$.
Let $i_0\in\{1,\ldots,p+1\}$.
\par When $i_0>1$ using the fact that $R(e_1,e_{2k-1})e_1=R(e_1,e_{2k-1})e_{2k-1}=0$ we get
\begin{align*}
R^{p+1}&\omega(e_1,e_{2k-1},\ldots,\overbrace{e_{2k-1},e_{2k}}^{(i_0)}\ldots,e_1,e_{2k-1},e_1,e_{2})\\
&=-R^{p}\omega(e_1,e_{2k-1},\ldots,\overbrace{e_{2k-1},R(e_1,e_{2k-1})e_{2k}}^{(i_0-1)}\ldots,e_1,e_{2k-1},e_1,e_{2})\\
&\phantom{=}-R^{p}\omega(e_1,e_{2k-1},\ldots,\overbrace{e_{2k-1},e_{2k}}^{(i_0-1)}\ldots,e_1,e_{2k-1},e_1,R(e_1,e_{2k-1})e_{2})\\
&=R^{p}\omega(e_1,e_{2k-1},\ldots,\overbrace{e_{2k-1},S_1 e_{2k-1}}^{(i_0-1)}\ldots,e_1,e_{2k-1},e_1,e_{2})\\
&\phantom{=}-R^{p}\omega(e_1,e_{2k-1},\ldots,\overbrace{e_{2k-1},e_{2k}}^{(i_0-1)}\ldots,e_1,e_{2k-1},e_1,S e_{1})\\
&=-\beta R^{p}\omega(e_1,e_{2k-1},\ldots,\overbrace{e_{2k-1},e_{2k}}^{(i_0-1)}\ldots,e_1,e_{2k-1},e_1,e_{2})\\
&\phantom{=}+\beta R^{p}\omega(e_1,e_{2k-1},\ldots,\overbrace{e_{2k-1},e_{2k}}^{(i_0-1)}\ldots,e_1,e_{2k-1},e_1,e_{2})\\
&\phantom{=}-R^{p}\omega(e_1,e_{2k-1},\ldots,\overbrace{e_{2k-1},e_{2k}}^{(i_0-1)}\ldots,e_1,e_{2k-1},e_1,e_{3})\\
&=-R^{p}\omega(e_1,e_{2k-1},\ldots,\overbrace{e_{2k-1},e_{2k}}^{(i_0-1)}\ldots,e_1,e_{2k-1},e_1,e_{3})=0
\end{align*}
 where the last equality follows from Lemma \ref{lm::complexblock4}.
\par If $i_0=1$, using the fact that $R(e_{2k-1},e_{2k})e_{2k-1}=0$ we obtain
\begingroup
\allowdisplaybreaks
\begin{align*}
R^{p+1}&\omega(\overbrace{e_{2k-1},e_{2k}}^{(1)},e_1,e_{2k-1},\ldots,e_1,e_{2k-1},e_1,e_{2})\\
&=-R^{p}\omega(R(e_{2k-1},e_{2k})e_1,e_{2k-1},e_1,e_{2k-1},\ldots,e_1,e_{2k-1},e_1,e_{2})\\
&\phantom{=}-R^{p}\omega(e_1,e_{2k-1},R(e_{2k-1},e_{2k})e_1,e_{2k-1},\ldots,e_1,e_{2k-1},e_1,e_{2})\\
&\phantom{=}\cdots\\
&\phantom{=}-R^{p}\omega(e_1,e_{2k-1},e_1,e_{2k-1},\ldots,R(e_{2k-1},e_{2k})e_1,e_{2k-1},e_1,e_{2})\\
&\phantom{=}-R^{p}\omega(e_1,e_{2k-1},e_1,e_{2k-1},\ldots,e_1,e_{2k-1},R(e_{2k-1},e_{2k})e_1,e_{2})\\
&\phantom{=}-R^{p}\omega(e_1,e_{2k-1},e_1,e_{2k-1},\ldots,e_1,e_{2k-1},e_1,R(e_{2k-1},e_{2k})e_{2})\\
&=\beta R^{p}\omega(\overbrace{e_{2k},e_{2k-1}}^{(1)},e_1,e_{2k-1},\ldots,e_1,e_{2k-1},e_1,e_{2})\\
&\phantom{=}+\beta R^{p}\omega(e_1,e_{2k-1},\overbrace{e_{2k},e_{2k-1}}^{(2)},\ldots,e_1,e_{2k-1},e_1,e_{2})\\
&\phantom{=}\cdots\\
&\phantom{=}+\beta R^{p}\omega(e_1,e_{2k-1},e_1,e_{2k-1},\ldots,\overbrace{e_{2k},e_{2k-1}}^{(p)},e_1,e_{2})\\
&\phantom{=}-R^{p}\omega(e_1,e_{2k-1},e_1,e_{2k-1},\ldots,e_1,e_{2k-1},S_1 e_{2k-1},e_{2})\\
&\phantom{=}+R^{p}\omega(e_1,e_{2k-1},e_1,e_{2k-1},\ldots,e_1,e_{2k-1},e_1,S_1 e_{2k})\\
&=-\beta R^{p}\omega(\overbrace{e_{2k-1},e_{2k}}^{(1)},e_1,e_{2k-1},\ldots,e_1,e_{2k-1},e_1,e_{2})\\
&\phantom{=}+R^{p}\omega(e_1,e_{2k-1},e_1,e_{2k-1},\ldots,e_1,e_{2k-1},e_{2},S_1 e_{2k-1})\\
&\phantom{=}+R^{p}\omega(e_1,e_{2k-1},e_1,e_{2k-1},\ldots,e_1,e_{2k-1},e_1,S_1 e_{2k}),\\
\end{align*}
\endgroup
since all but first and the last two terms are equal to 0 by (\ref{eq::Re2}).
Now Lemma \ref{lm::complexblock2} and Corollary \ref{cor::help1} imply that
\begin{align*}
&R^{p}\omega(e_1,e_{2k-1},e_1,e_{2k-1},\ldots,e_1,e_{2k-1},e_1,S_1 e_{2k})\\
&+R^{p}\omega(e_1,e_{2k-1},e_1,e_{2k-1},\ldots,e_1,e_{2k-1},e_{2},S_1 e_{2k-1})\\
&=\alpha(-\beta)^{p-1}\omega(e_1,S_1e_{2k-1})+(-1)^p\beta ^{p-1}\alpha(\alpha \omega (e_1,e_{2k-1})-\beta \omega(e_1,e_{2k}))\\
&=0.
\end{align*}
In consequence, using (\ref{eq::Re2}) we get
\begin{align*}
R^{p+1}&\omega(\overbrace{e_{2k-1},e_{2k}}^{(1)},e_1,e_{2k-1},\ldots,e_1,e_{2k-1},e_1,e_{2})\\
&=-\beta R^{p}\omega(\overbrace{e_{2k-1},e_{2k}}^{(1)},e_1,e_{2k-1},\ldots,e_1,e_{2k-1},e_1,e_{2})\\
&=(-\beta)^{p}(-\omega(S_1 e_{2k-1},e_2)+\omega(e_1,S_1 e_{2k})).
\end{align*}
Now the thesis follows from the induction principle.
\end{proof}
\end{lem}

To simplify further notation let us denote:
\begin{align}\label{eq::Aijp}
A_{ij}^p:=R^p\omega(e_1,e_{2k-1},e_1,e_{2k-1},\ldots,e_1,e_{2k-1},e_i,e_{2k-1},e_j,e_{2k-1})
\end{align}
where $i,j\in\{1,2,3\}$ and $p\geq 1$.
We have the following lemma:


\begin{lem}
Let $S_1$ be a $2k$-dimensional complex Jordan block, $k\geq 2$, $p\geq 1$. Then
\begin{align}
\label{eq::LemAijp::1} A_{11}^p&=A_{13}^p=A_{31}^p=A_{33}^p=0,\\
\label{eq::LemAijp::2} A_{12}^{p+1}&=\beta A_{12}^{p},\\
\label{eq::LemAijp::3} A_{21}^{p+1}&=\beta A_{21}^{p},\\
\label{eq::LemAijp::4} A_{23}^{p+1}&=\beta A_{23}^{p},\\
\label{eq::LemAijp::5} A_{32}^{p+1}&=\beta A_{32}^{p},\\
\label{eq::LemAijp::6} A_{22}^{p+1}&=-\alpha(A_{12}^{p}+A_{21}^{p})+2\beta A_{22}^{p}-(A_{32}^{p}+A_{23}^{p}).
\end{align}
\begin{proof}
In order to prove (\ref{eq::LemAijp::1}) first note that
\begin{align*}
A_{11}^p&=A_{13}^p=0
\end{align*}
immediately follows from Lemma \ref{lm::complexblock2}, formula (\ref{eq::q1}). By the Gauss equation we have
\begin{align*}
R(e_3&,e_{2k-1})e_1=R(e_3,e_{2k-1})e_{2k-1}=R(e_3,e_{2k-1})e_3\\
&=R(e_1,e_{2k-1})e_{1}=R(e_1,e_{2k-1})e_{2k-1}=R(e_1,e_{2k-1})e_{3}=0.
\end{align*}
Using the above equalities we easy obtain that
\begin{align*}
A_{31}^p=A_{33}^p=0.
\end{align*}
To prove (\ref{eq::LemAijp::2}) we compute
\begin{align*}
 A_{12}^{p+1}&=R^{p+1}\omega (e_1,e_{2k-1},\ldots,e_1,e_{2k-1},e_1,e_{2k-1},e_2,e_{2k-1})\\
 &=-R^p\omega (e_1,e_{2k-1},\ldots,e_1,e_{2k-1},e_1,e_{2k-1},R(e_1,e_{2k-1})e_2,e_{2k-1})\\
 &=-R^p\omega (e_1,e_{2k-1},\ldots,e_1,e_{2k-1},e_1,e_{2k-1},S_1e_1,e_{2k-1})\\
 &=-R^p\omega (e_1,e_{2k-1},\ldots,e_1,e_{2k-1},e_1,e_{2k-1},\alpha e_1-\beta e_2+e_3,e_{2k-1}),\\
 &=-\alpha A_{11}^{p}+\beta A_{12}^{p}-A_{13}^{p}=\beta A_{12}^{p},
\end{align*}
where the last equality follows from (\ref{eq::LemAijp::1}).
The formulas (\ref{eq::LemAijp::3})--(\ref{eq::LemAijp::5}) we prove in a similar way like (\ref{eq::LemAijp::2}).
To prove (\ref{eq::LemAijp::6}) we compute
\begin{align*}
 A_{22}^{p+1}&=R^{p+1}\omega (e_1,e_{2k-1},\ldots,e_1,e_{2k-1},e_2,e_{2k-1},e_2,e_{2k-1})\\
 &=-R^p\omega (e_1,e_{2k-1},\ldots,e_1,e_{2k-1},R(e_1,e_{2k-1})e_2,e_{2k-1},e_2,e_{2k-1})\\
 &\phantom{=}-R^p\omega (e_1,e_{2k-1},\ldots,e_1,e_{2k-1},R(e_1,e_2,e_{2k-1},R(e_1,e_{2k-1})e_2,e_{2k-1})\\
 &=-R^p\omega (e_1,e_{2k-1},\ldots,e_1,e_{2k-1},S_1e_1,e_{2k-1},e_2,e_{2k-1})\\
 &\phantom{=}-R^p\omega (e_1,e_{2k-1},\ldots,e_1,e_{2k-1},e_2,e_{2k-1},S_1e_1,e_{2k-1})\\
 &=-R^p\omega (e_1,e_{2k-1},\ldots,e_1,e_{2k-1},\alpha e_1-\beta e_2+e_3,e_{2k-1},e_2,e_{2k-1})\\
 &\phantom{=}-R^p\omega (e_1,e_{2k-1},\ldots,e_1,e_{2k-1},e_2,e_{2k-1},\alpha e_1-\beta e_2+e_3,e_{2k-1},e_2,e_{2k-1})\\
 &=-\alpha A_{12}^{p}+\beta A_{22}^{p}-A_{32}^{p}-\alpha A_{21}^{p}+\beta A_{22}^{p}-A_{23}^{p}\\
 &=-\alpha ( A_{12}^{p}+ A_{21}^{p})+2\beta A_{22}^{p}-(A_{32}^{p}+A_{23}^{p}).
\end{align*}
\end{proof}
\end{lem}

From the above lemma we obtain:


\begin{cor}\label{cor::ComplexCase}
Let $S_1$ be a $2k$-dimensional complex Jordan block, $k\geq 2$, $p\geq 1$.
If
\begin{align}\label{eq::omegaj2k2k-1Zero}
\omega(e_j,e_{2k-1})=\omega(e_j,e_{2k})=0
\end{align}
for $j\in \{3,\ldots, 2k\}$ then
\begin{align}
\label{eq::a12p}A_{12}^p&=\beta^{p-1}(-\alpha \omega (e_1,e_{2k-1})+\beta \omega (e_2,e_{2k-1})),\\
\label{eq::a21p}A_{21}^p&=\beta^{p-1}(\alpha \omega (e_1,e_{2k-1})-\beta \omega (e_1,e_{2k})),\\
\label{eq::a22p1} A_{22}^{p+1}&=-\alpha \beta ^p(\omega (e_2,e_{2k-1})-\omega (e_1,e_{2k}))+2\beta A_{22}^{p}.
\end{align}
\begin{proof}
From formulas (\ref{eq::LemAijp::2}), (\ref{eq::LemAijp::3})  and assumption (\ref{eq::omegaj2k2k-1Zero})  we obtain
\begin{align*}
A_{12}^p&=\beta ^{p-1}\cdot A_{12}^1=\beta ^{p-1}\cdot R \omega (e_1,e_{2k-1},e_2,e_{2k-1})\\
&=-\beta ^{p-1}\omega (R(e_1,e_{2k-1})e_2,e_{2k-1})=-\beta ^{p-1}\omega (S_1e_1,e_{2k-1})\\
&=-\alpha \beta^{p-1}\omega (e_1,e_{2k-1})+\beta ^p\omega (e_2,e_{2k-1}).\\
A_{21}^p&=\beta ^{p-1}\cdot A_{21}^1=\beta ^{p-1}\cdot R \omega (e_2,e_{2k-1},e_1,e_{2k-1})\\
&=-\beta ^{p-1}\omega (e_1,R(e_2,e_{2k-1})e_{2k-1})=\beta ^{p-1}\omega (e_1,S_1e_{2k-1})\\
&=\alpha \beta^{p-1}\omega (e_1,e_{2k-1})-\beta ^p\omega (e_1,e_{2k}),
\end{align*}
what proves (\ref{eq::a12p}) and (\ref{eq::a21p}).
Using (\ref{eq::a12p}) and (\ref{eq::a21p}) we obtain
\begin{align*}
A_{12}^p+A_{21}^p=\beta ^p (\omega (e_2,e_{2k-1})-\omega (e_1,e_{2k})).
\end{align*}
Moreover from formulas (\ref{eq::LemAijp::4}), (\ref{eq::LemAijp::5}) and (\ref{eq::omegaj2k2k-1Zero}) we get
\begin{align*}
A_{32}^p&=\beta ^{p-1}A_{32}^1=\beta ^{p-1}R \omega (e_3,e_{2k-1},e_2,e_{2k-1})\\
&=-\beta ^{p-1}\omega (R(e_3,e_{2k-1})e_2,e_{2k-1})=-\beta ^{p-1}\omega (S_1e_3,e_{2k-1})=0,\\
A_{23}^p&=\beta ^{p-1}\cdot A_{23}^1=0.
\end{align*}
Now (\ref{eq::a22p1}) immediately follows from (\ref{eq::LemAijp::6}).
\end{proof}
\end{cor}


\begin{lem}\label{lem::omegai2kkminus12k}
Let $S_1$ be a $2k$-dimensional complex Jordan block, $k\geq 2$. If $R^p\omega=0$ for some $p\geq 1$ then for $i\in \{3,\ldots, 2k\}$
\begin{align*}
\omega (e_i,e_{2k-1})=\omega (e_i,e_{2k})=0.
\end{align*}
\begin{proof}
From Lemma \ref{lm::complexblock2}, for $i=3,\ldots,2k-1$ we have
\begin{align}
\label{eq::complexblockq1} \omega (e_i,S_1e_{2k-1})=0.
\end{align}
In particular for $i=2k-1$ we obtain
\begin{align*}
0=\omega(e_{2k-1},S_1e_{2k-1})=\omega(e_{2k-1},\alpha e_{2k-1}-\beta e_{2k})=-\beta \omega(e_{2k-1},e_{2k}).
\end{align*}
Thus
\begin{align*}
\omega(e_{2k-1},e_{2k})=0,
\end{align*}
since $\beta \neq 0$.
Now we have
\begin{align*}
\omega(e_{2k},S_1e_{2k-1})=\omega(e_{2k},\alpha e_{2k-1}-\beta e_{2k})=\alpha \omega(e_{2k},e_{2k-1})=0.
\end{align*}
That is (\ref{eq::complexblockq1}) is also true for $i=2k$. From Lemma \ref{lm::complexblock3} for $i\in \{3,\ldots,2k\}\setminus \{2k-1\}$ we have
\begin{align}
\label{eq::complexblockq2} \omega(e_i,S_1e_{2k})=0.
\end{align}
Since $\omega(e_{2k-1},e_{2k})=0$ we also get that
\begin{align*}
\omega(e_{2k-1},S_1e_{2k})=\omega(e_{2k-1},\beta e_{2k-1}+\alpha e_{2k})=\alpha \omega(e_{2k-1},e_{2k})=0.
\end{align*}
That is (\ref{eq::complexblockq2}) also valid for $i=2k-1$. Now from (\ref{eq::complexblockq1}) and (\ref{eq::complexblockq2}), for $i\in \{3,\ldots,2k\}$ we get
\begin{align*}
0&=\alpha \omega(e_i,S_1e_{2k-1})+\beta \omega(e_i,S_1e_{2k})\\
&=\omega(e_i,\alpha \cdot (\alpha e_{2k-1}-\beta e_{2k})+\beta \cdot (\beta e_{2k-1}+\alpha e_{2k}))\\
&=(\alpha ^2+\beta ^2)\omega(e_i,e_{2k-1}).
\end{align*}
In a similar way we compute that
\begin{align*}
0&=-\beta \omega(e_i,S_1e_{2k-1})+\alpha \omega(e_i,S_1e_{2k})\\
&=(\alpha ^2+\beta ^2)\omega(e_i,e_{2k}).
\end{align*}
Since $\alpha ^2+\beta ^2\neq 0$ the above equations imply that
\begin{align*}
\omega(e_i,e_{2k-1})=\omega(e_i,e_{2k})=0
\end{align*}
for $i\in \{3,\ldots,2k\}$.
\end{proof}
\end{lem}

Now we are ready to prove that the decomposition from the Lemma \ref{lm::JordanDecompositionOfS} cannot contain complex Jordan blocks. Namely we have the following:
\begin{thm}\label{tw::NoComplex}
Let $f\colon \M\rightarrow\R^{2n+1}$ ($\dim M\geq 4$)  be a non-degenerate affine hypersurface with a locally equiaffine transversal vector field  $\xi$
and an almost symplectic form $\omega$. If $R^p\omega=0$ for some $p\geq 1$ and
\begin{align}\label{eq::SDecompositionComplex}
S=\left[\begin{matrix}
S_1 & 0 & \ldots & 0 \\
0 & S_2 & \ldots & 0 \\
\vdots & \vdots & \ddots & \ldots \\
0 & 0 & \ldots & S_{q+r}
\end{matrix}\right]
\end{align}
is the Jordan decomposition of $S$ as stated in the Lemma \ref{lm::JordanDecompositionOfS} then (\ref{eq::SDecompositionComplex}) does not contain complex Jordan blocks
(that is $r=0)$.
\begin{proof}
\par Let $\{e_1,\ldots,e_{2n}\}$ be the basis of $T_xM$ from Lemma \ref{lm::JordanDecompositionOfS}.
Without loss of generality, as described at the beginning of this section, we can
change order of $S_i$ and $H_i$ in such way that $S_{1},\ldots,S_{r}$ will be complex blocks and $S_{r+1},\ldots,S_{r+q}$ will be
real blocks. Moreover, we can assume that $\dim S_1\geq \dim S_i$ for $i=2,\ldots,r$.
\par First assume that $S_1$ is a complex block of dimension $2k$ and $k\geq 2$.
By Lemma \ref{lem::omegai2kkminus12k} we have that $\omega (e_i,e_{2k-1})=\omega (e_i,e_{2k})=0$ for $i\in \{3,\ldots, 2k\}$.
Now from Corollary \ref{cor::ComplexCase} (formulas (\ref{eq::a12p}) and (\ref{eq::a21p})) we get
\begin{align}\label{eq::lemRpOmega1}
\omega (e_2,e_{2k-1})=\frac{\alpha}{\beta} \omega (e_1,e_{2k-1})
\end{align}
and
\begin{align}\label{eq::lemRpOmega2}
 \omega (e_1,e_{2k})= \frac{\alpha}{\beta}\omega (e_1,e_{2k-1}),
\end{align}
since $R^p\omega=0$ and $\beta\neq 0$. In particular $\omega (e_2,e_{2k-1})=\omega (e_1,e_{2k})$ and (\ref{eq::a22p1}) simplify to the form
$$
A_{22}^{p+1}=2\beta A_{22}^{p}.
$$
Now one can easily find explicit formula for $A_{22}^{p}$. Namely we have
\begin{align*}
A_{22}^{p}=2^{p-1}\beta^{p-1}A_{22}^{1}&=2^{p-1}\beta^{p-1}R\omega(e_2,e_{2k-1},e_2,e_{2k-1})\\
&=2^{p-1}\beta^{p-1}(-\omega(S_1 e_2,e_{2k-1})+\omega(e_2,S_1 e_{2k-1}))\\
&=2^{p-1}\beta^{p-1}(-\beta\omega(e_1,e_{2k-1})-\alpha\omega(e_2,e_{2k-1})\\
&\phantom{=}+\alpha\omega(e_2,e_{2k-1})-\beta\omega(e_2,e_{2k}))\\
&=-2^{p-1}\beta^{p}(\omega(e_1,e_{2k-1})+\omega(e_2,e_{2k})).
\end{align*}
On the other hand we have $A_{22}^{p}=0$ (since $R^p\omega=0$) and in consequence
\begin{align}\label{eq::lemRpOmega3}
\omega(e_2,e_{2k})=-\omega(e_1,e_{2k-1}).
\end{align}
From Lemma \ref{lm::complexblock5} (formula (\ref{eq::Re2})) we obtain
\begin{align*}
\omega(S_1 e_{2k-1},e_2)-\omega(e_1,S_1 e_{2k})=0
\end{align*}
that is
\begin{align*}
\alpha\omega(e_{2k-1},e_2)-\beta\omega(e_{2k},e_2)-\beta\omega(e_1,e_{2k-1})-\alpha\omega(e_1,e_{2k})=0.
\end{align*}
Now using (\ref{eq::lemRpOmega1})--(\ref{eq::lemRpOmega3}) the above implies that
\begin{align*}
-\alpha\cdot\frac{\alpha}{\beta} \omega (e_1,e_{2k-1})-\beta\omega(e_1,e_{2k-1})-\beta\omega(e_1,e_{2k-1})-\alpha\cdot\frac{\alpha}{\beta}\omega (e_1,e_{2k-1})\\
=-2(\frac{\alpha^2}{\beta}+\beta)\omega (e_1,e_{2k-1})=-\frac{2}{\beta}(\alpha^2+\beta^2)\omega (e_1,e_{2k-1})=0.
\end{align*}
In this way we have shown that $\omega (e_1,e_{2k-1})=0$ and in consequence also $\omega (e_2,e_{2k-1})=0$.
Hence
$$
\omega (e_i,e_{2k-1})=0
$$
for $i\in \{1,\ldots, 2k\}$. From Corollary \ref{cor::omegaXej} we also have that
$$
\omega (e_i,e_{2k-1})=0
$$
for $i>2k$, that is $\omega$ is degenerate, what leads to contradiction and we must have $k<2$.
In this way we have shown that if Jordan decomposition of $S$ contains some complex Jordan blocks they all must be 2-dimensional.

\par It remained to show that also 2-dimensional complex Jordan blocks are not possible. In order to prove it let us assume that
$S_1$ is a 2-dimensional complex block. Since $R^p\omega=0$ then also $R^{2p}\omega=0$ and Lemma \ref{lm::ComplexBlocks1} implies that
$\omega (e_1,e_{l})=0$ for $l=3,\ldots,2n$. Now observe that since $\dim M\geq 4$ there must exist $i_0,j_0>2$ such that $h(e_{i_0},e_{j_0})\neq 0$
(otherwise $h$ would be degenerate).
Now from Lemma \ref{lm::OneComplexAndOther} we have
$$
2^{2p-1}\beta^2(\det S_1)^{p-1}h(e_{i_0},e_{j_0})\omega(e_1,e_2)=0.
$$
That is
$$
\omega(e_1,e_2)=0,
$$
since $h(e_{i_0},e_{j_0})\neq 0$. In this way we have shown that $\omega$ is degenerate since $\omega (e_1,e_{l})=0$ for $l=1,\ldots,2n$, what leads us again
to contradiction.
\end{proof}
\end{thm}



\section{Main Results}
Before we proceed with main results of this paper we need to recall the following two lemmas:
\begin{lem}[\cite{MSz}]\label{lm::NablaKOmega0::2}
Let $f\colon \M\rightarrow\R^{2n+1}$ be a non-degenerate affine hypersurface ($\dim M\geq 4$) with a transversal vector field  $\xi$ and an almost symplectic form $\omega$. Let $x\in M$.
If there exist a natural number $2 \leq s\leq 2n$ and a basis $\{e_1,\ldots,e_{2n}\}$ of $T_xM$ such that
$h(e_i,e_j)=\eps_i\delta_{ij}$, $\eps_i=\pm1$ for $i,j =1,\ldots, s$, $h(e_i,e_j)=0$ for $i=1,\ldots, s$, $j=s+1,\ldots, 2n$ and
$Se_i=\lambda_i e_i$ for $i=1,\ldots, s$, $\lambda_i\in\R$.
Then for every $k, j=1,2,\ldots,s$, $k\neq j$ and for every  $i=1,\ldots,2n$, $i\neq j$, $i\neq k$ we have
\begin{equation}\label{eq::R2ke_1e_j}
R^{2l}\omega (e_k,e_j,e_k,e_j,\ldots,e_k,e_j,e_k,e_i)=(-1)^l\eps_k^l\eps_j^l \lambda_k^l\lambda_j^l\omega(e_k,e_i)
\end{equation}
for every $l\in\N$.
\end{lem}

\begin{lem}[\cite{MSz}]\label{lm::R2kOmegaXZ1Z2Y}
Let $f\colon \M\rightarrow\R^{2n+1}$ be a non-degenerate affine hypersurface ($\dim M\geq 4$) with a transversal vector field  $\xi$ and an almost symplectic form $\omega$.
Let $x\in M$ and let $X, Y, Z_1, Z_2$ be vector fields from $T_xM$ such that $SX=\lambda X$, $SZ_1=0$, $SZ_2=0$, $h(Z_1,Z_2)=0$ and $h(Y,Z_2)=0$.
Then, for every  $l\geq 1$ we have
\begin{align}\label{eq::R2kOmegaXZ1Z2Y}
R^{2l}\omega (X&,\underbrace{Z_1,Z_2,Z_1,Z_2,\ldots,Z_1,Z_2}_{4l},Y)\\
\nonumber  &=(-1)^l\lambda^{2l}h(Z_1,Z_1)^lh(Z_2,Z_2)^l\omega(X,Y).
\end{align}
\end{lem}

We will need also below three lemmas.

\begin{lem}\label{lm::OneBock2x2andOne1x1}
Let $\{e_1,e_2,e_3\}$ be vectors such that $Se_1=\alpha e_1+e_2$, $Se_2=\alpha e_2$, $Se_3=\lambda e_3$ and $h(e_1,e_1)=h(e_2,e_2)=h(e_1,e_3)=h(e_2,e_3)=0$, $h(e_1,e_2)=\eta$, $h(e_3,e_3)=\varepsilon$, $\eta\neq 0, \varepsilon\neq 0$.
Then for every $p\geq 1$ we have
\begin{equation}\label{eq::OneBock2x2andOne1x1}
R^{p}\omega (e_1,e_2,e_1,e_2,\ldots,e_1,e_2,e_1,e_3,e_1,e_3)=(-1)^p\cdot (2\eta\alpha)^{p-1}\varepsilon\omega(e_1,e_2).
\end{equation}
\begin{proof}
By straightforward computations we obtain
\begin{align}
\label{eq::OneBock2x2andOne1x1::eq1}R\omega(e_1,e_3,e_1,e_3)&=-\varepsilon\omega(e_1,e_2)\\
\label{eq::OneBock2x2andOne1x1::eq2}R\omega(e_1,e_3,e_2,e_3)&=\varepsilon\alpha\omega(e_1,e_2)\\
\label{eq::OneBock2x2andOne1x1::eq3}R\omega(e_2,e_3,e_1,e_3)&=-\varepsilon\alpha\omega(e_1,e_2)\\
\label{eq::OneBock2x2andOne1x1::eq4}R\omega(e_2,e_3,e_2,e_3)&=0.
\end{align}
Now for $p\geq 1$ and $i,j\in\{1,2\}$ we get
\begin{align*}
R^{p+1}\omega &(e_1,e_2,e_1,e_2,\ldots,e_1,e_2,e_i,e_3,e_j,e_3)\\
&=-R^{p}\omega (e_1,e_2,e_1,e_2,\ldots,e_1,e_2,R(e_1,e_2)e_i,e_3,e_j,e_3)\\
&\phantom{=}-R^{p}\omega (e_1,e_2,e_1,e_2,\ldots,e_1,e_2,e_i,e_3,R(e_1,e_2)e_j,e_3)\\
&=-\alpha h(e_2,e_i)R^{p}\omega (e_1,e_2,e_1,e_2,\ldots,e_1,e_2,e_1,e_3,e_j,e_3)\\
&\phantom{=}+(\alpha h(e_1,e_i)-h(e_2,e_i))R^{p}\omega (e_1,e_2,e_1,e_2,\ldots,e_1,e_2,e_2,e_3,e_j,e_3)\\
&\phantom{=}-\alpha h(e_2,e_j)R^{p}\omega (e_1,e_2,e_1,e_2,\ldots,e_1,e_2,e_i,e_3,e_1,e_3)\\
&\phantom{=}+(\alpha h(e_1,e_j)-h(e_2,e_j))R^{p}\omega (e_1,e_2,e_1,e_2,\ldots,e_1,e_2,e_i,e_3,e_2,e_3).
\end{align*}
Using different configurations of $i$ and $j$ we obtain the following four relations:
\begin{align}
\label{eq::OneBock2x2andOne1x1::eq5}R^{p+1}\omega &(e_1,e_2,e_1,e_2,\ldots,e_1,e_2,e_2,e_3,e_2,e_3)\\
\nonumber&=2\alpha h(e_1,e_2) R^{p}\omega (e_1,e_2,e_1,e_2,\ldots,e_1,e_2,e_2,e_3,e_2,e_3)\\
\label{eq::OneBock2x2andOne1x1::eq6}R^{p+1}\omega &(e_1,e_2,e_1,e_2,\ldots,e_1,e_2,e_1,e_3,e_2,e_3)\\
\nonumber&=-h(e_1,e_2) R^{p}\omega (e_1,e_2,e_1,e_2,\ldots,e_1,e_2,e_2,e_3,e_2,e_3)\\
\label{eq::OneBock2x2andOne1x1::eq7}R^{p+1}\omega &(e_1,e_2,e_1,e_2,\ldots,e_1,e_2,e_2,e_3,e_1,e_3)\\
\nonumber&=-h(e_1,e_2) R^{p}\omega (e_1,e_2,e_1,e_2,\ldots,e_1,e_2,e_2,e_3,e_2,e_3)\\
\label{eq::OneBock2x2andOne1x1::eq8}R^{p+1}\omega &(e_1,e_2,e_1,e_2,\ldots,e_1,e_2,e_1,e_3,e_1,e_3)\\
\nonumber&=-2\alpha h(e_1,e_2) R^{p}\omega (e_1,e_2,e_1,e_2,\ldots,e_1,e_2,e_1,e_3,e_1,e_3)\\
\nonumber&\phantom{=}-h(e_1,e_2) R^{p}\omega (e_1,e_2,e_1,e_2,\ldots,e_1,e_2,e_2,e_3,e_1,e_3)\\
\nonumber&\phantom{=}-h(e_1,e_2) R^{p}\omega (e_1,e_2,e_1,e_2,\ldots,e_1,e_2,e_1,e_3,e_2,e_3).
\end{align}
Now from (\ref{eq::OneBock2x2andOne1x1::eq4}) and (\ref{eq::OneBock2x2andOne1x1::eq5}), by the induction principle we easily obtain that
$$
R^{p}\omega (e_1,e_2,e_1,e_2,\ldots,e_1,e_2,e_2,e_3,e_2,e_3)=0
$$
for all $p\geq 1$. In consequence from (\ref{eq::OneBock2x2andOne1x1::eq6}) and (\ref{eq::OneBock2x2andOne1x1::eq7}) we get that also
\begin{align*}
R^{p}\omega (e_1,e_2,e_1,e_2,\ldots,e_1,e_2,e_1,e_3,e_2,e_3)=0
\end{align*}
and
\begin{align*}
R^{p}\omega (e_1,e_2,e_1,e_2,\ldots,e_1,e_2,e_2,e_3,e_1,e_3)=0
\end{align*}
for all $p\geq 2$. Now (\ref{eq::OneBock2x2andOne1x1::eq8}) simplify to the form
\begin{align}\label{eq::OneBock2x2andOne1x1::eq9}
R^{p+1}\omega &(e_1,e_2,e_1,e_2,\ldots,e_1,e_2,e_1,e_3,e_1,e_3)\\
\nonumber&=-2\alpha h(e_1,e_2) R^{p}\omega (e_1,e_2,e_1,e_2,\ldots,e_1,e_2,e_1,e_3,e_1,e_3)
\end{align}
for $p\geq 2$. Note that using (\ref{eq::OneBock2x2andOne1x1::eq1})--(\ref{eq::OneBock2x2andOne1x1::eq3}) one may easily show that (\ref{eq::OneBock2x2andOne1x1::eq9}) is also true for $p=1$.
From (\ref{eq::OneBock2x2andOne1x1::eq1}) we see that (\ref{eq::OneBock2x2andOne1x1}) holds for $p=1$.
Assume now that (\ref{eq::OneBock2x2andOne1x1}) is true for some $p\geq 1$. From (\ref{eq::OneBock2x2andOne1x1::eq9}) we compute that
\begin{align*}
R^{p+1}\omega &(e_1,e_2,e_1,e_2,\ldots,e_1,e_2,e_1,e_3,e_1,e_3)\\
&=-2\alpha h(e_1,e_2)\cdot (-1)^p\cdot (2\eta\alpha)^{p-1}\varepsilon\omega(e_1,e_2)\\
&=(-1)^{p+1}\cdot (2\eta\alpha)^{p}\varepsilon\omega(e_1,e_2).
\end{align*}
Now by the induction principle (\ref{eq::OneBock2x2andOne1x1}) is true for all $p\geq 1$.
\end{proof}
\end{lem}


\begin{lem}\label{lm::OneBock2x2andOne1x1andAlpha0}
Let $\{e_1,e_2,\ldots,e_{2n}\}$ be vectors such that $Se_1=e_2$, $Se_2=0$, $h(e_1,e_1)=h(e_2,e_2)=0$, $h(e_1,e_2)=\eta$ and $Se_i=\lambda_{i-2} e_i$, $h(e_1,e_i)=h(e_2,e_i)=0$, $h(e_i,e_i)=\varepsilon_{i-2}$ for $i=3,\ldots,2n$ where $\eta\neq 0, \varepsilon\neq 0$.
Then for every $p\geq 1$ and for every $i=3,\ldots,2n$ we have
\begin{align}\label{eq::e1e2nminus1eie2n}
R^p\omega (e_{3},e_{1},e_{1},e_{2},\ldots,e_{1},e_{2},e_i,e_{2})=-\eta^p{\lambda_1}^p\omega(e_i,e_{3}).
\end{align}
\begin{proof}
\par First we shall show that for every $p\geq 1$ and for every $i=3,\ldots,2n$ the following formula holds:
\begin{align}\label{eq::e1e2neie2n}
R^p\omega (e_{3},e_{2},e_{1},e_{2},\ldots,e_{1},e_{2},e_i,e_{2})=0.
\end{align}
For $p=1$ we have
\begin{align*}
R\omega (e_{3},e_{2},e_i,e_{2})=-\omega (R(e_{3},e_{2})e_i,e_{2})-\omega (e_i,R(e_{3},e_{2})e_{2})=0,
\end{align*}
since
\begin{align}\label{eq::5::Re3e2ei}
R(e_{3},e_{2})e_i=0, \quad R(e_{3},e_{2})e_{2}=0.
\end{align}
Now assume that (\ref{eq::e1e2neie2n}) holds for some $p\geq 1$, we shall show that it also holds for $p+1$.
From the Gauss equation we have
\begin{align}\label{eq::5::Re3e2e1}
R(e_{3},e_{2})e_{1}=\eta Se_{3}=\eta \lambda_1e_{3}.
\end{align}
Now, using (\ref{eq::5::Re3e2ei}) and (\ref{eq::5::Re3e2e1}) we obtain
\begin{align*}
R^{p+1}\omega &(e_{3},e_{2},e_{1},e_{2},\ldots,e_{1},e_{2},e_i,e_{2})\\
&=-R^{p}\omega (\eta\lambda_1e_{3},e_{2},e_{1},e_{2},\ldots,e_{1},e_{2},e_i,e_{2})\\
&\phantom{=}-R^{p}\omega (e_{1},e_{2},\eta\lambda_1e_{3},e_{2},e_{1},e_{2},\ldots,e_{1},e_{2},e_i,e_{2})\\
&\phantom{=}-\ldots \\
&\phantom{=}-R^{p}\omega (e_{1},e_{2},e_{1},e_{2},\ldots,\eta \lambda_1e_{3},e_{2},e_i,e_{2})\\
&=-\eta\lambda_1R^{p}\omega (\overbrace{e_{3},e_{2}}^{1},e_{1},e_{2},\ldots,e_{1},e_{2},e_i,e_{2})\\
&\phantom{=}-\eta\lambda_1R^{p}\omega (e_{1},e_{2},\overbrace{e_{3},e_{2}}^{2},e_{1},e_{2},\ldots,e_{1},e_{2},e_i,e_{2})\\
&\phantom{=}-\ldots \\
&\phantom{=}-\eta\lambda_1R^{p}\omega (e_{1},e_{2},e_{1},e_{2},\ldots,\overbrace{e_{3},e_{2}}^{p},e_i,e_{2})\\
&=-\eta\lambda_1R^{p}\omega (e_{3},e_{2},e_{1},e_{2},\ldots,e_{1},e_{2},e_i,e_{2}),
\end{align*}
since all terms but first are equal $0$ due to the following identities
\begin{align*}
R(e_{1},e_{2})e_{3}&=R(e_{1},e_{2})e_i=R(e_{1},e_{2})e_{2}=0,\\
R(e_{1},e_{2})e_{1}&=\eta Se_{1}=\eta e_{2}.
\end{align*}
Now by the induction principle (\ref{eq::e1e2neie2n}) is true for all $p\geq 1$.
\par Now we can prove (\ref{eq::e1e2nminus1eie2n}).
First we check that (\ref{eq::e1e2nminus1eie2n}) is satisfied for $p=1$ and $p=2$. Indeed, since
\begin{align}
\label{eq::gauss::1}R(e_{3},e_{1})e_i=-h(e_{3},e_i)e_{2}, \quad R(e_{3},e_{1})e_{2}=\eta\lambda_1e_{3}, \quad R(e_{3},e_{1})e_{1}=0,\\
\label{eq::gauss::2}R(e_{1},e_{2})e_{3}=R(e_{1},e_{2})e_{2}=R(e_{1},e_{2})e_i=0,\quad R(e_{1},e_{2})e_{1}=\eta e_{2}
\end{align}
by straightforward computations we easily obtain that
\begin{align*}
R\omega &(e_{3},e_{1},e_i,e_{2})=-\omega (R(e_{3},e_{1})e_i,e_{2})-\omega (e_i,R(e_{3},e_{1})e_{2})\\
&=h(e_{3},e_i)\omega(e_{2},e_{2})-\eta\lambda_1\omega(e_i,e_{3})=-\eta\lambda_1\omega(e_i,e_{3})
\end{align*}
and
\begin{align*}
R^2\omega&(e_{3},e_{1},e_{1},e_{2},e_i,e_{2})\\
&=-R\omega(R(e_{3},e_{1})e_{1},e_{2},e_i,e_{2})-R\omega(e_{1},R(e_{3},e_{1})e_{2},e_i,e_{2})\\
&\phantom{=}-R\omega(e_{1},e_{2},R(e_{3},e_{1})e_i,e_{2})-R\omega(e_{1},e_{2},e_i,R(e_{3},e_{1})e_{2})\\
&=-R\omega(e_{1},\eta\lambda_1e_{3},e_i,e_{2})-R\omega(e_{1},e_{2},e_i,\eta\lambda_1e_{3})\\
&=\eta\lambda_1(-\eta\lambda_1\omega(e_i,e_{3}))-0=-\eta^2\lambda_{1}^2\omega(e_i,e_{3}).
\end{align*}
Let us assume that (\ref{eq::e1e2nminus1eie2n}) holds for some $p\geq 2$. Then
using (\ref{eq::gauss::1}) we compute
\begin{align*}
R^{p+1}\omega &(e_{3},e_{1},e_{1},e_{2},\ldots,e_{1},e_{2},e_i,e_{2})\\
&=-R^{p}\omega (e_{1},\eta\lambda_1e_{3},e_{1},e_{1},e_{2},\ldots,e_{1},e_{2},e_i,e_{2})\\
&\phantom{=}-R^{p}\omega (e_{1},e_{2},e_{1},\eta\lambda_1e_{3},e_{1},e_{2},\ldots,e_{1},e_{2},e_i,e_{2})\\
&\phantom{=}\ldots \\
&\phantom{=}-R^{p}\omega (e_{1},e_{2},e_{1},e_{2},\ldots,e_{1},\eta\lambda_1e_{3},e_i,e_{2})\\
&\phantom{=}-R^{p}\omega (e_{1},e_{2},e_{1},e_{2},\ldots,e_{1},e_{2},-h(e_{3},e_i)e_{2},e_{2})\\
&\phantom{=}-R^{p}\omega (e_{1},e_{2},e_{1},e_{2},\ldots,e_{1},e_{2},e_i,\eta\lambda_1e_{3})\\
&=-\eta\lambda_1R^{p}\omega (\overbrace{e_{1},e_{3}}^{1},e_{1},e_{2},\ldots,e_{1},e_{2},e_i,e_{2})\\
&\phantom{=}-\eta\lambda_1R^{p}\omega (e_{1},e_{2},\overbrace{e_{1},e_{3}}^{2},e_{1},e_{2},\ldots,e_{1},e_{2},e_i,e_{2})\\
&\phantom{=}\ldots \\
&\phantom{=}-\eta\lambda_1R^{p}\omega (e_{1},e_{2},\ldots,e_{1},e_{2},\overbrace{e_{1},e_{3}}^{p},e_i,e_{2})\\
&\phantom{=}-\eta\lambda_1R^{p}\omega (e_{1},e_{2},\ldots,e_{1},e_{2},e_i,e_{3}).
\end{align*}
From (\ref{eq::gauss::2}) it follows that for $j=2,\ldots,p$ and $i\in\{1,\ldots,2n\}\setminus\{2\}$
\begin{align}\label{eq::jotapoz}
R^{p}\omega& (e_{1},e_{2},\ldots,\overbrace{e_{1},e_{3}}^{j},\ldots,e_{1},e_{2},e_i,e_{2})\\
\nonumber &=-R^{p-1}\omega(e_{1},e_{2},\ldots,\overbrace{\eta e_{2},e_{3}}^{j-1},\ldots,e_{1},e_{2},e_i,e_{2}).
\end{align}
We also have
\begin{align*}
R^{p}\omega (e_{1},e_{2},\ldots,e_{1},e_{2},e_i,e_{3})=0
\end{align*}
for $i=3,2,\ldots,2n$ thanks to (\ref{eq::gauss::2}) and since $p\geq 2$.
Now we obtain
\begin{align*}
R^{p+1}\omega &(e_{3},e_{1},e_{1},e_{2},\ldots,e_{1},e_{2},e_i,e_{2})\\
&=\eta\lambda_1R^{p}\omega (e_{3},e_{1},e_{1},e_{2},\ldots,e_{1},e_{2},e_i,e_{2})\\
&\phantom{=}+\eta\lambda_1R^{p-1}\omega(\overbrace{\eta e_{2},e_{3}}^{1},e_{1},e_{2},\ldots,e_{1},e_{2},e_i,e_{2})\\
&\phantom{=}\ldots \\
&\phantom{=}+\eta\lambda_1R^{p-1}\omega(e_{1},e_{2},\ldots,e_{1},e_{2},\overbrace{\eta e_{2},e_{3}}^{p-1},e_i,e_{2})\\
&=\eta\lambda_1R^{p}\omega (e_{3},e_{1},e_{1},e_{2},\ldots,e_{1},e_{2},e_i,e_{2})\\
&\phantom{=}+\eta^2\lambda_1R^{p-1}\omega(e_{2},e_{3},e_{1},e_{2},\ldots,e_{1},e_{2},e_i,e_{2})\\
&=\eta\lambda_1R^{p}\omega (e_{3},e_{1},e_{1},e_{2},\ldots,e_{1},e_{2},e_i,e_{2}),
\end{align*}
where the last two equalities follow from (\ref{eq::gauss::2}) and (\ref{eq::e1e2neie2n}) respectively.
By the induction principle (\ref{eq::e1e2nminus1eie2n}) is true for all $p\geq 1$.
\end{proof}
\end{lem}


\begin{lem}\label{lm::OneBock2x2andOne1x1andAlpha0Next}
Let $\{e_1,e_2,\ldots,e_{2n}\}$ be vectors with properties like in the Lemma \ref{lm::OneBock2x2andOne1x1andAlpha0}.
Then for every $p\geq 1$ we have
\begin{align}\label{eq::TBA}
R^{p}\omega (e_{3},e_{2},e_{1},e_{2},\ldots,e_{1},e_{2})&=(-1)^p\eta^p\lambda_1^p\omega(e_{3},e_{2})\\
\label{eq::TBA2}
R^{p}\omega (e_{3},e_{1},e_{1},e_{2},\ldots,e_{1},e_{2})&=-\eta^p\lambda_1^p\omega(e_{1},e_{3})\\
\nonumber&\phantom{=}+\frac{1}{2}\Big((-1)^p+1\Big)\eta^p\lambda_1^{p-1}\omega(e_{2},e_{3}).
\end{align}
\begin{proof}
Since basis $\{e_{1},\ldots,e_{2n}\}$ satisfy conditions of Lemma \ref{lm::OneBock2x2andOne1x1andAlpha0}.
in particular we have (\ref{eq::gauss::2}) and (\ref{eq::jotapoz}).
We also have
\begin{align}\label{eq::12n2n-12n}
R(e_{3},e_{2})e_{1}=\eta\lambda_1e_{3}, \quad R(e_{3},e_{2})e_{2}=0.
\end{align}
By direct computations we check that
\begin{align*}
R\omega (e_{3},e_{2},e_{1},e_{2})&=-\eta\lambda_1\omega(e_{3},e_{2}),\\
R\omega (e_{3},e_{1},e_{1},e_{2})&=-\eta\lambda_1\omega(e_{1},e_{3}),\\
R^2\omega(e_{3},e_{1},e_{1},e_{2},e_{1},e_{2})&=-\eta^2{\lambda_1}^2\omega(e_{1},e_{3})+\eta^2\lambda_1\omega(e_{2},e_{3}).
\end{align*}
It means that (\ref{eq::TBA}) is true for $p=1$ and (\ref{eq::TBA2}) is true for $p=1,2$.
\par In order to show that (\ref{eq::TBA}) is true for all $p\geq 1$, using (\ref{eq::12n2n-12n}) we compute that
\begin{align*}
R^{p+1}\omega &(e_{3},e_{2},e_{1},e_{2},\ldots,e_{1},e_{2})\\
&=-R^p\omega(\eta\lambda_1e_{3},e_{2},e_{1},e_{2},\ldots,e_{1},e_{2})\\
&\phantom{=}-R^p\omega(e_{1},e_{2},\eta\lambda_1e_{3},e_{2},\ldots,e_{1},e_{2})\\
&\phantom{=}\ldots \\
&\phantom{=}-R^p\omega(e_{1},e_{2},\ldots,e_{1},e_{2},\eta\lambda_1e_{3},e_{2})\\
&=-\eta\lambda_1R^p\omega(e_{3},e_{2},e_{1},e_{2},\ldots,e_{1},e_{2}),
\end{align*}
since all terms but first are equal to $0$ due to (\ref{eq::gauss::2}). Now by the induction principle (\ref{eq::TBA}) is true for all $p\geq 1$.

\par To prove (\ref{eq::TBA2}) let us assume that (\ref{eq::TBA2}) is true for some $p\geq 2$ we compute
\begin{align*}
R^{p+1}\omega &(e_{3},e_{1},e_{1},e_{2},\ldots,e_{1},e_{2})\\
&=-R^p\omega(e_{1},\eta\lambda_1e_{3},e_{1},e_{2},\ldots,e_{1},e_{2})\\
&\phantom{=}-R^p\omega(e_{1},e_{2},e_{1},\eta\lambda_1e_{3},\ldots,e_{1},e_{2})\\
&\phantom{=}\ldots \\
&\phantom{=}-R^p\omega(e_{1},e_{2},\ldots,e_{1},e_{2},e_{1},\eta\lambda_1e_{3})\\
&=-\eta\lambda_1R^p\omega(e_{1},e_{3},e_{1},e_{2},\ldots,e_{1},e_{2})\\
&\phantom{=}-\eta\lambda_1R^p\omega(e_{1},e_{2},e_{1},e_{3},\ldots,e_{1},e_{2})\\
&\phantom{=}\ldots \\
&\phantom{=}-\eta\lambda_1R^p\omega(e_{1},e_{2},\ldots,e_{1},e_{2},e_{1},e_{3}).
\end{align*}
Now using (\ref{eq::jotapoz}) (for $i=1$) we obtain
\begin{align*}
R^{p+1}\omega &(e_{3},e_{1},e_{1},e_{2},\ldots,e_{1},e_{2})\\
&=-\eta\lambda_1R^p\omega(e_{1},e_{3},e_{1},e_{2},\ldots,e_{1},e_{2})\\
&\phantom{=}+\eta^2\lambda_1R^{p-1}\omega (e_{2},e_{3},e_{1},e_{2},\ldots,e_{1},e_{2})\\
&\phantom{=}+\eta^2\lambda_1R^{p-1}\omega (e_{1},e_{2},e_{2},e_{3},\ldots,e_{1},e_{2})\\
&\phantom{=}\ldots\\
&\phantom{=}+\eta^2\lambda_1R^{p-1}\omega (e_{1},e_{2},\ldots,e_{2},e_{3},e_{1},e_{2})\\
&\phantom{=}-\eta\lambda_1R^p\omega(e_{1},e_{2},\ldots,e_{1},e_{2},e_{1},e_{3})\\
&=-\eta\lambda_1R^p\omega(e_{1},e_{3},e_{1},e_{2},\ldots,e_{1},e_{2})\\
&\phantom{=}+\eta^2\lambda_1R^{p-1}\omega (e_{2},e_{3},e_{1},e_{2},\ldots,e_{1},e_{2})\\
&=\eta\lambda_1 \Big(-\eta^p\lambda_1^p\omega(e_{1},e_{3})
+\frac{1}{2}((-1)^p+1)\eta^p\lambda_1^{p-1}\omega(e_{2},e_{3})\Big)\\
&\phantom{=}-\eta^2\lambda_1(-1)^{p-1}\eta^{p-1}\lambda_1^{p-1}\omega(e_{3},e_{2})\\
&=-\eta^{p+1}\lambda_1^{p+1}\omega(e_{1},e_{3})+\frac{1}{2}((-1)^p+1+2\cdot(-1)^{p-1})\eta^{p+1}\lambda_1^{p}\omega(e_{2},e_{3})\\
&=-\eta^{p+1}\lambda_1^{p+1}\omega(e_{1},e_{3})+\frac{1}{2}((-1)^{p+1}+1)\eta^{p+1}\lambda_1^{p}\omega(e_{2},e_{3})
\end{align*}
where the last equalities are consequence of (\ref{eq::TBA2})  and (\ref{eq::TBA}) (for $p-1$). By the induction principle (\ref{eq::TBA2}) is true for all $p\geq 1$.
\end{proof}
\end{lem}

\begin{thm}\label{tw::FinalFormOfS}
Let $f\colon \M\rightarrow\R^{2n+1}$ ($\dim M\geq 4$)  be a non-degenerate affine hypersurface with a locally equiaffine transversal vector field  $\xi$
and an almost symplectic form $\omega$. If $R^p\omega=0$ for some $p\geq 1$ then for every point $x\in M$ either $S=0$ in $x$ or
there exists a basis ${e_1,\ldots,e_{2n}}$ of $T_xM$
such that $S$ in this basis has the form
\begin{align}\label{eq::SDecompositionFinal}
S=\left[\begin{matrix}
0 & 0 & 0 & \ldots & 0 \\
1 & 0 & 0 & \ldots & 0 \\
0 & 0 & 0 & \ldots & 0 \\
\vdots & \vdots & \ddots & \ldots \\
0 & 0 &0  & \ldots & 0
\end{matrix}\right]
\end{align}
\begin{proof}
From Theorem \ref{tw::NoComplex} we have that Jordan block decomposition of $S$ do not contain complex blocks.
From Theorem \ref{tw::MaxOneRealOfDim2} we also know that $S$ contains at most one real Jordan block of dimension $2$
and remaining blocks are all of dimension 1.
\par Now (rearranging vectors ${e_1,\ldots,e_{2n}}$ if needed) $S$ and $h$ can be represented in one of the following two forms:
\begin{align}\label{eq::SDecompositionDiag}
S=\left[\begin{matrix}
\lambda_1 & 0 & 0 & \ldots & 0 \\
0 & \lambda_2 & 0 & \ldots & 0 \\
0 & 0 & \lambda_3 & \ldots & 0 \\
\vdots & \vdots & \vdots  & \ddots & \vdots\\
0 & 0 &0  & \ldots & \lambda_{2n}
\end{matrix}\right],
\quad
h=\left[\begin{matrix}
\varepsilon_1 & 0 & 0 & \ldots & 0 \\
0 & \varepsilon_2 & 0 & \ldots & 0 \\
0 & 0 & \varepsilon_3 & \ldots & 0 \\
\vdots & \vdots & \vdots  & \ddots & \vdots\\
0 & 0 &0  & \ldots & \varepsilon_{2n}
\end{matrix}\right]
\end{align}
or
\begin{align}\label{eq::SDecompositionWithOneBlock}
S=\left[\begin{matrix}
\alpha & 0  & 0 & \ldots   &0  \\
1 & \alpha  & 0 & \ldots   &0 \\
0 & 0 & \lambda_1 & \ldots &0 \\
\vdots & \vdots& \vdots  & \ddots  &\vdots\\
0 & 0 & 0& \ldots & \lambda_{2n-2} \\
\end{matrix}\right],
\quad
h=\left[\begin{matrix}
0 & \eta  & 0 & \ldots   &0  \\
\eta & 0  & 0 & \ldots   &0 \\
0 & 0 & \varepsilon_1 & \ldots &0 \\
\vdots & \vdots& \vdots  & \ddots  &\vdots\\
0 & 0 & 0& \ldots & \varepsilon_{2n-2} \\
\end{matrix}\right]
\end{align}

where $\varepsilon_{i}\in\{-1,1\}$ for $i=1,\ldots,2n$, $\eta\in\{-1,1\}$ and
\begin{align}\label{eq::LambdaInequality}
|\lambda_1|\geq\cdots\geq |\lambda_{2n}|.
\end{align}
We need to show that $\lambda_1=\cdots=\lambda_{2n}=0$ (respectively $\alpha=0$ and $\lambda_1=\cdots=\lambda_{2n-2}=0$).
\par First assume that (\ref{eq::SDecompositionDiag}) holds. Since $\omega$ is non-degenerate there exists $i_0$ such that $\omega(e_1,e_{i_0})\neq 0$.
If $i_0>2$ then using Lemma \ref{lm::NablaKOmega0::2} ($s=2n$, $k=1$, $j=2$, $i=i_0$) we get
\begin{align*}
R^{2p}\omega (e_1,e_2,e_1,e_2,\ldots,e_1,e_2,e_1,e_{i_0})=(-1)^p\eps_1^p\eps_2^p \lambda_1^p\lambda_2^p\omega(e_1,e_{i_0}).
\end{align*}
If $i_0=2$ then from Lemma \ref{lm::NablaKOmega0::2} ($s=2n$, $k=1$, $j=3$, $i=i_0$) we get
\begin{align*}
R^{2p}\omega (e_1,e_3,e_1,e_3,\ldots,e_1,e_3,e_1,e_{i_0})=(-1)^p\eps_1^p\eps_3^p \lambda_1^p\lambda_3^p\omega(e_1,e_{i_0}).
\end{align*}
Since $R^p\omega=0$ then also $R^{2p}\omega=0$ and the above implies that
$$
(-1)^p\eps_1^p\eps_2^p \lambda_1^p\lambda_2^p\omega(e_1,e_{i_0})=(-1)^p\eps_1^p\eps_3^p \lambda_1^p\lambda_3^p\omega(e_1,e_{i_0})=0.
$$
Taking into account (\ref{eq::LambdaInequality}) we deduce that $\lambda_2=\lambda_3=0$ if $i_0>2$ and $\lambda_3=\lambda_4=0$ if $i_0=2$.
Now, if $i_0\neq 4$ using Lemma \ref{lm::R2kOmegaXZ1Z2Y} ($X=e_1$, $Y=e_{i_0}$, $Z_1=e_3$, $Z_4=e_4$) we get
\begin{align*}
R^{2p}\omega (e_1,\underbrace{e_3,e_4,e_3,e_4,\ldots,e_3,e_4}_{4p},e_{i_0})=(-1)^p\lambda_1^{2p}\varepsilon_3^p \varepsilon_4^p\omega(e_1,e_{i_0}).
\end{align*}
If $i_0=4$ then we have that also $\lambda_2=0$ and in this case from Lemma \ref{lm::R2kOmegaXZ1Z2Y} ($X=e_1$, $Y=e_{i_0}$, $Z_1=e_2$, $Z_2=e_3$) we get
\begin{align*}
R^{2p}\omega (e_1,\underbrace{e_2,e_3,e_2,e_3,\ldots,e_2,e_3}_{4p},e_{i_0})=(-1)^p\lambda_1^{2p}\varepsilon_2^p \varepsilon_3^p\omega(e_1,e_{i_0}).
\end{align*}
The above implies that $\lambda_1=0$ and thanks to (\ref{eq::LambdaInequality}) we get that $S=0$.
\par Now assume that $S$ and $h$ have the form (\ref{eq::SDecompositionWithOneBlock}).
First we shall show that $\alpha=0$. For this purpose note that from Corollary  \ref{cor::rpomegaeiek} ($k=2$) we have
\begin{align*}
R^p\omega (e_1,e_2,\ldots,e_1,e_2,e_i,e_2)=\eta^p\alpha ^p\omega (e_i,e_2)
\end{align*}
for $i\in\{2,\ldots,2n\}$. From Lemma \ref{lm::OneBock2x2andOne1x1} ($\varepsilon=\varepsilon_1$) we have
\begin{align*}
R^{p}\omega (e_1,e_2,e_1,e_2,\ldots,e_1,e_2,e_1,e_3,e_1,e_3)=(-1)^p\cdot (2\eta\alpha)^{p-1}\varepsilon_1\omega(e_1,e_2).
\end{align*}
Since $R^p\omega=0$ and $\eta,\varepsilon_1\neq 0$ we obtain
$$
\alpha ^p\omega (e_2,e_2)=\cdots=\alpha ^p\omega (e_{2n},e_2)=\alpha^{p-1}\omega(e_1,e_2)=0
$$
and in consequence $\alpha=0$ (since $\omega$ is non-degenerate).

\par Now we are able to show that $\lambda_1=\cdots=\lambda_{2n-2}=0$. Indeed, from (\ref{eq::LambdaInequality}) it follows that it is enough to show that $\lambda_1=0$.
Since $\alpha=0$, the basis $\{e_{1},\ldots,e_{2n}\}$ satisfy conditions of Lemma \ref{lm::OneBock2x2andOne1x1andAlpha0} and Lemma \ref{lm::OneBock2x2andOne1x1andAlpha0Next}.
Thus using formulas (\ref{eq::e1e2nminus1eie2n}), (\ref{eq::TBA}) and (\ref{eq::TBA2}) and taking into account that $R^p\omega=0$ we obtain
\begin{align}
\label{eq::lambda1::1}-\eta^p{\lambda_1}^p\omega(e_i,e_{3})&=0 \qquad \text{for $i=3,\ldots,2n$},\\
\label{eq::lambda1::2}(-1)^p\eta^p\lambda_1^p\omega(e_{3},e_{2})&=0, \\
\label{eq::lambda1::3}-\eta^p\lambda_1^p\omega(e_{1},e_{3})+\frac{1}{2}\Big((-1)^p+1\Big)&\eta^p\lambda_1^{p-1}\omega(e_{2},e_{3})=0.
\end{align}
Since $\omega$ is non-degenerate there must exist $i_0\in\{1,\ldots,2n\}\setminus\{3\}$ such that $\omega(e_{i_0},e_{3})\neq 0$. If $i_0>3$ then from (\ref{eq::lambda1::1})
we immediately get that $\lambda_1=0$. If $i_0=2$ then (\ref{eq::lambda1::2}) implies that $\lambda_1=0$. Eventually, if $i_0=1$ and $\omega(e_{2},e_{3})=0$ we obtain that $\lambda_1=0$ from (\ref{eq::lambda1::3}).
The proof of the theorem is completed.
\end{proof}
\end{thm}

As a consequence of Theorem \ref{tw::FinalFormOfS} we obtain

\begin{thm}\label{tw::RKomega0General}
Let $f\colon \M\rightarrow\R^{2n+1}$ ($\dim M\geq 4$)  be a non-degenerate affine hypersurface with a locally equiaffine transversal vector field  $\xi$
and an almost symplectic form $\omega$. If $R^k\omega=0$ for some $k\geq 1$ then the shape operator $S$ has the rank $\leq 1$.
\begin{proof}
From Theorem \ref{tw::FinalFormOfS} it follows that for every $x\in M$ either $S_x=0$ or $S_x$ is of the form (\ref{eq::SDecompositionFinal}) thus  $\rank S_x\leq 1$ for every $x\in M$.
\end{proof}
\end{thm}

\par Recall that we have the following lemma (\cite{MSz}).

\begin{lem}[\cite{MSz}]\label{lm::WzorNaRkT}
Let $T$ be a tensor of  type $(0,p)$ and let $\nabla$ be an affine torsion-free connection. Then for every $k\geq 1$ and for any $2k+p$ vector fields
$X_{\pm 1}^1,\ldots,X_{\pm 1}^k$, $Y_1,\ldots,Y_p$ the following identity holds:
\begin{align}\label{eq::WzorNaRkT}
(R^k&\cdot T)(X_1^1,X_{-1}^1,\ldots,X_1^k,X_{-1}^k,Y_1,\ldots,Y_p)\\
\nonumber &=\sum_{a\in\mathcal{J}}\sgn a(\nabla^{2k}T)(X_{a(1)}^1,X_{-a(1)}^1,\ldots,X_{a(k)}^k,X_{-a(k)}^k,Y_1,\ldots,Y_p),
\end{align}
where $\mathcal{J}=\{a\colon I_k\rightarrow \{-1,1\}\}$ and $\sgn a:=a(1)\cdot\ldots\cdot a(k)$.
\end{lem}

From Theorem \ref{tw::RKomega0General} and Lemma \ref{lm::WzorNaRkT} we have the following
\begin{thm}\label{tw::NablaKOmega0General}
Let $f\colon \M\rightarrow\R^{2n+1}$ ($\dim M\geq 4$) be a non-degenerate affine hypersurface with a locally equiaffine transversal vector field  $\xi$
and an almost symplectic form $\omega$. If $\nabla^p\omega=0$ for some $p\geq 1$ then the shape operator $S$ has the rank $\leq 1$.
\end{thm}

We conclude this section with the following example
\begin{ex}
Let $n\geq 2$ and let $\gamma,\alpha_i\colon\R\rightarrow\R^{2n+1}$ be curves given as follows:
\begin{align*}
  \gamma(t)&:=(\cos t,\sin t,0,\ldots,0)^T; \\
  \alpha_i(t)&:=(0,\ldots,0,\overbrace{\cos t}^{(i+3)},\overbrace{\sin t}^{(i+4)},0,\ldots,0)^T
\end{align*}
for $i=0,\ldots,2n-3$. Let $\varepsilon_i\in\{-1,1\}$ for $i=1,\ldots,2n-3$.
\par Now let us consider an immersion  $f\colon \R\setminus\{0\}\times(0,\infty)\times(0,\frac{\pi}{2})^{2n-2}\rightarrow\R^{2n+1}$ given by the formula:
$$
f(x,y,z_0,\ldots,z_{2n-3}):=y\gamma'(x)+x\alpha_0(z_0)+\sum_{i=1}^{2n-3}\varepsilon_i\alpha_i(z_i)
$$
together with the transversal vector field
$$
\xi\colon\R\setminus\{0\}\times(0,\infty)\times(0,\frac{\pi}{2})^{2n-2}\ni (x,y,z_0,\ldots,z_{2n-3})\mapsto -\gamma(x)\in\R^{2n+1}.
$$

By straightforward computations we get
$$
h=\left[\begin{matrix}
    0      & 1      & 0      & 0                             & \cdots & 0  \\
    1      & 0      & 0      &  0                            & \cdots & 0  \\
    0      & 0      & xy     & 0                             & \cdots & 0  \\
    0      & 0      & 0      & \varepsilon_1 y\frac{\sin(z_0)}{\cos(z_1)} & \cdots & 0  \\
    \vdots & \vdots & \vdots & \vdots                        & \ddots & \vdots \\
    0      & 0      & 0      & 0                             &\cdots  & \varepsilon_{2n-3}y\frac{\sin(z_0)}{\cos(z_{2n-3})}\prod_{i=1}^{2n-4}{\tan{z_i}} \\ 
  \end{matrix}\right]
$$
and
$$
S=\left[\begin{matrix}
    0 & 0 & 0  & \cdots & 0 \\
    1 & 0 & 0 & \cdots & 0\\
    0 & 0 & 0 & \cdots & 0\\
    \vdots & \vdots & \vdots &  \ddots &\vdots \\
    0 & 0 & 0 & \cdots & 0
  \end{matrix}\right],\quad \tau=0.
$$

Thus $f$ is a non-degenerate equiaffine hypersurface with the second fundamental form of signature $(1,-1,\varepsilon_0,\varepsilon_{1},\ldots,\varepsilon_{2n-3})$ where $\varepsilon_0=1$ for $x>0$ and $\varepsilon_0=-1$ for $x<0$.
Note also that $R\neq 0$ since $S\neq 0$.

Now let us define
$$
\Omega=[\omega_{i,j}]_{i,j=1\ldots 2n}
$$
such that $\omega_{i,j}=-\omega_{j,i}$ and $\det \Omega\neq 0$. That is $\Omega$ is a symplectic form.
We easily check that
$$
R\Omega\Big(\pszm{x},\pszm{z_0},\pszm{x},\pszm{z_0}\Big)=-xy\omega_{1,2}
$$
and
$$
R^2\Omega\Big(\pszm{x},\pszm{z_0},\pszm{x},\pszm{z_0},\pszm{x},\pszm{z_0}\Big)=xy\omega_{2,3}
$$
thus $R\Omega\neq 0$ and $R^2\Omega\neq 0$ if only $\omega_{1,2}\neq0, \omega_{2,3}\neq 0$. On the other hand one may show that $R^3\Omega = 0$.
\end{ex}

\emph{This Research was financed by the Ministry of Science and Higher Education of the Republic of Poland.}

\end{document}